\documentclass[a4paper,11pt]{article}
\usepackage[utf8]{inputenc}

\usepackage{graphicx}
\usepackage{amsthm}
\newtheorem{remark}{\bf Remark}[section]

\newtheorem{example}{Example}[section]

\usepackage[nottoc]{tocbibind}
\usepackage{color}
\usepackage{slantsc}
\usepackage{framed,multirow}

\usepackage{amssymb}
\usepackage{latexsym}

\usepackage{amsmath}
\usepackage{caption}

\newcommand{\FD}{\text{FD}}
\newcommand{\LSZ}{\text{\tiny LSZ}}
\newcommand{\MWENO}{\text{\tiny MWENO}}
\newcommand{\CWENO}{\text{\tiny CWENO}}
\newcommand{\OPT}{\text{\fontsize{5pt}{6pt}\selectfont OPT}}
\newcommand{\Le}{\text{\tiny L}}
\newcommand{\Mi}{\text{\tiny M}}
\newcommand{\Ri}{\text{\tiny R}}
\newcommand{\C}{\text{\tiny C}}
\newcommand{\cfl}{\text{CFL}}

\usepackage{url}
\usepackage{xcolor}
\definecolor{newcolor}{rgb}{.8,.349,.1}

\begin{document}
\title{A sixth-order central WENO scheme for nonlinear degenerate parabolic equations}
\date{}
\author{
{ Samala Rathan\thanks{Email: rathans.math@iipe.ac.in, Faculty of Mathematics, Department of Humanities and Sciences, Indian Institute of Petroleum and Energy-Visakhapatnam, 530003, India},}
{ Jiaxi Gu\thanks{Email: jiaxigu@postech.ac.kr, Department of Mathematics \& POSTECH MINDS (Mathematical Institute for Data Science), Pohang University of Science and Technology, Pohang 37673, Korea}}}
\maketitle
\begin{abstract}
In this paper we develop a new sixth-order finite difference central weighted essentially non-oscillatory (WENO) scheme with Z-type nonlinear weights for nonlinear degenerate parabolic equations.
The centered polynomial is introduced for the WENO reconstruction in order to avoid the negative linear weights.
We choose the Z-type nonlinear weights based on the $L^2$-norm smoothness indicators, yielding the new WENO scheme with more accurate resolution.
It is also confirmed that the proposed central WENO scheme with the devised nonlinear weights achieves sixth order accuracy in smooth regions.
One- and two-dimensional numerical examples are presented to demonstrate the improved performance of the proposed central WENO scheme.
\end{abstract}
{\small \textbf{Keywords:} Finite difference, Central WENO scheme, Z-type nonlinear weights, Nonlinear degenerate parabolic equation.}

\section{Introduction} \label{sec:intro}
In this paper we are interested in solving the one-dimensional parabolic equation
\begin{equation} \label{eq:parabolic}
 u_t = b(u)_{xx},
\end{equation}
where $u = u(x,t)$ is a scalar quantity and $b'(u) \geqslant 0$. 

One example of such nonlinear equation is the porous medium equation (PME) of the degenerate parabolic type,
\begin{equation} \label{eq:pme}
 u_t = (u^m)_{xx},~~m>1,
\end{equation}
which describes the flow of an isentropic gas through a porous medium \cite{Aronson, Muskat}, the heat radiation in plasmas \cite{ZeldovichRaizer}, and various physical processes.
The classical linear heat equation can be considered to be the limit of PME \eqref{eq:pme} as $m \to 1$.
Assuming $u \geqslant 0$, the PME could be written in the form
$$
   u_t = \left( m u^{m-1} u_x \right)_x.
$$
Then the PME is parabolic only at those points where $u \ne 0$, while it degenerates as the vanishing of the term $m u^{m-1}$ wherever $u = 0$.
In other words, the PME is a degenerate parabolic equation.
One important property of the PME is the finite propagation, which is different from the infinite speed of propagation in the classical heat equation.
This property implies the appearance of free boundaries that separate the regions where the solution is positive from those where $u = 0$, giving rise to the sharp interfaces \cite{Vazquez}.
Since the free boundaries move with respect to time, their behavior looks similar to the behavior of shocks in the hyperbolic conservation laws.
So it is reasonable to apply the weighted essentially non-oscillatory (WENO) philosophy to the PME, enabling the free boundaries to be well resolved. 

Before we discuss the WENO schemes, we would like to mention several different schemes that specialize in nonlinear degenerate parabolic equations in the literature.
For example, the explicit diffusive kinetic schemes have been designed in \cite{Aregba}.
Also, the high-order relaxation scheme has been introduced in \cite{Cavali}.
A local discontinuous Galerkin finite element method for the PME was studied in \cite{ZhangWu}.
Other approaches based on the finite volume method were investigated in \cite{Bess, Arbogast}.
In the more general nonlinear degenerate convection-diffusion case, the entropy stable finite difference schemes were proposed in \cite{Jerez}. 

Our focus of this paper is on the WENO schemes. 
In \cite{Liu}, Liu et al. constructed the finite difference WENO (WENO-LSZ) schemes for the equation \eqref{eq:parabolic}, which approximate the second derivative term directly by a conservative flux difference.
However, unlike the positive linear weights of WENO schemes for hyperbolic conservation laws \cite{JiangShu, Shu}, the negative linear weights exist so that some special care, such as the technique in \cite{Shi}, was applied to guarantee the non-oscillatory performance in regions of sharp interfaces.
Following the definition of the smoothness indicators in \cite{JiangShu, Shu} and invoking the mapped function in \cite{Henrick}, the resulting nonlinear weights meet the requirement of sixth order accuracy.
In \cite{Hajipour}, Hajipour and Malek proposed the modified WENO (MWENO) scheme with Z-type nonlinear weights \cite{Borges} and nonstandard Runge–Kutta (NRK) schemes.
Further, the hybrid scheme, based on the spatial MWENO and the temporal NRK schemes, was employed to solve the equation \eqref{eq:parabolic} numerically.
Recently, Abedian et al. \cite{AbedianAdibiDehghan, Abedian} aimed at avoiding negative linear weights and presented some modifications to the numerical flux. 
In \cite{Rathan}, Rathan et al. showed a new type of local and global smoothness indicators in $L^1$ norm via undivided differences and subsequently constructed the new Z-type nonlinear weights. 
Christelieb et al. employed a kernel based approach with the philosophy of the method of lines transpose, giving a high-order WENO method with a nonlinear filter in \cite{Christlieb}. 
In \cite{Jiang}, Jiang designed an alternative formulation to approximate the second derivatives in a conservative form, where the odd order derivatives at half points were used to construct the numerical flux.

In this paper we present the central WENO (CWENO) scheme based on the point values for the diffusion term, following the notion of compact CWENO schemes based on the cell avarages for the convection term proposed by Levy et al. in \cite{Levy}.
The negative linear weights in \cite{Liu}, which require some special care, can be circumvented in our scheme.
We further devise the Z-type nonlinear weights \cite{Borges}, which are dependent on the smoothness indicators in \cite{Liu} of $L^2$ norms.
The global smoothness indicator can be designed to attain higher order so that we do not need the power to maintain the order of accuracy.
In our scheme, not only is the computational cost reduced without estimating the mapped function in WENO-LSZ or applying the splitting technique to treat the negative weights in both WENO-LSZ and MWENO, but the non-oscillatory performance is improved since there exist small-scale oscillations around the sharp interfaces for WENO-LSZ in some cases as the time advances whereas those oscillations are largely damped by our scheme.   
We also provide the sufficient conditions for sixth order accuracy in smooth regions and an analysis of nonlinear weights shows that the proposed WENO scheme is in compliance with those criteria.
The implementation of WENO schemes for the equation \eqref{eq:parabolic} could be extended to the convection–diffusion equations with the WENO schemes for convection terms \cite{Borges, Henrick, JiangShu, Shu, Gu} and to multi-space dimensions in a dimension-by-dimension approach.

The paper is organized as follows. 
In Section \ref{sec:weno}, the sixth-order WENO scheme for the parabolic equation \eqref{eq:parabolic} and some of the relevant analytical results are reviewed. 
Section \ref{sec:cweno} presents the CWENO approximation for the diffusion term with the new Z-type nonlinear weights based on $L^2$ norms, and the sufficient conditions for sixth order accuracy in smooth regions, which the devised nonlinear weights satisfy.
In Section \ref{sec:nr}, the proposed CWENO scheme and the WENO-LSZ and MWENO schemes are compared with the simulation of one- and two-dimensional numerical experiments, including 1D and 2D heat equations for the sixth-order verification; 1D and 2D PMEs with various initial conditions; 1D and 2D Buckley–Leverett equations; 1D and 2D strongly degenerate parabolic convection-diffusion equations.
A brief concluding remark is presented in Section \ref{sec:conclusion}.

\section{WENO approximation to the second derivative} \label{sec:weno}
In this section, we review the direct WENO discretization to the second derivative in the conservation form \cite{Liu}.
Consider a uniform grid defined by the points $x_0 < x_1 < \cdots < x_{N-1} < x_N$ with $x_{i \pm 1/2} = x_i \pm \Delta x/2,~i = 0, \dots N$, where $\Delta x = x_{i+1} - x_i$ is the uniform grid spacing.
Then the spatial domain is discretized by this uniform grid.
The semi-discrete form of Equation \eqref{eq:parabolic} with respect to $t$, yields
\begin{equation} \label{eq:parabolic_discrete}
 \frac{du_{i}(t)}{dt} = \left. \frac{\partial^2 b}{\partial x^2} \right|_{x=x_i},
\end{equation} 
where $u_i(t)$ is the numerical approximation to the point value $u(x_i,t)$.
Define the function $h(x)$ implicitly by
\begin{equation} \label{eq:primitive}
 b(u(x)) = \frac{1}{\Delta x^2} \int^{x+\Delta x/2}_{x-\Delta x/2} \left( \, \int^{\eta+\Delta x/2}_{\eta-\Delta x/2} h(\xi) d\xi \right) d\eta. 
\end{equation} 
Differentiating both sides twice with respect to $x$, we obtain
$$
   b(u)_{xx} = \frac{h(x+\Delta x) - 2 h(x) + h(x-\Delta x)}{\Delta x^2}.
$$
Setting $g(x) = h(x+\Delta x/2) - h(x-\Delta x/2)$ gives the equation
$$
   \left. \frac{\partial^2 b}{\partial x^2} \right|_{x=x_i} = \frac{g_{i+1/2} - g_{i-1/2}}{\Delta x^2}, 
$$
where $g_{i \pm 1/2} = g(x_{i \pm 1/2})$.
Then the equation \eqref{eq:parabolic_discrete} becomes 
\begin{equation} \label{eq:parabolic_discrete_h}
 \frac{du_{i}(t)}{dt} = \frac{g_{i+1/2} - g_{i-1/2}}{\Delta x^2}.
\end{equation} 

\begin{figure}[h!]
\centering
\includegraphics[width=0.8\textwidth]{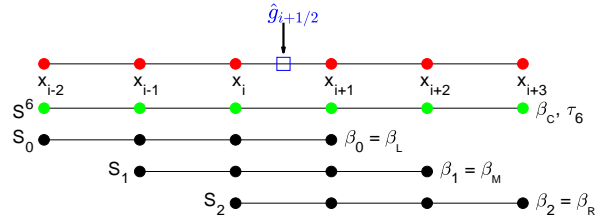}
\vspace{-4.5cm}
\captionof{figure}{The numerical flux $\hat{g}_{i+1/2}$ is constructed on the stencil $S^6 = \{ x_{i-2}, \cdots, x_{i+3} \}$ with six uniform points, as well as three 4-point substencils $S_0, S_1, S_2$.}
\label{fig:stencil}
\end{figure}

In order to approximate $g_{i+1/2}$, a polynomial approximation $q(x)$ to $h(x)$ of degree at most 5,
$$
   q(x) = a_0 + a_1 x + a_2 x^2 + a_3 x^3 + a_4 x^4 + a_5 x^5,
$$
can be constructed on the 6-point stencil $S^6$, as shown in Figure \ref{fig:stencil}. 
The polynomial $q(x)$ interpolates $b_{i+j} = b(u(x_{i+j},t)),~j = -2,\cdots,3$, in the sense of \eqref{eq:primitive}, which gives
\begin{align*}
 a_0 &= \frac{2 b_{i-2} - 23 b_{i-1} + 222 b_i - 23 b_{i+1} + 2 b_{i+2}}{180}, \\
 a_1 &= \frac{8 b_{i-2} - 55 b_{i-1} - 70 b_i + 160 b_{i+1} - 50 b_{i+2} + 7 b_{i+3}}{120 \Delta x}, \\
 a_2 &= - \frac{b_{i-2} - 10 b_{i-1} + 18 b_i - 10 b_{i+1} + b_{i+2}}{120 \Delta x^2}, \\
 a_3 &= - \frac{b_{i-2} + 4 b_{i-1} - 20 b_i + 26 b_{i+1} - 13 b_{i+2} + 2 b_{i+3}}{36 \Delta x^3}, \\
 a_4 &= \frac{b_{i-2} - 4 b_{i-1} + 6 b_i - 4 b_{i+1} + b_{i+2}}{24 \Delta x^4}, \\
 a_5 &= - \frac{b_{i-2} - 5 b_{i-1} + 10 b_i - 10 b_{i+1} + 5 b_{i+2} - b_{i+3}}{120 \Delta x^5}.
\end{align*}
The polynomial $p(x)$ of degree at most 4 approximating $g(x)$ is obtained by taking the difference of $q(x+\Delta x/2)$ and $q(x-\Delta x/2)$,
$$
   g(x) = h(x+\Delta x/2) - h(x-\Delta x/2) \approx q(x+\Delta x/2) - q(x-\Delta x/2) = p(x).
$$
Then we have
\begin{equation} \label{eq:polynomial}
 p(x) = c_0 + c_1 x + c_2 x^2 + c_3 x^3 + c_4 x^4,
\end{equation}
where
\begin{align*}
 c_0 &= \frac{341 b_{i-2} - 2785 b_{i-1} - 2590 b_i + 6670 b_{i+1} - 1895 b_{i+2} + 259 b_{i+3}}{5760}, \\
 c_1 &= - \frac{b_{i-2} - 12 b_{i-1} + 22 b_i - 12 b_{i+1} + b_{i+2}}{8 \Delta x}, \\
 c_2 &= - \frac{5 b_{i-2} + 11 b_{i-1} - 70 b_i + 94 b_{i+1} - 47 b_{i+2} + 7 b_{i+3}}{48 \Delta x^2}, \\
 c_3 &= \frac{b_{i-2} - 4 b_{i-1} + 6 b_i - 4 b_{i+1} + b_{i+2}}{6 \Delta x^3}, \\
 c_4 &= - \frac{b_{i-2} - 5 b_{i-1} + 10 b_i - 10 b_{i+1} + 5 b_{i+2} - b_{i+3}}{24 \Delta x^4}.
\end{align*}
Evaluating $p(x)$ at $x = x_{i+1/2}$ yields the finite difference numerical flux
\begin{equation} \label{eq:numerical_flux_plus}
 \hat{g}_{i+1/2}^{\FD} = p(x_{i+1/2}) = - \frac{1}{90} b_{i-2} + \frac{5}{36} b_{i-1} - \frac{49}{36} b_i + \frac{49}{36} b_{i+1} - \frac{5}{36} b_{i+2} + \frac{1}{90} b_{i+3}.
\end{equation}
The numerical flux $\hat{g}_{i-1/2}^{\FD}$ is obtained directly by shifting one grid to the left,
\begin{equation} \label{eq:numerical_flux_minus}
 \hat{g}_{i-1/2}^{\FD} = - \frac{1}{90} b_{i-3} + \frac{5}{36} b_{i-2} - \frac{49}{36} b_{i-1} + \frac{49}{36} b_i - \frac{5}{36} b_{i+1} + \frac{1}{90} b_{i+2}.
\end{equation}
Applying the Taylor expansions to $\hat{g}_{i \pm 1/2}^{\FD}$ \eqref{eq:numerical_flux_plus} and \eqref{eq:numerical_flux_minus} would give
\begin{align}
 \hat{g}_{i+1/2}^{\FD} &= g_{i+1/2} + \frac{1}{560} h_i^{(7)} \Delta x^7 + O(\Delta x^8), \label{eq:numerical_flux_plus_Taylor} \\
 \hat{g}_{i-1/2}^{\FD} &= g_{i-1/2} + \frac{1}{560} h_i^{(7)} \Delta x^7 + O(\Delta x^8). \label{eq:numerical_flux_minus_Taylor}
\end{align}
Replacing $g_{i \pm 1/2}$ in \eqref{eq:parabolic_discrete_h} by \eqref{eq:numerical_flux_plus_Taylor} and \eqref{eq:numerical_flux_minus_Taylor}, respectively, we have the sixth-order approximation
\begin{equation} \label{eq:FD_scheme}
 \frac{du_{i}(t)}{dt} = \frac{\hat{g}_{i+1/2}^{\FD} - \hat{g}_{i-1/2}^{\FD}}{\Delta x^2} + O(\Delta x^6).
\end{equation}

A similar approach can be used to obtain a polynomial $p_k(x)$ of degree at most 2 on each 4-point substentil $S_k = \{ x_{i-2+k}, \cdots, x_{i+1+k} \}$ with $k = 0,1,2$, 
where
\begin{equation} \label{eq:subpolynomials}
\begin{aligned}
 p_0(x) &= \frac{5 b_{i-2} - 27 b_{i-1} + 15 b_i + 7 b_{i+1}}{24} + \frac{b_{i-1} - 2 b_i + b_{i+1}}{\Delta x} x - \frac{b_{i-2} - 3 b_{i-1} + 3 b_i - b_{i+1}}{2 \Delta x^2} x^2, \\
 p_1(x) &= - \frac{7 b_{i-1} + 15 b_i - 27 b_{i+1} + 5 b_{i+2}}{24} + \frac{b_{i-1} - 2 b_i + b_{i+1}}{\Delta x} x - \frac{b_{i-1} - 3 b_i + 3 b_{i+1} - b_{i+2}}{2 \Delta x^2} x^2, \\
 p_2(x) &= - \frac{43 b_i - 69 b_{i+1} + 33 b_{i+2} - 7 b_{i+3}}{24} + \frac{2 b_i - 5 b_{i+1} + 4 b_{i+2} - b_{i+3}}{\Delta x} x - \frac{b_i - 3 b_{i+1}+ 3 b_{i+2} - b_{i+3}}{2 \Delta x^2} x^2.
\end{aligned}
\end{equation}
This results in the numerical fluxes $\hat{g}^k_{i+1/2}$,
\begin{equation} \label{eq:numerical_flux_substencil}
\begin{aligned}
 \hat{g}^0_{i+1/2} &=   \frac{1}{12} b_{i-2} - \frac{1}{4} b_{i-1} - \frac{3}{4} b_i     + \frac{11}{12} b_{i+1}, \\
 \hat{g}^1_{i+1/2} &=   \frac{1}{12} b_{i-1} - \frac{5}{4} b_i     + \frac{5}{4} b_{i+1} - \frac{1}{12} b_{i+2}, \\
 \hat{g}^2_{i+1/2} &= - \frac{11}{12} b_i    + \frac{3}{4} b_{i+1} + \frac{1}{4} b_{i+2} - \frac{1}{12} b_{i+3}.
\end{aligned}
\end{equation}
We could obtain the numerical fluxes $\hat{g}^k_{i-1/2}$ after shifting each index by -1.
Hence the Taylor series expansion gives
\begin{equation} \label{eq:numerical_flux_substencil_Taylor}
\begin{aligned}
 \hat{g}^0_{i \pm 1/2} &= g_{i \pm 1/2} + \frac{1}{12} h_i^{(4)} \Delta x^4 + O(\Delta x^5), \\
 \hat{g}^1_{i \pm 1/2} &= g_{i \pm 1/2} - \frac{1}{90} h_i^{(5)} \Delta x^5 + O(\Delta x^6), \\
 \hat{g}^2_{i \pm 1/2} &= g_{i \pm 1/2} - \frac{1}{12} h_i^{(4)} \Delta x^4 + O(\Delta x^5).
\end{aligned}
\end{equation}
It is clear that the linear combination of all $\hat{g}^k_{i+1/2}$ can produce $\hat{g}_{i+1/2}^{\FD}$ that approximates the flux $g_{i+1/2}$ in \eqref{eq:parabolic_discrete_h}, that is, there are linear weights $d_0 = d_2 = - \frac{2}{15}$ and $d_1 = \frac{19}{15}$ such that
\begin{equation} \label{eq:numerical_flux_linear_weights}
 \hat{g}_{i+1/2}^{\FD} = \sum_{k=0}^2 d_k \hat{g}^k_{i+1/2}.
\end{equation}
Similarly, an index shift by -1 returns the corresponding relation between $\hat{g}_{i-1/2}^{\FD}$ and $\hat{g}^k_{i-1/2}$.

Since \eqref{eq:numerical_flux_linear_weights} is not a convex combination of \eqref{eq:numerical_flux_substencil} as the linear weights $d_0$ and $d_2$ are negative, the WENO procedure cannot be applied directly to obtain a stable scheme.
The test cases in \cite{Shi} showed that WENO schemes without special treatment to the negative weights may lead to the blow-up of the numerical solution.
Thus the splitting technique in \cite{Shi} could be utilized to treat the negative weights $d_0$ and $d_2$. 
The linear weights are split into positive and negative parts,
$$
   \tilde{\gamma}^+_k = \frac{1}{2} \left( d_k + 3 | d_k | \right),~~\tilde{\gamma}^-_k = \tilde{\gamma}^+_k - d_k,~~k=0,1,2.
$$
Then $d_k = \tilde{\gamma}^+_k - \tilde{\gamma}^-_k$ and 
\begin{equation*}
\begin{aligned}
 \tilde{\gamma}^+_0 &= \frac{2}{15}, & \tilde{\gamma}^+_1 &= \frac{38}{15}, & \tilde{\gamma}^+_2 &= \frac{2}{15}; \\
 \tilde{\gamma}^-_0 &= \frac{4}{15}, & \tilde{\gamma}^-_1 &= \frac{19}{15}, & \tilde{\gamma}^+_2 &= \frac{4}{15}.
\end{aligned}
\end{equation*}
We scale them by
$$
   \sigma^{\pm} = \sum_{k=0}^2 \tilde{\gamma}_k^{\pm},~~\gamma_k^{\pm} = \tilde{\gamma}_k^{\pm} / \sigma^{\pm},~~k=0,1,2.
$$
Then the linear positive and negative weights $\gamma_k^{\pm}$ are given by
\begin{equation} \label{eq:gamma}
\begin{aligned}
 \gamma^+_0 &= \frac{1}{21}, & \gamma^+_1 &= \frac{19}{21}, & \gamma^+_2 &= \frac{1}{21}; \\
 \gamma^-_0 &= \frac{4}{27}, & \gamma^-_1 &= \frac{19}{27}, & \gamma^+_2 &= \frac{4}{27},
\end{aligned}
\end{equation}
which satisfy 
\begin{equation} \label{eq:linear_weight_relation}
 d_k = \sigma^+ \gamma_k^+ - \sigma^- \gamma_k^-.
\end{equation}

Following the definition of the smoothness indicators in \cite{JiangShu, Shu}, which measure the regularity of the polynomial approximation $p_k(x)$ over some interval, the smoothness indicators are defined as 
$$
   \beta_k = \sum_{l=1}^2 \Delta x^{2l-1} \int_{x_i}^{x_{i+1}} \left( \frac{d^l}{d x^l} p_k(x) \right)^2 dx,
$$
which gives
\begin{equation} \label{eq:smooth_indicator_Liu}
\begin{aligned}
 \beta_0 &= \frac{13}{12} \left( b_{i-2} - 3 b_{i-1} + 3 b_i - b_{i+1} \right)^2 + \frac{1}{4} \left( b_{i-2} - 5 b_{i-1} + 7 b_i - 3 b_{i+1} \right)^2, \\
 \beta_1 &= \frac{13}{12} \left( b_{i-1} - 3 b_i + 3 b_{i+1} - b_{i+2} \right)^2 + \frac{1}{4} \left( b_{i-1} - b_i - b_{i+1} + b_{i+2} \right)^2,\\
 \beta_2 &= \frac{13}{12} \left( b_i - 3 b_{i+1} + 3 b_{i+2} - b_{i+3} \right)^2 + \frac{1}{4} \left( -3 b_i + 7 b_{i+1} - 5 b_{i+2} + b_{i+3} \right)^2.
\end{aligned}
\end{equation}
The integration over the interval $[x_i, x_{i+1}]$ is performed to satisfy the symmetry property of the parabolic equation, and the factor $\Delta x^{2l-1}$ is introduced to remove any $\Delta x$ dependency in the derivatives. 

Depending on the linear weights \eqref{eq:gamma} and the smoothness indicators \eqref{eq:smooth_indicator_Liu}, the nonlinear weights could be defined for the WENO approximation.
In \cite{Liu}, Liu et al. derived the sufficient conditions to attain sixth order accuracy in smooth regions,
\begin{gather} 
 \omega_0 - \omega_2 = O(\Delta x^4), \label{eq:nonlinear_weights_condition} \\
 \omega_k - d_k = O(\Delta x^3).      \label{eq:linear_nonlinear_condition}
\end{gather}

In \cite{Liu}, the nonlinear positive and negative weights $\omega^{\pm}_k$ are defined as
\begin{equation} \label{eq:weights_Liu_pm}
 \omega^{\pm}_k = \frac{\alpha^{\pm}_k}{\sum^2_{l=0} \alpha^{\pm}_l},~~\alpha^{\pm}_k = \frac{\gamma^{\pm}_k}{(\beta_k + \epsilon)^2},~~k=0,1,2,
\end{equation}
where $\epsilon>0$ is known to avoid the denominator becoming zero. 
Based on the relation \eqref{eq:linear_weight_relation} for the linear weights, the nonlinear weights are defined by
\begin{equation} \label{eq:weights}
 \omega_k = \sigma^+ \omega^+_k - \sigma^- \omega^-_k.
\end{equation}
However, the nonlinear weights defined in \eqref{eq:weights_Liu_pm} and \eqref{eq:weights} give
$$
   \omega_k - d_k = O(\Delta x),
$$
where the condition \eqref{eq:linear_nonlinear_condition} is not satisfied.
To increase the accuracy of the nonlinear weights, the mapped function in \cite{Henrick} is employed: 
$$
   g_k(\omega) = \frac{\omega (d_k + d_k^2 - 3 d_k \omega + \omega^2)}{d_k^2 + \omega (1 - 2d_k)},~~k = 0,1,2.
$$
The final nonlinear weights are formulated as
$$
   \omega^{\LSZ}_k = \frac{\alpha_k}{\sum^2_{l=0} \alpha_l},~~\alpha_k = g_k(\omega_k),~~k=0,1,2.
$$
It is shown in \cite{Liu} with Taylor expansion that 
\begin{gather*} 
 \omega^{\LSZ}_0 - \omega^{\LSZ}_2 = O(\Delta x^4), \\
 \omega^{\LSZ}_k - d_k = O(\Delta x^3). 
\end{gather*}
So both conditions \eqref{eq:nonlinear_weights_condition} and \eqref{eq:linear_nonlinear_condition} are satisfied.
\begin{remark} \label{rmk:epsilon_LSZ}
{\rm As $\epsilon = 10^{-6}$ is used in \cite{Liu}, it causes some small-scale oscillations around the sharp interfaces, e.g., in Examples \ref{ex:Barenblatt}, \ref{ex:two_box} and \ref{ex:PME_2d}, Section \ref{sec:nr}, and the NAN values in some computer systems, e.g., in Example \ref{ex:PME_2d}, Section \ref{sec:nr}. 
The value of $\epsilon$ is thus replaced by $10^{-10}$ in the WENO-LSZ scheme for all numerical experiments in Section \ref{sec:nr}, so that some oscillations are smoothed and there is none NAN value for all tested computer systems.}
\end{remark}

In \cite{Hajipour}, the MWENO scheme were proposed with Z-type nonlinear weights \cite{Borges}, where the global smoothness indicator is supposed to give higher order, which implies that the lower order terms happen to cancel out if the function is smooth in the stencil.
The global smoothness indicator $\tau$ is simply the absolute difference between $\beta_0$ and $\beta_2$,
$$
   \tau = |\beta_0 - \beta_2|,
$$
and the nonlinear positive and negative weights are defined as 
\begin{equation} \label{eq:weights_MWENO_pm}
 \omega^{\pm}_k = \frac{\alpha^{\pm}_k}{\sum^2_{l=0} \alpha^{\pm}_l},~~\alpha^{\pm}_k = \gamma^{\pm}_k \left( 1 + \left( \frac{\tau}{\beta_k + \epsilon} \right)^2 \right),~~k=0,1,2,
\end{equation}
with $\gamma^{\pm}_k$ in \eqref{eq:gamma} and $\epsilon = 10^{-30}$.
As defined in \eqref{eq:weights}, the MWENO nonlinear weights are 
$$
   \omega^{\MWENO}_k = \sigma^+ \omega^+_k - \sigma^- \omega^-_k,
$$
which satisfy the sufficient conditions for sixth order accuracy in smooth regions as shown in \cite{Hajipour}.

Hence the WENO numerical flux is
$$
   \hat{g}_{i+1/2} = \sum_{k=0}^2 \omega_k \hat{g}^k_{i+1/2},
$$
where $\hat{g}^k_{i+1/2},~k=0,1,2$ are given by \eqref{eq:numerical_flux_substencil}.
Then the semi-discrete finite difference WENO scheme of the conservation form is 
\begin{equation} \label{eq:FD_WENO_scheme}
 \frac{du_{i}(t)}{dt} = \frac{\hat{g}_{i+1/2} - \hat{g}_{i-1/2}}{\Delta x^2}.
\end{equation}
\begin{remark}
{\rm A major difference between the WENO approximations to the first and second derivatives is that in one stencil, two numerical fluxes $\hat{f}^+_{i-1/2}$ and $\hat{f}^-_{i+1/2}$ at the respective points $x_{i-1/2}$ and $x_{i+1/2}$ need evaluating for the first derivative, while only one numerical flux $\hat{g}_{i+1/2}$ at $x_{i+1/2}$ is estimated for the second derivative.}
\end{remark}

\section{Central WENO approximation to the second derivative} \label{sec:cweno}

Our goal is to obtain a convex combination of the numerical fluxes $\hat{g}^k_{i+1/2}$ in \eqref{eq:numerical_flux_substencil} as a new approximation to $g_{i+1/2}$ in \eqref{eq:parabolic_discrete_h}.
Motivated by the compact central WENO schemes for hyperbolic conservation laws \cite{Levy}, a centered polynomial $p_C(x)$ is introduced for the WENO approximation, such that all linear weights are positive without the concern about dealing with negative linear weights.
To conform to the notation in \cite{Levy}, we set $p_{\OPT}(x) = p(x),~p_\Le (x) = p_0(x),~p_\Mi (x) = p_1(x),~p_\Ri (x) = p_2(x)$ with $p(x)$ \eqref{eq:polynomial} and $p_0(x),~p_1(x),~p_2(x)$ \eqref{eq:subpolynomials}, and apply this setting to every related term in the previous section accordingly.
The centered polynomial $p_C(x)$ is constructed by the following relation
\begin{equation} \label{eq:cpolynomial_relation}
 p_{\OPT}(x) = C_\Le p_\Le(x) + C_\Mi p_\Mi(x) + C_\Ri p_\Ri(x) + C_\C p_\C(x),~~\sum_k C_k = 1,~~k \in \{ \Le, \Mi, \Ri, \C \},
\end{equation}
where $C_\Le,~C_\Mi,~C_\Ri$ and $C_\C$ are positive constants.
It is required that $C_\Le=C_\Ri$, as discussed further in the following subsection.
With different combinations of $(C_\Le,~C_\Mi,~C_\Ri)$ attempted, we pick out $(C_\Le,~C_\Mi,~C_\Ri) = \left( \frac{1}{6},~\frac{1}{3},~\frac{1}{6} \right)$.
Then
$$
   p_\C(x) = 3 p_{\OPT}(x) - \frac{1}{2} p_\Le(x) - p_\Mi(x) - \frac{1}{2} p_\Ri(x).
$$
The evaluation of the polynomial at $x = x_{i+1/2}$ gives rise to the central numerical flux $\hat{g}^\C_{i+1/2}$
\begin{equation} \label{eq:numerical_flux_central}
 \hat{g}^\C_{i+1/2} = p_\C(x_{i+1/2}) = - \frac{3}{40} b_{i-2} + \frac{11}{24} b_{i-1} - 2b_i + 2b_{i+1} - \frac{11}{24} b_{i+2} + \frac{3}{40} b_{i+3},
\end{equation}
and the Taylor series expansion shows that
\begin{equation} \label{eq:numerical_flux_central_plus_Taylor}
 \hat{g}^\C_{i+1/2} = g_{i+1/2} + \frac{23}{360} h_i^{(5)} \Delta x^5 + O(\Delta x^6).
\end{equation}
From \eqref{eq:cpolynomial_relation}, we have
\begin{equation} \label{eq:numerical_flux_central_linear_weights}
 \hat{g}_{i+1/2}^{\FD} = \sum_k C_k \hat{g}^k_{i+1/2},~~k \in \{ \Le, \Mi, \Ri, \C \}.
\end{equation}
It is straightforward to obtain the corresponding $\hat{g}^\C_{i-1/2}$ and the relation between $\hat{g}_{i-1/2}^{\FD}$ and $\hat{g}^k_{i-1/2}$ with every index shifted by -1. 
Similarly, expanding in Taylor series gives
\begin{equation} \label{eq:numerical_flux_central_minus_Taylor}
 \hat{g}^\C_{i-1/2} = g_{i-1/2} + \frac{23}{360} h_i^{(5)} \Delta x^5 + O(\Delta x^6).
\end{equation}

We define the central smoothness indicator as  
$$
   \beta_\C = \sum_{l=1}^4 \Delta x^{2l-1} \int_{x_i}^{x_{i+1}} \left( \frac{d^l}{d x^l} p_{\OPT}(x) \right)^2 dx,
$$
where the polynomial $p_\C (x)$ is not used but replaced by $p_{\OPT}(x)$, and the order of the derivative is up to 4.
After some algebra, the central smoothness indicator could be written as
\begin{equation} \label{eq:smooth_indicator_central}
\begin{aligned}
 \beta_\C = & \frac{4273}{20160} \left( b_{i-2} -5	b_{i-1} + 10 b_i - 10 b_{i+1} + 5 b_{i+2} - b_{i+3} \right)^2 + \\
            & \frac{29}{345600} \left( 5 b_{i-2} + 11 b_{i-1} - 70 b_i + 94 b_{i+1} - 47 b_{i+2} + 7 b_{i+3} \right)^2 + \\
            & \frac{1}{3600} \left( 35 b_{i-2} - 139 b_{i-1} + 230 b_i - 206 b_{i+1} + 103 b_{i+2} - 23 b_{i+3} \right)^2 + \\
            & \frac{1}{576} \left( 7 b_{i-2} - 51 b_{i-1} + 134 b_i - 166 b_{i+1} + 99 b_{i+2} - 23 b_{i+3} \right)^2 + \\           
            & \frac{1}{2304} \left( 7 b_{i-2} - 56 b_{i-1} + 106 b_i - 76 b_{i+1} + 23 b_{i+2} - 4 b_{i+3} \right)^2 + \\
            & \frac{1}{9216} \left( 65 b_{i-2} - 353 b_{i-1} + 690 b_i - 602 b_{i+1} + 221 b_{i+2} - 21 b_{i+3} \right)^2 + \\
            & \frac{1}{9216} \left( 23 b_{i-2} - 63 b_{i-1} - 34 b_i + 186 b_{i+1} - 133 b_{i+2} + 21 b_{i+3} \right)^2 + \\
            & \frac{1}{2304} \left( 13 b_{i-2} - 28 b_{i-1} + 30 b_i - 28 b_{i+1} + 13 b_{i+2} \right)^2 + \\
            & \frac{2}{15} \left( b_{i-2} - 4 b_{i-1} + 6 b_i - 4 b_{i+1} + b_{i+2} \right)^2 +           
              \frac{1}{1152} \left( b_{i-2} - 12 b_{i-1} + 22 b_i - 12 b_{i+1} + b_{i+2} \right)^2.
\end{aligned}
\end{equation}
We set the new global smoothness indicator $\tau_6$ on the stencil $S^6$ as
$$
   \tau_6 = \left| \beta_\C - \frac{1}{24} ( 5 \beta_\Le + 14 \beta_\Mi + 5 \beta_\Ri) \right|.
$$
The nonlinear weights are defined by
\begin{equation} \label{eq:weights_CWENO}
 \omega^{\CWENO}_k = \frac{\alpha_k}{\sum_l \alpha_l},~~\alpha_k = C_k \left( 1 + \left( \frac{\tau_6}{\beta_k + \epsilon} \right)^p \right),~~k,~l \in \{ \Le, \Mi, \Ri, \C \},
\end{equation}
with $C_\Le = C_\Ri = \frac{1}{6},~C_\Mi = C_\C= \frac{1}{3}$ and $\epsilon = 10^{-40}$.
The free parameter $p$ is important to achieve sixth order accuracy in smooth regions, as well as control the amount of numerical dissipation.
The choice of $p$ will be discussed below.
We end up with the CWENO numerical flux
\begin{equation} \label{eq:CWENO_numerical_flux}
 \hat{g}_{i+1/2} = \sum_k \omega^{\CWENO}_k \hat{g}^k_{i+1/2},~~k \in \{ \Le, \Mi, \Ri, \C \}.
\end{equation}
Note that the central numerical flux $\hat{g}^\C_{i-1/2}$ is used in smooth regions.
Otherwise its contribution vanishes and the WENO numerical flux is determined by the nonlinear weight(s) corresponding to the smoothness indicators of smaller magnitude.

\subsection{Spatial sixth order accuracy in smooth regions}

We next consider the sufficient conditions of the finite difference WENO scheme \eqref{eq:FD_WENO_scheme} with the new numerical flux \eqref{eq:CWENO_numerical_flux} so as to maintain sixth order accuracy in smooth regions.
Let
$$
   \hat{g}_{i \pm 1/2} = \sum_k \omega^{\pm}_k \hat{g}^k_{i \pm 1/2},~~k \in \{ \Le, \Mi, \Ri, \C \}.
$$
Here we drop the superscript CWENO in \eqref{eq:CWENO_numerical_flux} to simplify the notation.
The superscripts $\pm$ in the nonlinear weights $\omega^{\pm}_k$ represent two different stencils, with $+$ for $\{ x_{i-2}, \cdots, x_{i+3} \}$ and $-$ for $\{ x_{i-1}, \cdots, x_{i+2} \}$.
The nonlinear weights $\omega^{\pm}_k$ in this subsection are not the nonlinear positive and negative weights $\omega^{\pm}_k$ \eqref{eq:weights_Liu_pm} and \eqref{eq:weights_MWENO_pm} in Section \ref{sec:weno}.
From the relation \eqref{eq:numerical_flux_central_linear_weights}, the numerical flux \eqref{eq:CWENO_numerical_flux} can be rewritten as
$$
   \hat{g}_{i+1/2} = \sum_k C_k \hat{g}^k_{i+1/2} + \sum_k ( \omega^+_k - C_k ) \hat{g}^k_{i+1/2} 
                   = \hat{g}^{\FD}_{i+1/2} + \sum_k ( \omega^+_k - C_k ) \hat{g}^k_{i+1/2}.
$$
We expand the last term by using \eqref{eq:numerical_flux_substencil_Taylor} and \eqref{eq:numerical_flux_central_plus_Taylor},
\begin{align*}
 \sum_k ( \omega^+_k - C_k ) \hat{g}^k_{i+1/2} = 
 {} & ( \omega^+_\Le - C_\Le ) \left[ g_{i+1/2} + \frac{1}{12} h_i^{(4)} \Delta x^4 + O(\Delta x^5) \right] + \\
 {} & ( \omega^+_\Mi - C_\Mi ) \left[ g_{i+1/2} - \frac{1}{90} h_i^{(5)} \Delta x^5 + O(\Delta x^6) \right] + \\
 {} & ( \omega^+_\Ri - C_\Ri ) \left[ g_{i+1/2} - \frac{1}{12} h_i^{(4)} \Delta x^4 + O(\Delta x^5) \right] + \\
 {} & ( \omega^+_\C - C_\C ) \left[ g_{i+1/2} + \frac{23}{360} h_i^{(5)} \Delta x^5 + O(\Delta x^6) \right] \\
 = {} & g_{i+1/2} \sum_k(\omega^+_k - C_k)+\frac{1}{12} h_i^{(4)} \Delta x^4 (\omega^+_\Le - \omega^+_\Ri)+\sum_k (\omega^+_k - C_k) O(\Delta x^5)\\
 = {} & \frac{1}{12} h_i^{(4)} \Delta x^4 (\omega^+_\Le - \omega^+_\Ri) + \sum_k ( \omega^+_k - C_k ) O(\Delta x^5).
\end{align*}
Then
$$
  \hat{g}_{i+1/2} = \hat{g}^{\FD}_{i+1/2}+ \frac{1}{12} h_i^{(4)} \Delta x^4 (\omega^+_\Le - \omega^+_\Ri)+ \sum_k ( \omega^+_k - C_k ) O(\Delta x^5). 
$$
Similarly, with the help of \eqref{eq:numerical_flux_substencil_Taylor} and \eqref{eq:numerical_flux_central_minus_Taylor}, we find that
$$
  \hat{g}_{i-1/2} = \hat{g}^{\FD}_{i-1/2}+ \frac{1}{12} h_i^{(4)} \Delta x^4 (\omega^-_\Le - \omega^-_\Ri)+ \sum_k ( \omega^-_k - C_k ) O(\Delta x^5). 
$$
By \eqref{eq:numerical_flux_plus_Taylor} and \eqref{eq:numerical_flux_minus_Taylor}, we have
\begin{align*}
 \frac{\hat{g}_{i+1/2} - \hat{g}_{i-1/2}}{\Delta x^2} = {} & \frac{g_{i+1/2} - g_{i-1/2}}{\Delta x^2} + O(\Delta x^6) + \frac{1}{12} h_i^{(4)} \Delta x^2 (\omega^+_\Le - \omega^+_\Ri) - \frac{1}{12} h_i^{(4)} \Delta x^2 (\omega^-_\Le - \omega^-_\Ri) \\
                                                        {} & + \sum_k ( \omega^+_k - C_k ) O(\Delta x^3) - \sum_k ( \omega^-_k - C_k ) O(\Delta x^3).
\end{align*} 
Thus the sufficient conditions for sixth order accuracy are given by
\begin{gather} 
 \omega_\Le - \omega_\Ri = O(\Delta x^4), \label{eq:CWENO_nonlinear_weights_condition} \\
 \omega_k - C_k = O(\Delta x^3),          \label{eq:CWENO_linear_nonlinear_condition}
\end{gather}
where the superscripts are dropped, meaning that the nonlinear weights $\omega_k$ for each stencil $S^6$ are supposed to satisfy both conditions in smooth regions for sixth order accuracy.

Expanding the smoothness indicators $\beta_k,~k \in \{ \Le, \Mi, \Ri, \C \}$ \eqref{eq:smooth_indicator_Liu} and \eqref{eq:smooth_indicator_central} in Taylor series at $x = x_{i+1/2}$, we obtain
\begin{align*}
 \beta_\Le = b''^2_{i+1/2} \Delta x^4 &+ \left( \frac{13}{12} b'''^2_{i+1/2} - \frac{7}{12} b''_{i+1/2} b^{(4)}_{i+1/2} \right) \Delta x^6 \\
             &+ \left( - \frac{13}{6} b'''_{i+1/2} b^{(4)}_{i+1/2} + \frac{1}{2} b''_{i+1/2} b^{(5)}_{i+1/2} \right) \Delta x^7 + O (\Delta x^8), \\
 \beta_\Mi = b''^2_{i+1/2} \Delta x^4 &+ \left( \frac{13}{12} b'''^2_{i+1/2} + \frac{5}{12} b''_{i+1/2} b^{(4)}_{i+1/2} \right) \Delta x^6 + O (\Delta x^8), \\
 \beta_\Ri = b''^2_{i+1/2} \Delta x^4 &+ \left( \frac{13}{12} b'''^2_{i+1/2} - \frac{7}{12} b''_{i+1/2} b^{(4)}_{i+1/2} \right) \Delta x^6 \\
               &+ \left( \frac{13}{6} b'''_{i+1/2} b^{(4)}_{i+1/2} - \frac{1}{2} b''_{i+1/2} b^{(5)}_{i+1/2} \right) \Delta x^7 + O (\Delta x^8), \\
 \beta_\C  = b''^2_{i+1/2} \Delta x^4 &+ \frac{13}{12} b'''^2_{i+1/2} \Delta x^6 + O (\Delta x^8).
\end{align*}

If there is no inflection (or undulation) point at $x_{i+1/2}$, i.e., the second derivative is nonzero, then $\tau_6 = O (\Delta x^8)$.
Since $\tau_6$ is of order $O (\Delta x^8)$ and each $\beta_k$ is of order $O (\Delta x^4)$, one can find that
$$
   \left( \frac{\tau_6}{\beta_k} \right)^p = O (\Delta x^{4p}),
$$
by setting $\epsilon = 0$ in the Taylor expansion analysis.
From the definitions \eqref{eq:weights_CWENO},
\begin{equation} \label{eq:linear_nonlinear_relation}
 \omega^{\CWENO}_k = C_k + O (\Delta x^{4p}). 
\end{equation}
The minimum value $p$ to satisfy both conditions \eqref{eq:CWENO_nonlinear_weights_condition} and \eqref{eq:CWENO_linear_nonlinear_condition} is $p = 1$.
Note that the condition \eqref{eq:CWENO_nonlinear_weights_condition} combined with \eqref{eq:linear_nonlinear_relation} explains the requirement $C_\Le = C_\Ri$.

Now we consider the convergence behavior of the nonlinear weights when there exists an inflection point at $x_{i+1/2}$, that is, the second derivative is zero but the third derivative is nonzero.
Then it can be verified through the Taylor expansion analysis above that  
$$
   \omega^{\CWENO}_k = C_k + O (\Delta x^{2p}),
$$
and it is clear that $p=2$ is the minimum value to maintain sixth order accuracy.

As pointed out by Borges et al. \cite{Borges}, increasing the value of $p$ amplifies the numerical dissipation around the discontinuities.
We then choose $p = 1$ in this paper even if it does not satisfy the sufficient conditions \eqref{eq:CWENO_nonlinear_weights_condition} and \eqref{eq:CWENO_linear_nonlinear_condition} at the inflection points.
However, our numerical experiments in the next section show that it still provide sixth order accuracy overall.

\section{Numerical results} \label{sec:nr}
This section presents some numerical experiments to demonstrate the performance of the proposed central WENO scheme and compare with the WENO-LSZ and MWENO schemes.
We examine the accuracy of the WENO schemes for one- and two-dimensional heat equations in terms of $L^1,~L^2$ and $L^{\infty}$ error norms:
\begin{align*}
 & L^1 = \frac{1}{N+1} \sum_{i=0}^N \left| u_i(T) - u(x_i, T) \right|, \\
 & L^2 = \sqrt{\frac{1}{N+1} \sum_{i=0}^N \left( u_i(T) - u(x_i, T) \right)^2}, \\
 & L^{\infty} = \max_{0 \leqslant i \leqslant N} \left| u_i(T) - u(x_i, T) \right|,
\end{align*}
where $u(x_i, T)$ denotes the exact solution and $u_i(T)$ is the numerical approximation at the final time $t = T$.
The rest numerical experiments show the resolution of the numerical solutions with the WENO-LSZ, MWENO and central WENO schemes.
For time discretization, we use the explicit third-order total variation diminishing Runge-Kutta method \cite{ShuOsherI}
\begin{align*}
 u^{(1)} &= u^n + \Delta t L(u^n), \\
 u^{(2)} &= \frac{3}{4} u^n + \frac{1}{4} u^{(1)} + \frac{1}{4} \Delta t L \left( u^{(1)} \right), \\
 u^{n+1} &= \frac{1}{3} u^n + \frac{2}{3} u^{(2)} + \frac{2}{3} \Delta t L \left( u^{(2)} \right),
\end{align*}
where $L$ is the spatial operator.
We follow the CFL condition in \cite{Liu} to set $\cfl = 0.4$ unless otherwise stated.
The central WENO scheme in Section \ref{sec:cweno} is termed as CWENO-DZ with $p=1$.
We choose $\epsilon = 10^{-40}$ for the CWENO-DZ scheme whereas $\epsilon = 10^{-15}$ is set for WENO-LSZ as explained in Remark \ref{rmk:epsilon_LSZ} and $\epsilon = 10^{-30}$ for MWENO as in \cite{Hajipour}.

\subsection{One-dimensional numerical examples}

\begin{example} \label{ex:heat_1d}
We test the accuracy of those WENO schemes for the one-dimensional heat equation
$$
   u_t = u_{xx},~~-\pi \leqslant x \leqslant \pi,~~t>0
$$
with the following initial data
$$
   u(x,0) = \sin(x), 
$$
and the periodic boundary condition. 
The exact solution is given by
$$
   u(x,t)= e^{-t} \sin(x).
$$
The numerical solution is computed up to the time $T=2$ with the time step $\Delta t = \cfl \cdot \Delta x^2$. 
We present the $L_1, L_2$ and $L_{\infty}$ errors versus $N$, as well as the order of accuracy, for the WENO-LSZ, MWENO, CWENO-DZ schemes in Tables \ref{tab:heat_1d_L1}, \ref{tab:heat_1d_L2} and \ref{tab:heat_1d_Linf}, respectively. 
It is clear that the expected order of accuracy is achieved for all schemes.
Although the errors of the CWENO-DZ scheme are larger than WENO-LSZ for $N=10$, CWENO-DZ yields the most accurate results as $N$ increases.
\end{example}

\begin{table}[h!]
\centering
\caption{$L_1$ error and order of accuracy for Example \ref{ex:heat_1d}.}      
\begin{tabular}{clcrlcrlc} 
\hline  
N & \multicolumn{2}{l}{WENO-LSZ} & & \multicolumn{2}{l}{MWENO} & & \multicolumn{2}{l}{CWENO-DZ}  \\ 
    \cline{2-3}                      \cline{5-6}                     \cline{8-9}                    
  & Error & Order                & & Error & Order               & & Error & Order              \\
\hline 
10  & 6.31E-6  & --     & & 3.17E-5  & --     & & 4.15E-5   & --      \\  
20  & 1.41E-7  & 5.4883 & & 2.16E-7  & 7.1985 & & 1.77E-8   & 11.1951  \\  
40  & 2.27E-9  & 5.9514 & & 2.36E-9  & 6.5124 & & 1.94E-9   & 3.1896  \\
80  & 3.54E-11 & 6.0028 & & 3.55E-11 & 6.0562 & & 3.47E-11  & 5.8050  \\ 
160 & 5.70E-13 & 5.9582 & & 5.70E-13 & 5.9613 & & 5.69E-13  & 5.9304  \\  
\hline
\end{tabular}
\label{tab:heat_1d_L1}
\end{table}

\begin{table}[h!]
\centering
\caption{$L_2$ error and order of accuracy for Example \ref{ex:heat_1d}.}      
\begin{tabular}{clcrlcrlc} 
\hline  
N & \multicolumn{2}{l}{WENO-LSZ} & & \multicolumn{2}{l}{MWENO} & & \multicolumn{2}{l}{CWENO-DZ} \\ 
    \cline{2-3}                      \cline{5-6}                     \cline{8-9}                    
  & Error & Order                & & Error & Order               & & Error & Order              \\
\hline 
10  & 7.50E-6  & --     & & 3.79E-5  & --     & & 4.91E-5   & --     \\  
20  & 1.61E-7  & 5.5422 & & 2.47E-7  & 7.2580 & & 2.11E-8   & 11.1843 \\  
40  & 2.56E-9  & 5.9742 & & 2.66E-9  & 6.5387 & & 2.21E-9   & 3.2551 \\
80  & 3.96E-11 & 6.0136 & & 3.97E-11 & 6.0664 & & 3.89E-11  & 5.8281 \\ 
160 & 6.35E-13 & 5.9633 & & 6.35E-13 & 5.9663 & & 6.34E-13  & 5.9391 \\   
\hline
\end{tabular}
\label{tab:heat_1d_L2}
\end{table}

\begin{table}[h!]
\centering
\caption{$L_\infty$ error and order of accuracy for Example \ref{ex:heat_1d}.}      
\begin{tabular}{clcrlcrlc} 
\hline  
N & \multicolumn{2}{l}{WENO-LSZ} & & \multicolumn{2}{l}{MWENO} & & \multicolumn{2}{l}{CWENO-DZ} \\ 
    \cline{2-3}                      \cline{5-6}                     \cline{8-9}                    
  & Error & Order                & & Error & Order               & & Error & Order              \\
\hline 
10  & 1.01E-5  & --     & & 5.22E-5  & --     & & 6.43E-5   & --     \\  
20  & 2.31E-7  & 5.4501 & & 3.54E-7  & 7.2038 & & 3.74E-8   & 10.7476 \\  
40  & 3.66E-9  & 5.9780 & & 3.80E-9  & 6.5398 & & 3.21E-9   & 3.5424 \\
80  & 5.64E-11 & 6.0217 & & 5.65E-11 & 6.0722 & & 5.54E-11  & 5.8565 \\ 
160 & 9.01E-13 & 5.9677 & & 9.02E-13 & 5.9702 & & 8.99E-13  & 5.9454 \\   
\hline
\end{tabular}
\label{tab:heat_1d_Linf}
\end{table}

\begin{figure}[h!]
\centering
\includegraphics[width=0.32\textwidth]{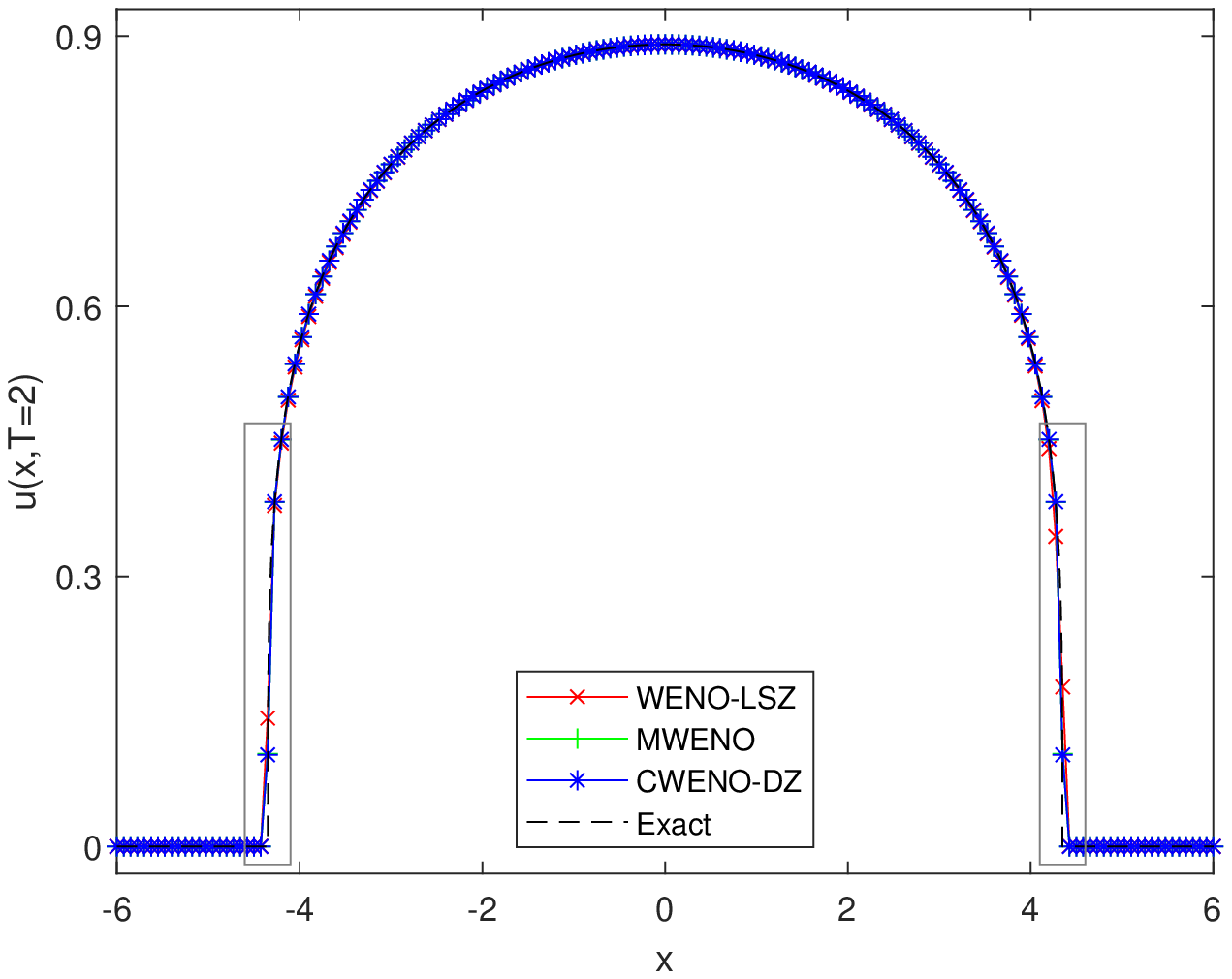}
\includegraphics[width=0.32\textwidth]{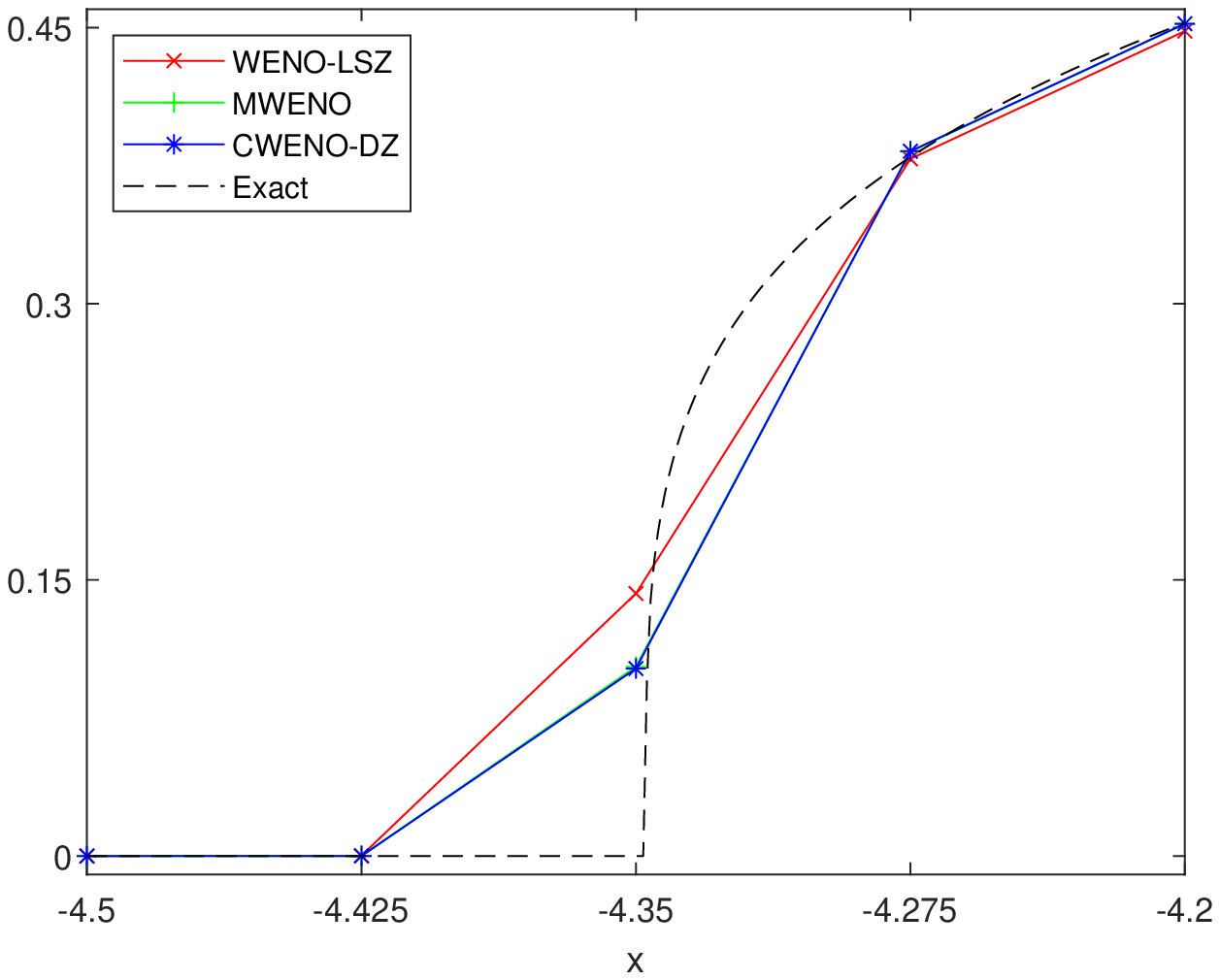}
\includegraphics[width=0.32\textwidth]{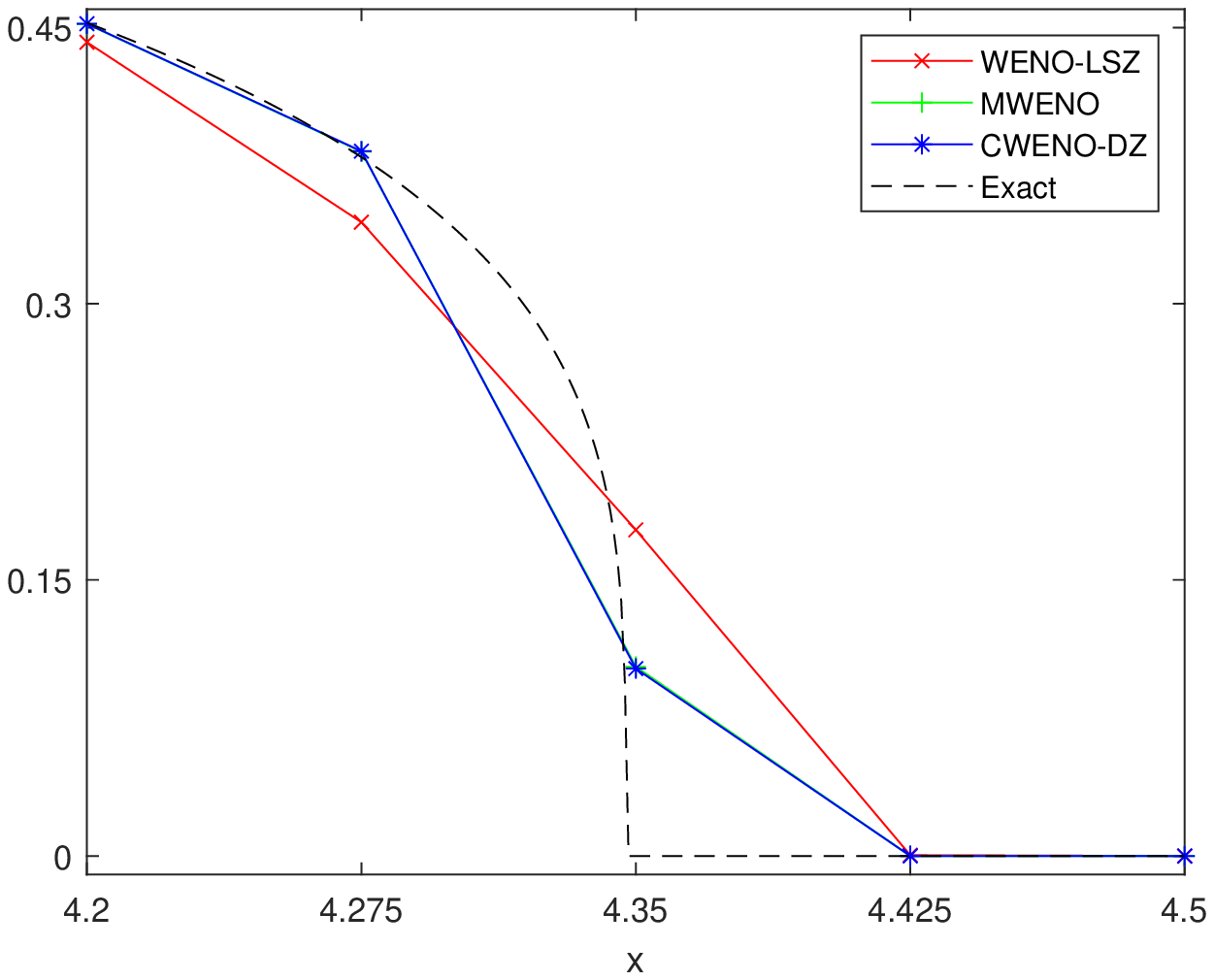}
\captionof{figure}{Barenblatt solution profiles for Example \ref{ex:Barenblatt} with $m=5$ at $T=2$ (left), close-up view of the solutions in the boxes on the left/right (middle/right) computed by WENO-LSZ (red), MWENO (green) and CWENO-DZ (blue) with $N=160$. The dashed black lines are the exact solution.}
\label{fig:Barenblatt_m5}
\end{figure}

\begin{figure}[h!]
\centering
\includegraphics[width=0.32\textwidth]{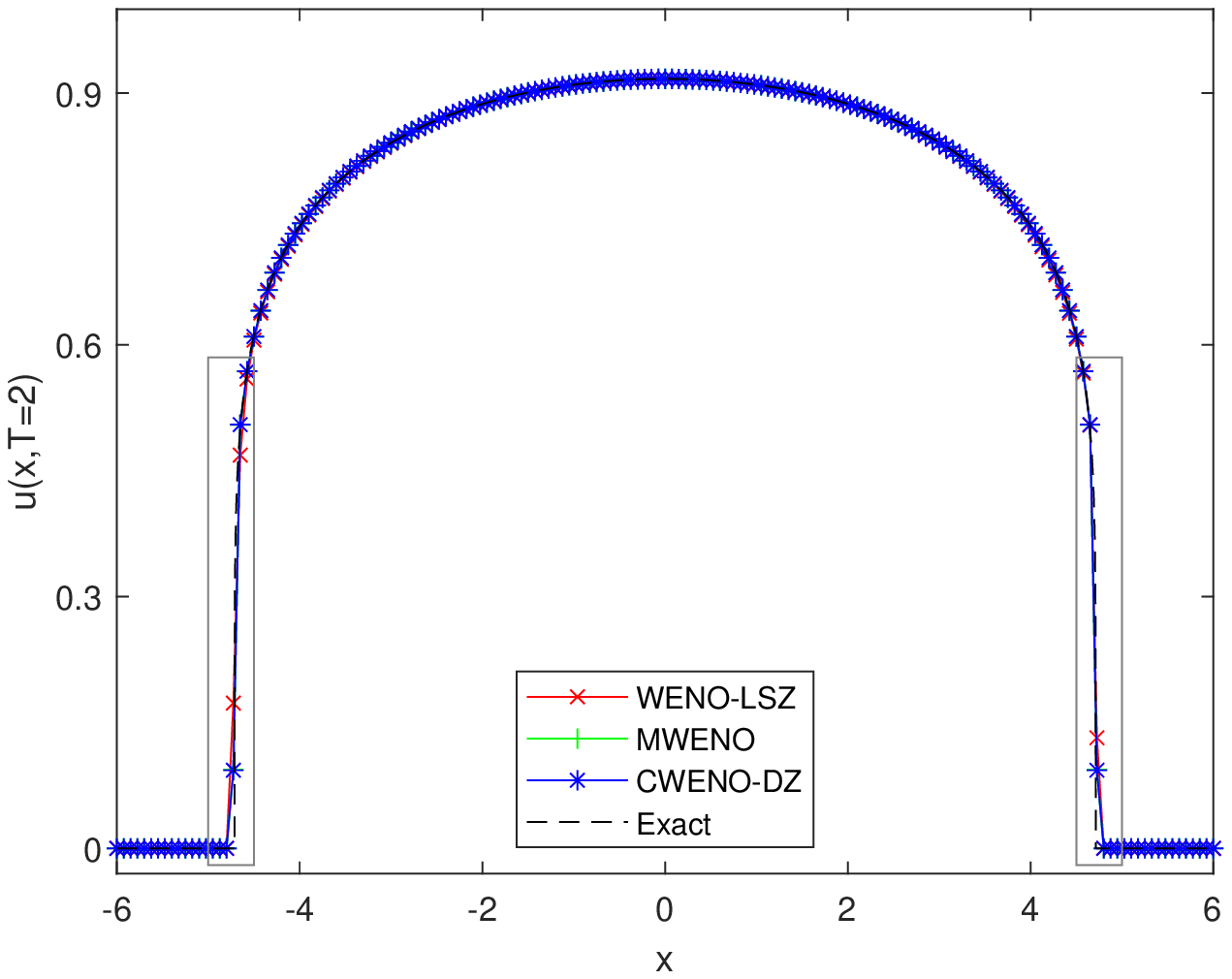}
\includegraphics[width=0.32\textwidth]{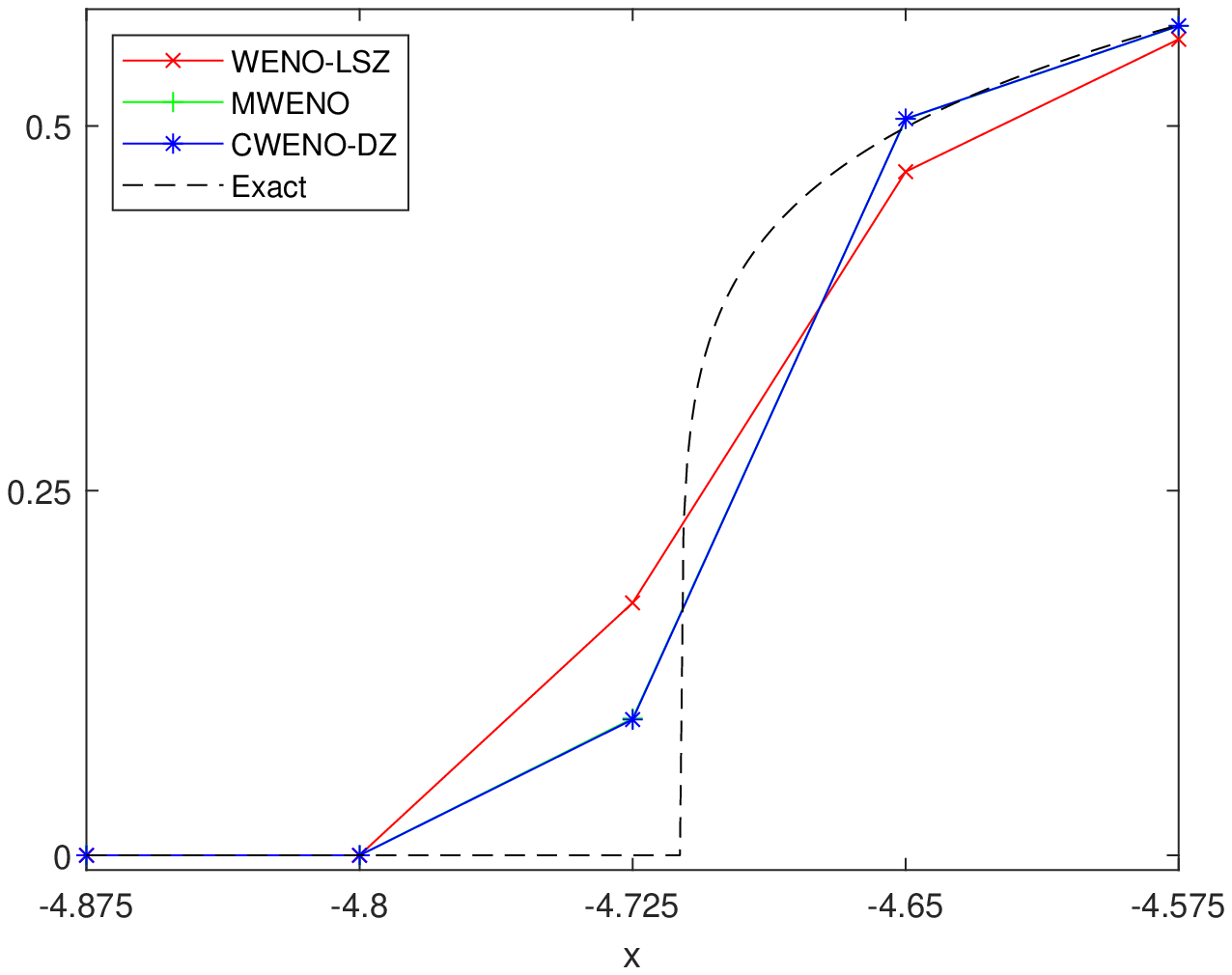}
\includegraphics[width=0.32\textwidth]{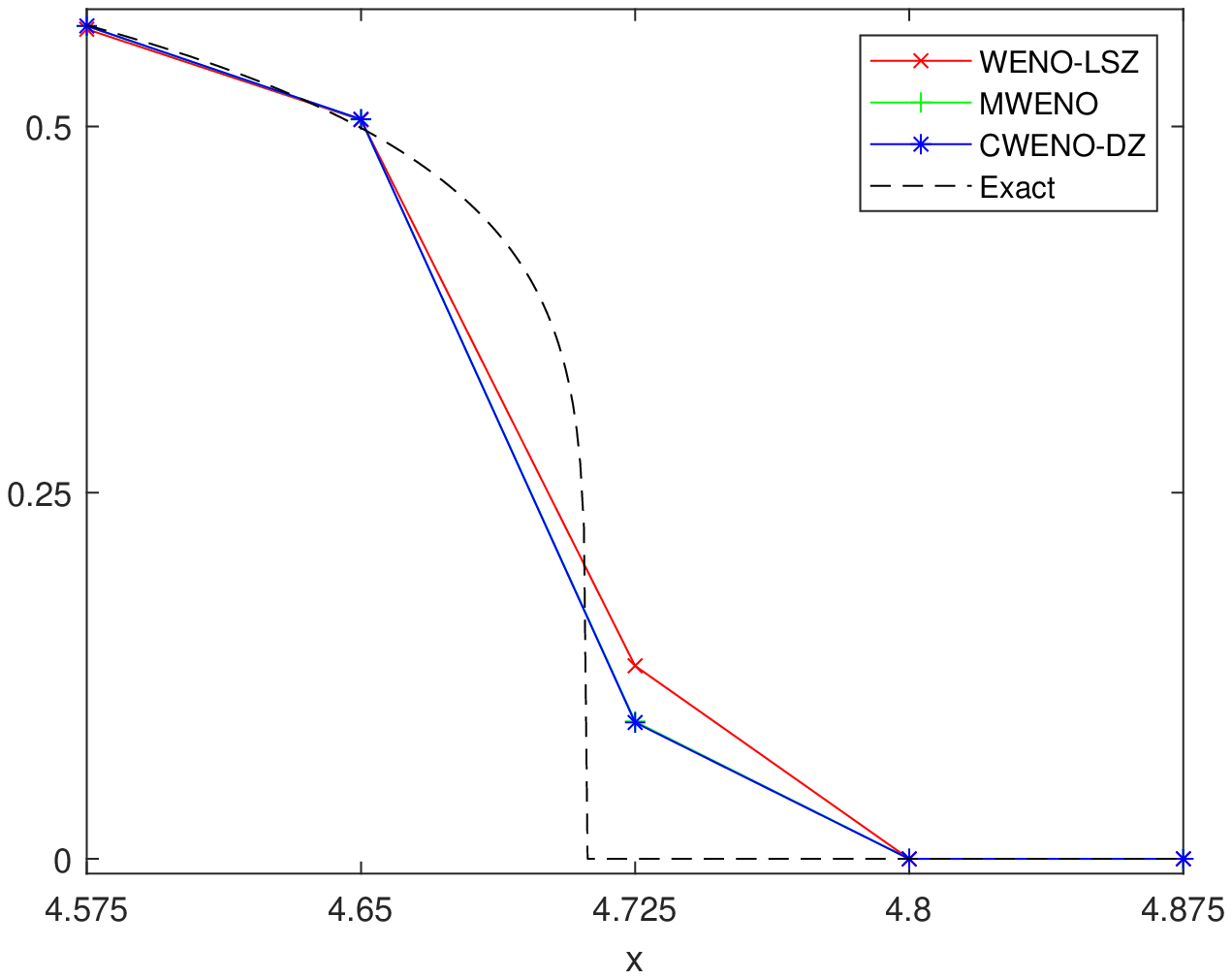}
\captionof{figure}{Barenblatt solution profiles for Example \ref{ex:Barenblatt} with $m=7$ at $T=2$ (left), close-up view of the solutions in the boxes on the left/right (middle/right) computed by WENO-LSZ (red), MWENO (green) and CWENO-DZ (blue) with $N=160$. The dashed black lines are the exact solution.}
\label{fig:Barenblatt_m7}
\end{figure}

\begin{figure}[h!]
\centering
\includegraphics[width=0.32\textwidth]{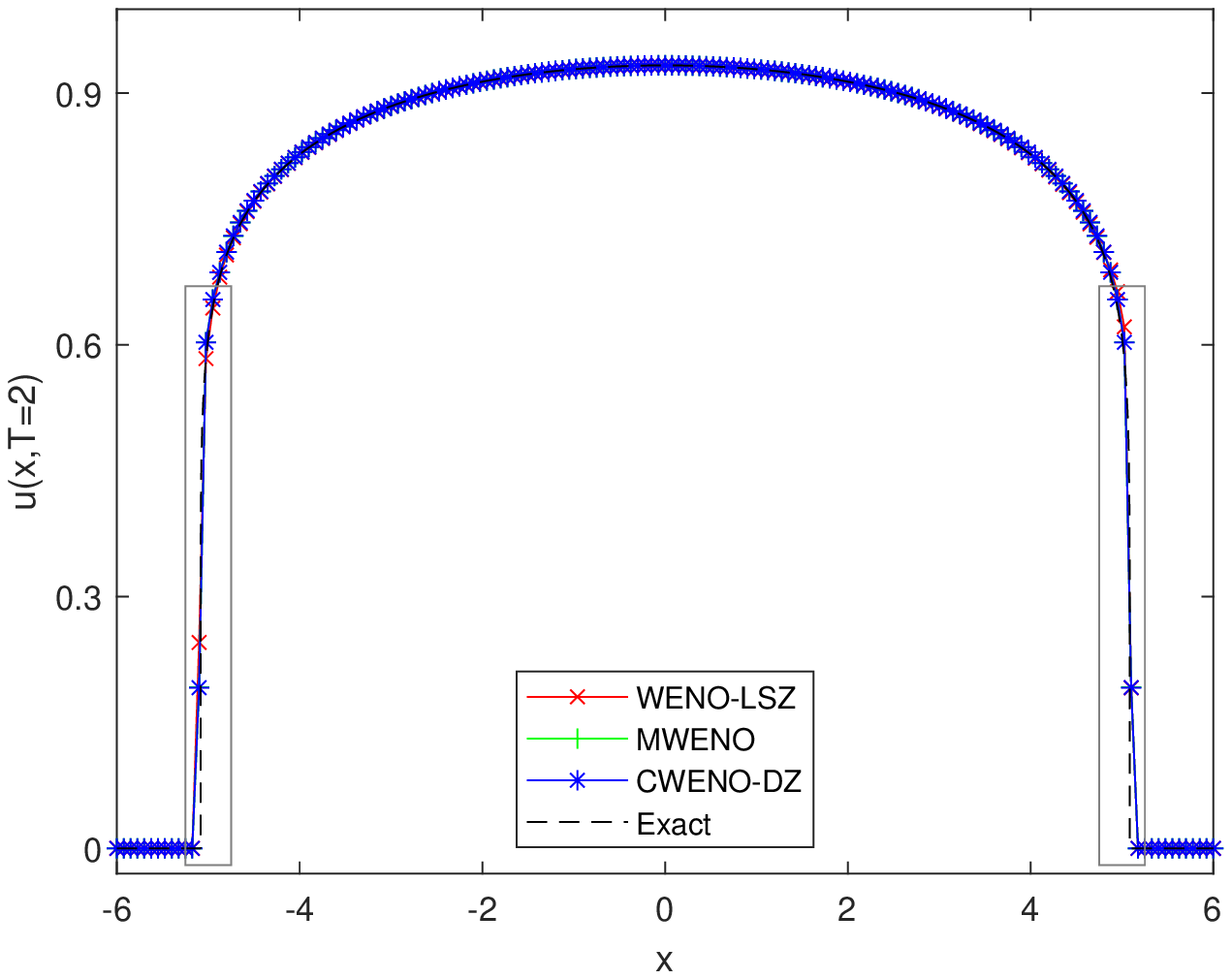}
\includegraphics[width=0.32\textwidth]{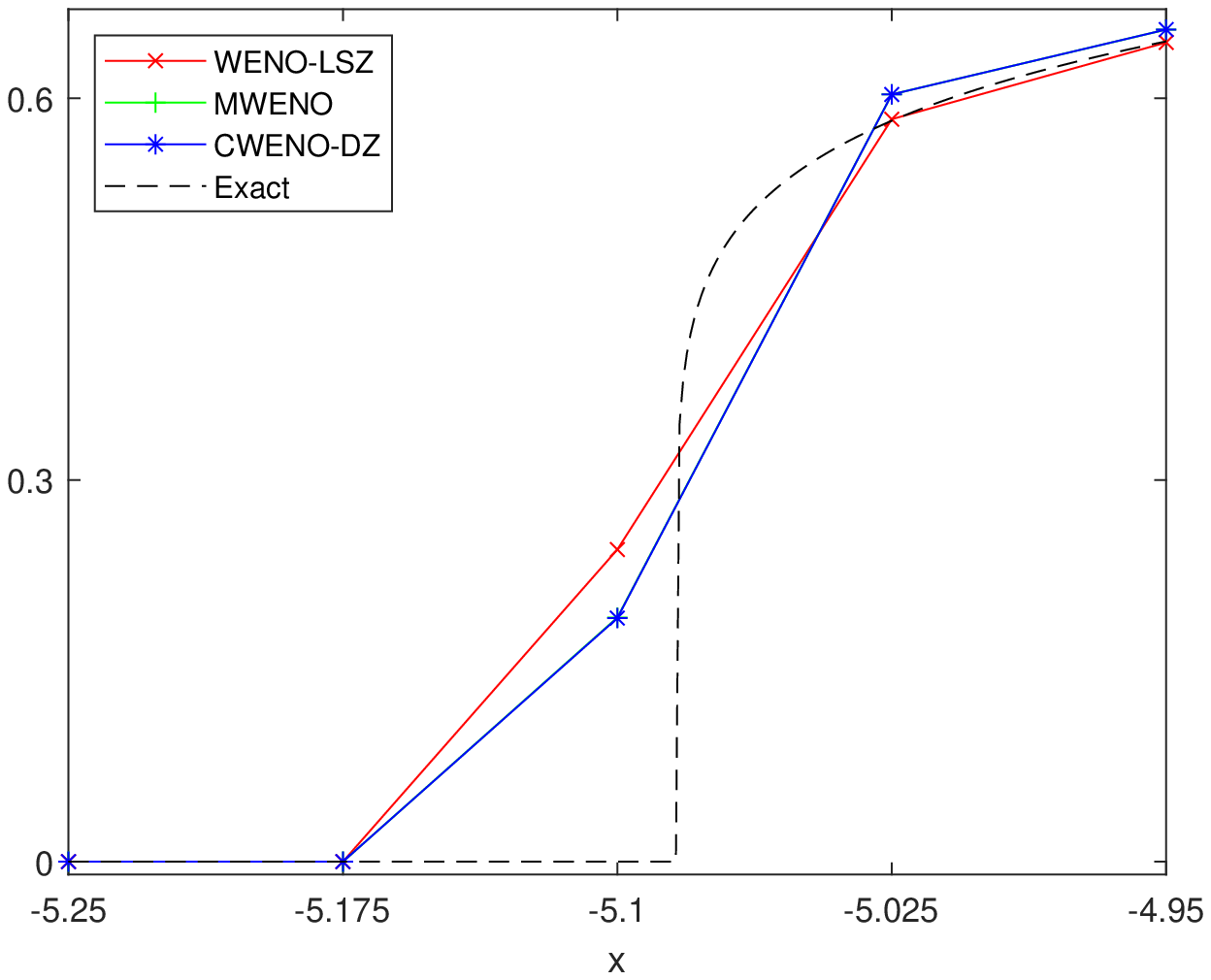}
\includegraphics[width=0.32\textwidth]{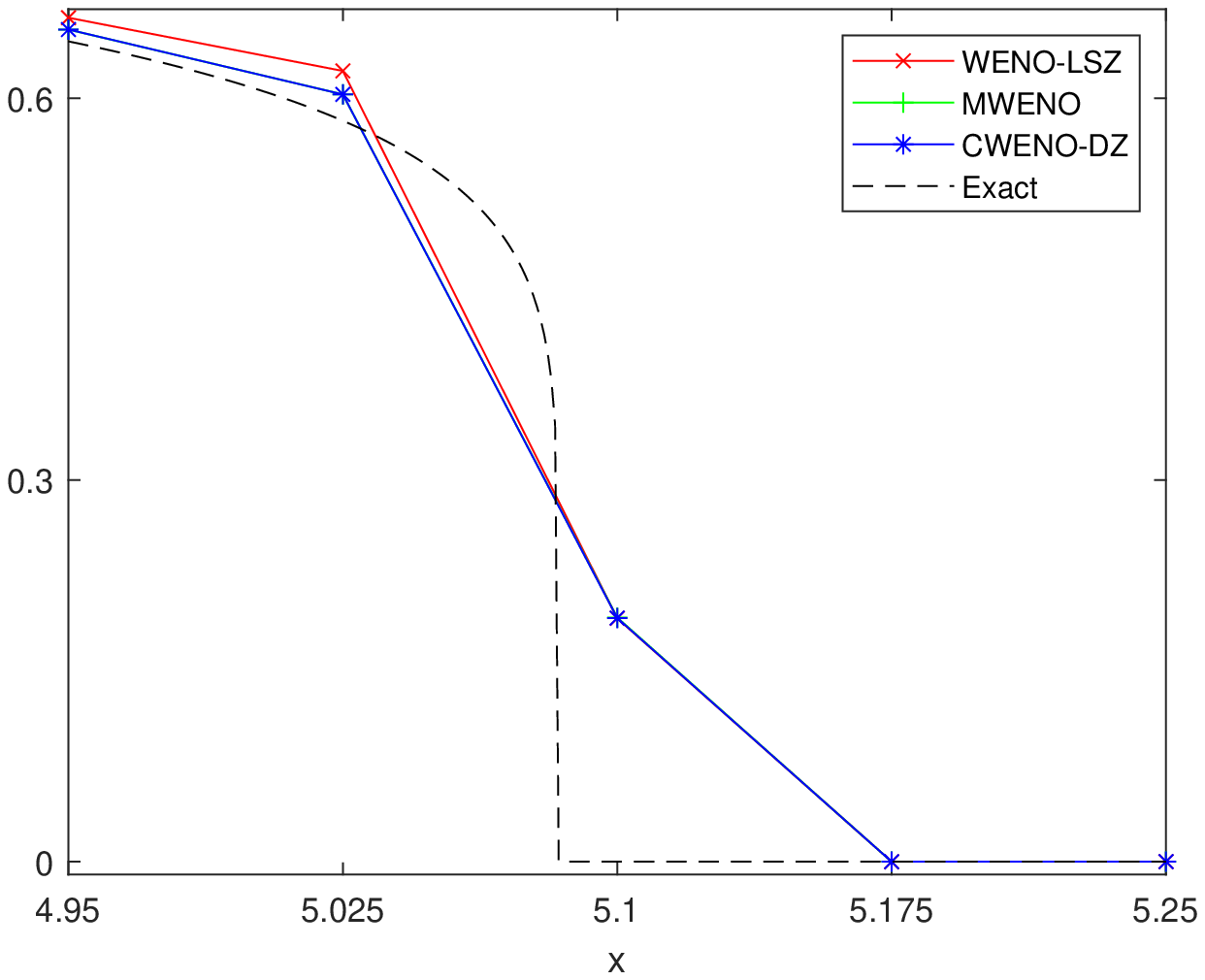}
\captionof{figure}{Barenblatt solution profiles for Example \ref{ex:Barenblatt} with $m=9$ at $T=2$ (left), close-up view of the solutions in the boxes on the left/right (middle/right) computed by WENO-LSZ (red), MWENO (green) and CWENO-DZ (blue) with $N=160$. The dashed black lines are the exact solution.}
\label{fig:Barenblatt_m9}
\end{figure}

\begin{example} \label{ex:Barenblatt} 
Consider the PME \eqref{eq:pme}. 
If the initial condition is set as the Dirac delta, the Barenblatt solution $B_m(x,t)$ \cite{Barenblatt, ZeldovichKompaneetz}, representing the heat release from a point source, takes the explicit formula 
\begin{equation} \label{eq:Barenblatt}
 B_m(x,t) = t^{-q} \left[ \left( 1 - \frac{q(m-1)}{2m} \frac{x^2}{t^{2q}} \right)_+ \right]^{1/(m-1)},~~m>1
\end{equation}
where $s_+ = \max (s,0)$ and $q = (m + 1)^{-1}$. 
For $t > 0$, the solution has a compact support $[ -a_m(t), a_m(t) ]$, where
$$
   \alpha_m(t) = \sqrt{\frac{2m}{k(m-1)}}~t^k,
$$
and the interfaces $|x| = a_m(t)$ move outward at a finite speed.
Moreover, the larger the value of $p$, the sharper the interfaces that separate the compact support and the zero solution.

We simulate the Barenblatt solution \eqref{eq:Barenblatt} of the PME \eqref{eq:pme} with the initial condition as the Barenblatt solution at $t=1$, $u(x,0) = B_m(x,1)$, and the boundary conditions $u(\pm 6,t)=0$ for $t>0$. 
The final time is $T = 2$ and the time step is $\Delta t = \cfl \cdot \Delta x^2/m$. 
We take $N=160$ and plot the numerical solutions at the final time for $m = 5, \, 7$ and $9$, in Figures \ref{fig:Barenblatt_m5}, \ref{fig:Barenblatt_m7} and \ref{fig:Barenblatt_m9}, respectively.
We can see that the solution of the proposed CWENO-DZ almost overlaps the one of MWENO but both give more accurate solution profiles around the interfaces than WENO-LSZ.
This is also demonstrated by Table \ref{tab:Barenblatt}, which provides the $L_1, L_2$ and $L_{\infty}$ errors for the WENO-LSZ, MWENO and CWENO-DZ schemes.
\end{example}

\begin{table}[h!]
\centering
\caption{$L_1,~L_2$ and $L_\infty$ errors for Example \ref{ex:Barenblatt}.}      
\begin{tabular}{ccccc} 
\hline  
m & error & WENO-LSZ & MWENO & CWENO-DZ \\ 
\hline 
\multirow{3}{*}{5} & $L_1$      & 2.81E-3 & 1.47E-3 & 1.45E-3 \\  
                   & $L_2$      & 1.82E-2 & 1.15E-2 & 1.14E-2 \\  
                   & $L_\infty$ & 1.77E-1 & 1.03E-1 & 1.02E-1 \\
\hline 
\multirow{3}{*}{7} & $L_1$      & 2.77E-3 & 1.39E-3 & 1.37E-3 \\  
                   & $L_2$      & 1.74E-2 & 1.05E-2 & 1.04E-2 \\
                   & $L_\infty$ & 1.73E-1 & 9.38E-2 & 9.31E-2 \\
\hline
\multirow{3}{*}{9} & $L_1$      & 3.25E-3 & 3.19E-3 & 3.19E-3 \\  
                   & $L_2$      & 2.48E-2 & 2.16E-2 & 2.15E-2 \\  
                   & $L_\infty$ & 2.45E-1 & 1.92E-1 & 1.91E-1 \\
\hline
\end{tabular}
\label{tab:Barenblatt}
\end{table}

\begin{figure}[h!]
\centering
\includegraphics[width=0.32\textwidth]{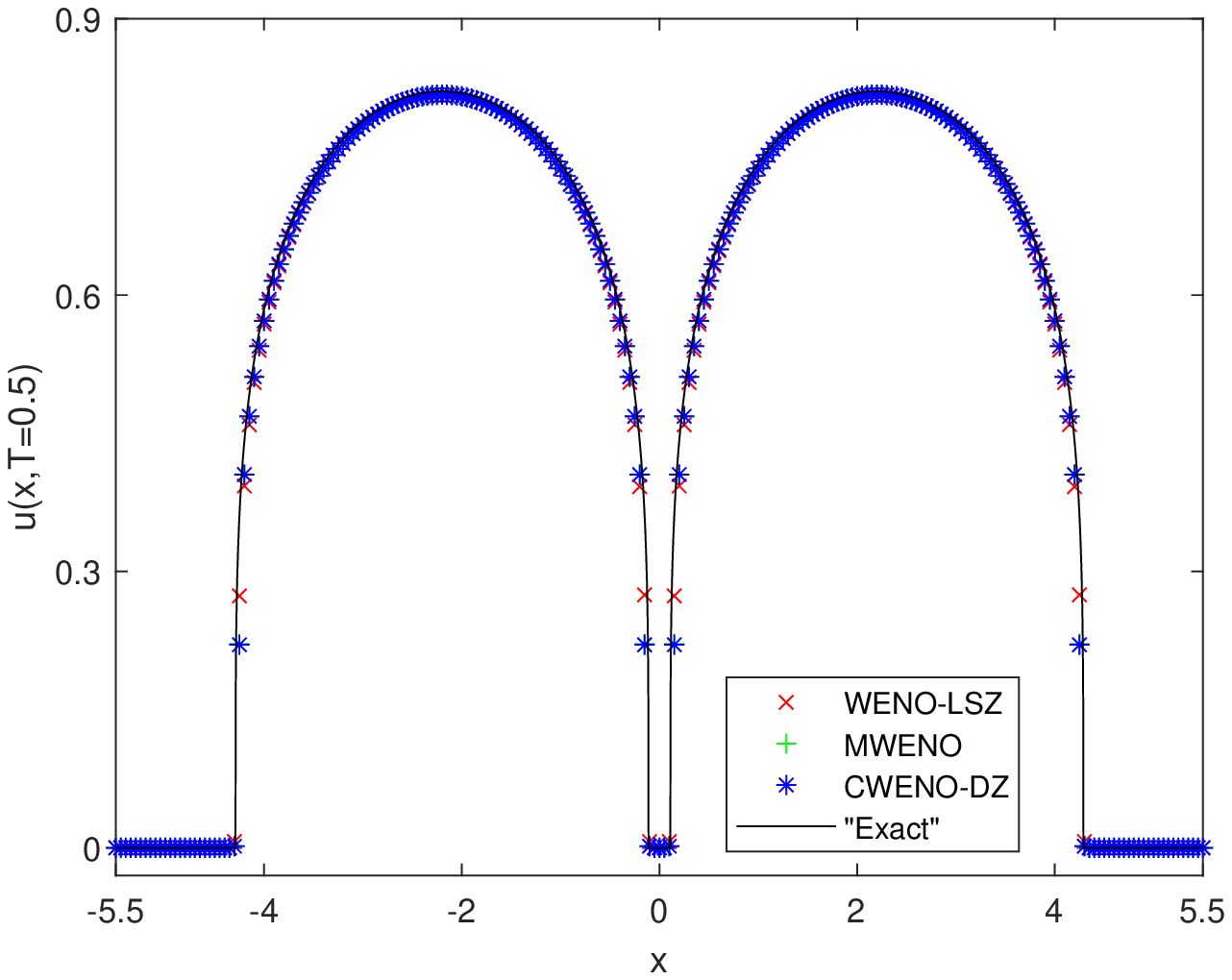}
\includegraphics[width=0.32\textwidth]{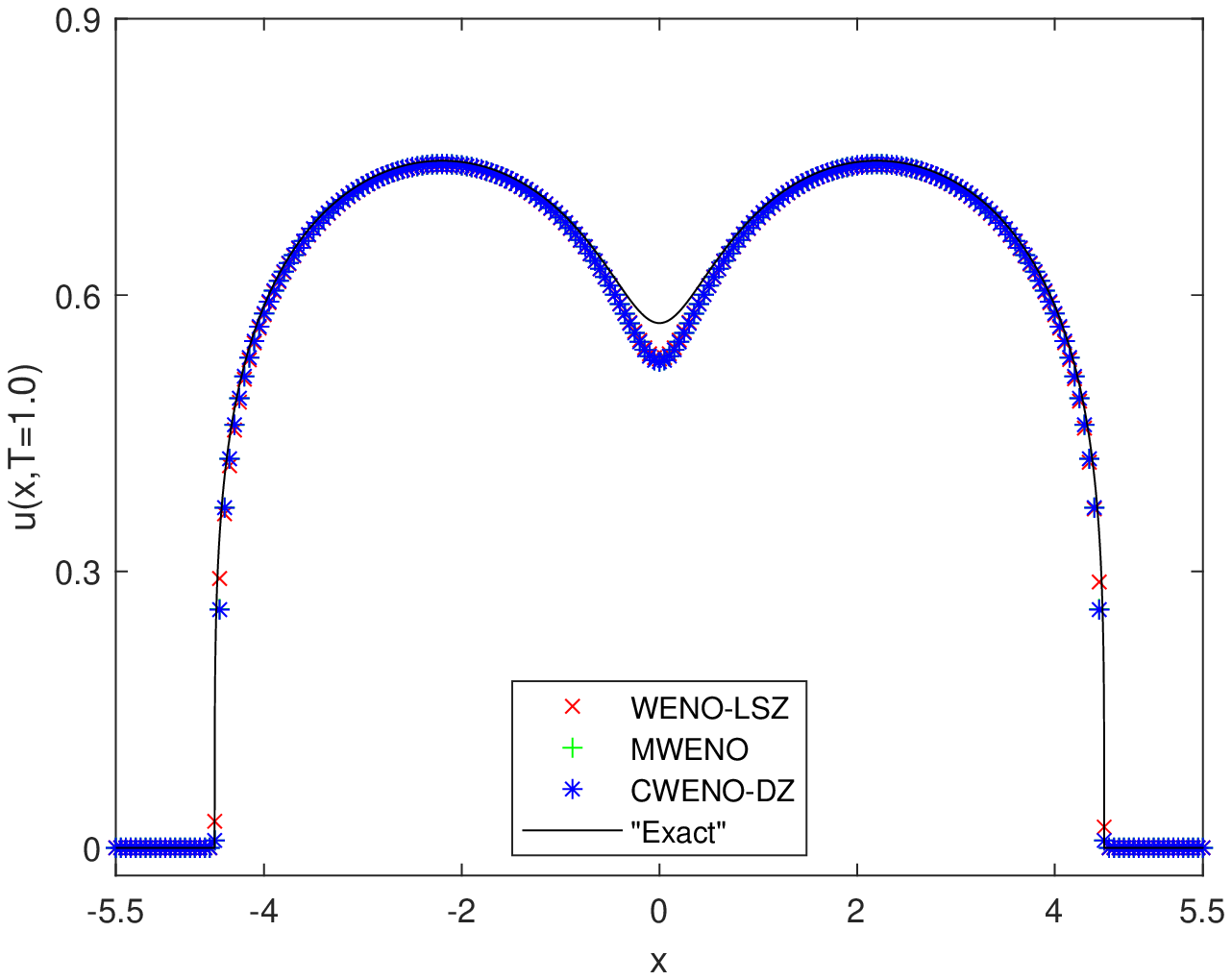}
\includegraphics[width=0.32\textwidth]{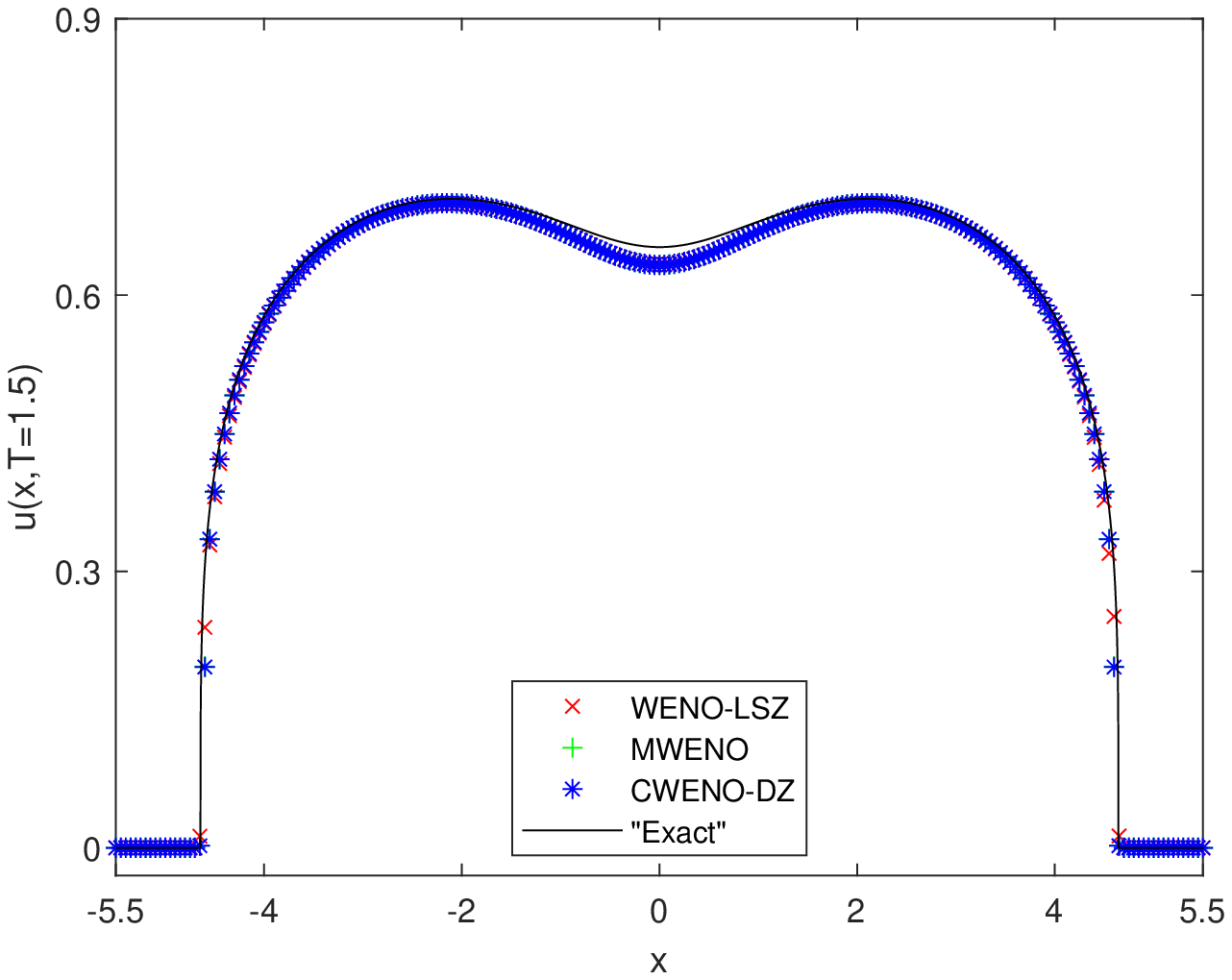}
\captionof{figure}{Solution profiles for PME \eqref{eq:pme} ($m=5$) with the initial condition \eqref{eq:two_box_init_sh} at $t = 0.5$ (left), $1.0$ (middle) and $1.5$ (right) approximated by WENO-LSZ (red), MWENO (green) and CWENO-DZ (blue) with $N = 220$.
The black lines are generated by MWENO with $N = 11000$.}
\label{fig:two_box_sh}
\end{figure}

\begin{figure}[h!]
\centering
\includegraphics[width=0.32\textwidth]{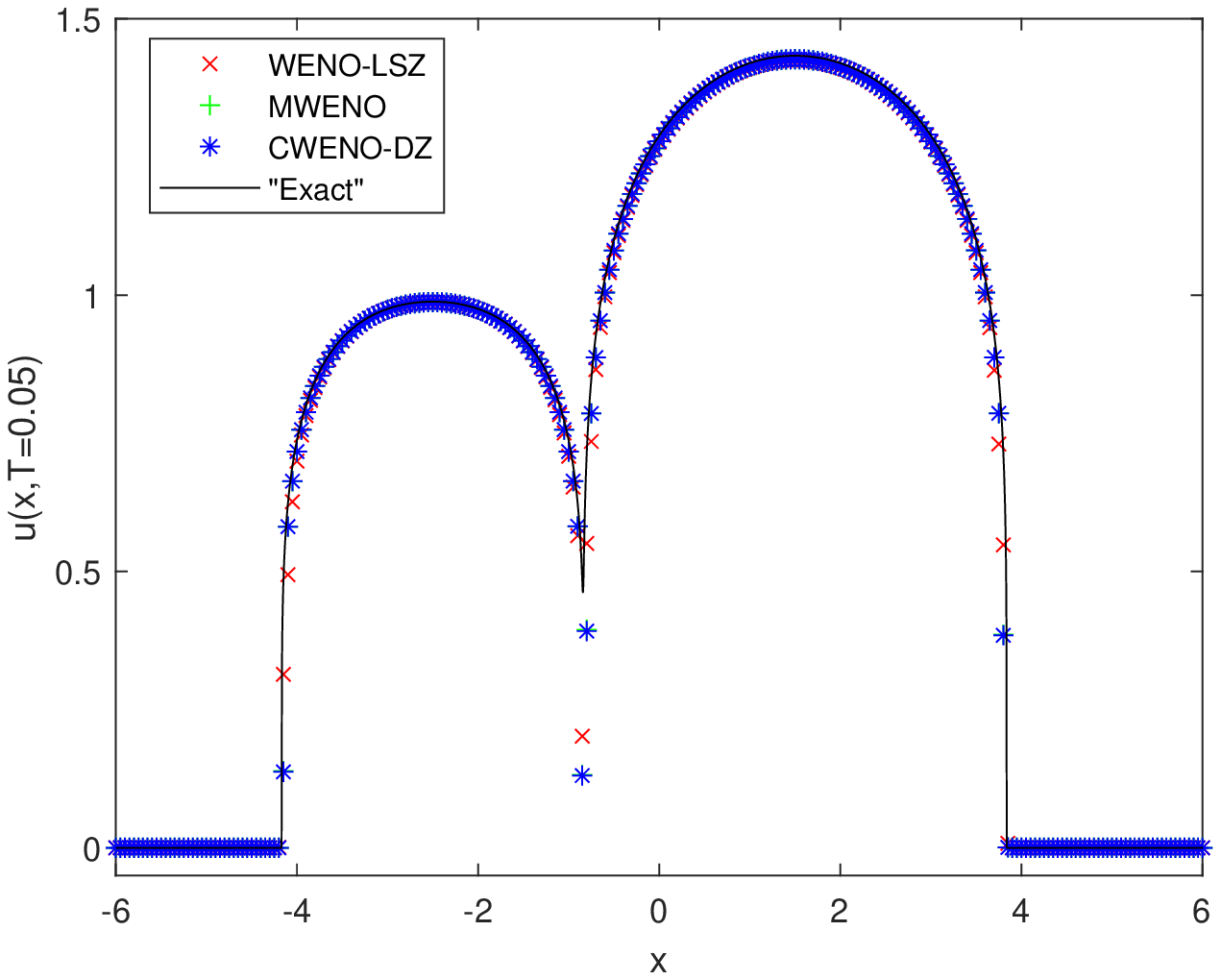}
\includegraphics[width=0.32\textwidth]{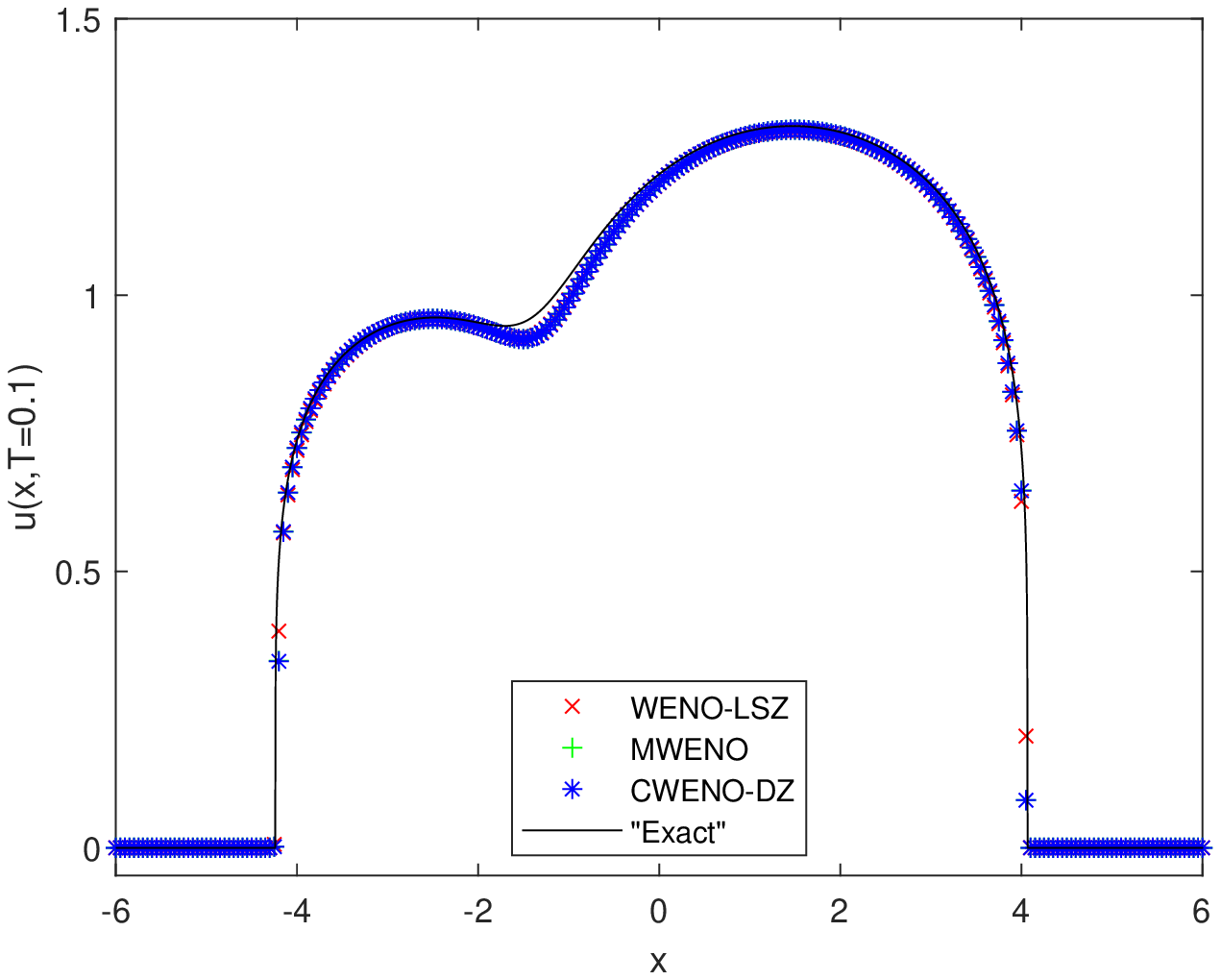}
\includegraphics[width=0.32\textwidth]{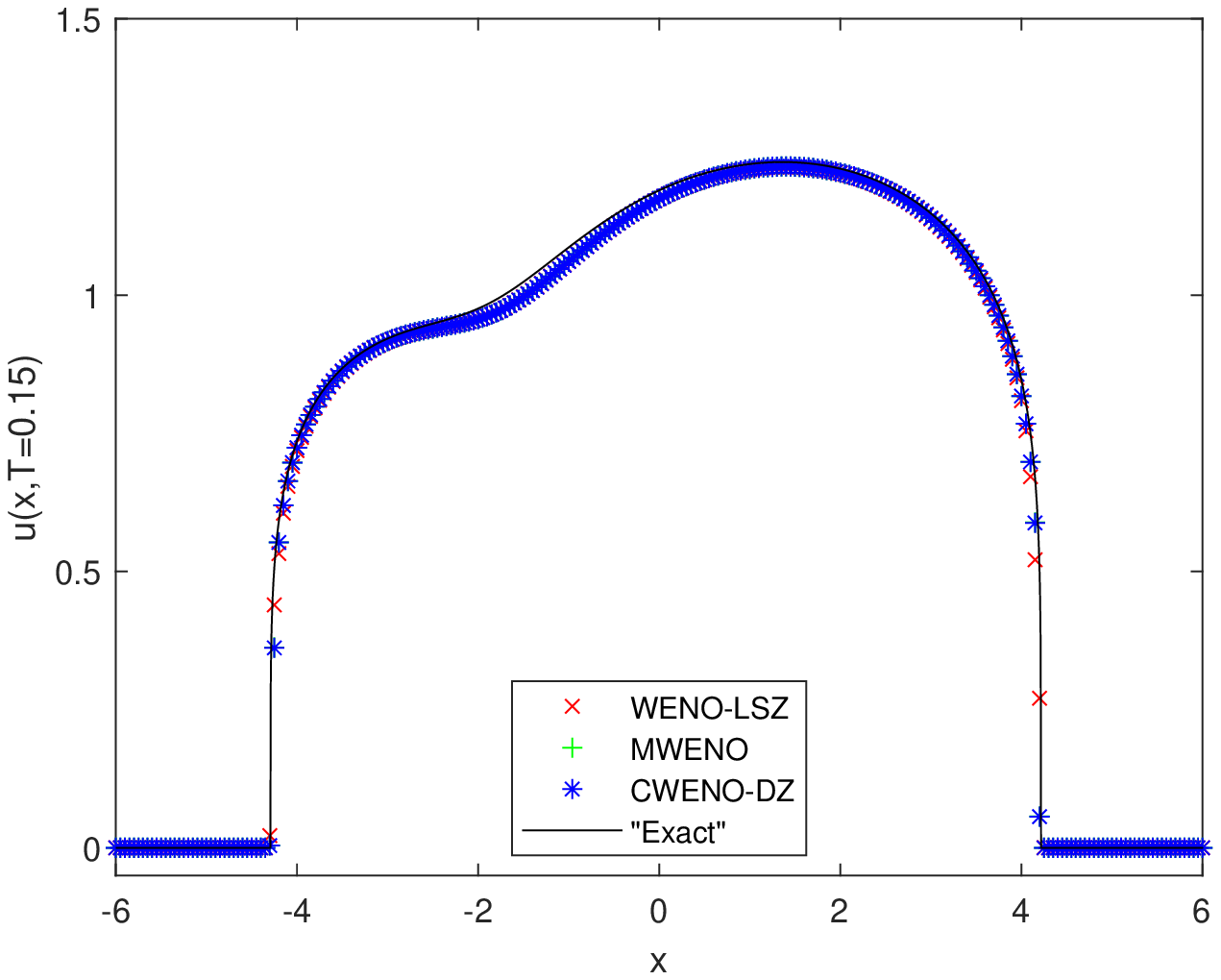}
\captionof{figure}{Solution profiles for PME \eqref{eq:pme} ($m=6$) with the initial condition \eqref{eq:two_box_init_dh} at $t = 0.05$ (left), $0.1$ (middle) and $0.15$ (right) approximated by WENO-LSZ (red), MWENO (green) and CWENO-DZ (blue) with $N = 240$.
The black lines are generated by MWENO with $N = 6000$.}
\label{fig:two_box_dh}
\end{figure}

\begin{example} \label{ex:two_box}
We continue to consider the PME \eqref{eq:pme}, where the shape of the initial condition is two separate boxes.
If the solution $u$ represents the temperature, the PME models the variations in temperature when two hot spots are situated in the domain.

We first consider the PME with $m=5$, where the initial condition is given by
\begin{equation} \label{eq:two_box_init_sh}
 u(x,0) = \begin{cases}
          1, & x \in (-3.7,-0.7) \cup (0.7, 3.7), \\
          0, & \mbox{otherwise},
         \end{cases}
\end{equation}
in which the two boxes have the same height, and the boundary conditions are $u(\pm 5.5,t)=0$ for $t>0$. 
We divide the computational domain $[-5.5,~5.5]$ into $N = 220$ uniform cells.
The final time is $T = 1.5$ and the time step is $\Delta t = \cfl \cdot \Delta x^2/m$.
We present the numerical solutions at $t = 0.5, \; 1.0, \; 1.5$, as shown in Figure \ref{fig:two_box_sh}.
The numerical solution, computed by MWENO with a high resolution of $N = 11000$ points, will be referred to as the ``exact'' solution.

Now we consider the PME with $m=6$. 
The initial condition in this case is
\begin{equation} \label{eq:two_box_init_dh}
 u(x,0) = \begin{cases}
          1, & -4 < x < -1, \\
          2, & 0 < x < 3, \\
          0, & \mbox{otherwise},
         \end{cases}
\end{equation}
and the boundary conditions are $u(\pm 6,t) = 0$ for $t>0$. 
We select $N = 240$ for the computational domain $[-6,~6]$.
Figure \ref{fig:two_box_dh} shows the approximate results obtained when solving PME up to the final time $T = 0.15$ with the time step $\Delta t = \cfl \cdot \left. \Delta x^2 \right/ \left( m 2^{m-1} \right)$.
We still take the solution computed by MWENO with $N = 6000$ points as the ``exact'' solution.

As seen in Figures \ref{fig:two_box_sh} and \ref{fig:two_box_dh}, all schemes are able to capture the sharp interfaces, and MWENO and CWENO-DZ yield very similar solution profiles.
\end{example}

Next, we solve the one-dimensional scalar convection-diffusion equation of the form
$$
   u_t + f(u)_x = g(u)_{xx}.
$$
For the convection term, the fifth-order finite difference Lax–Friedrichs flux splitting WENO scheme, WENO-JS \cite{JiangShu, Shu}, is employed as we want to see how those WENO schemes for the diffusion term affect the numerical solutions.
The numerical solution, computed by WENO-M \cite{Henrick} and MWENO for the respective convection and diffusion terms with a high resolution, will be referred to as the ``exact'' solution.

\begin{figure}[h]
\centering
\includegraphics[width=0.32\textwidth]{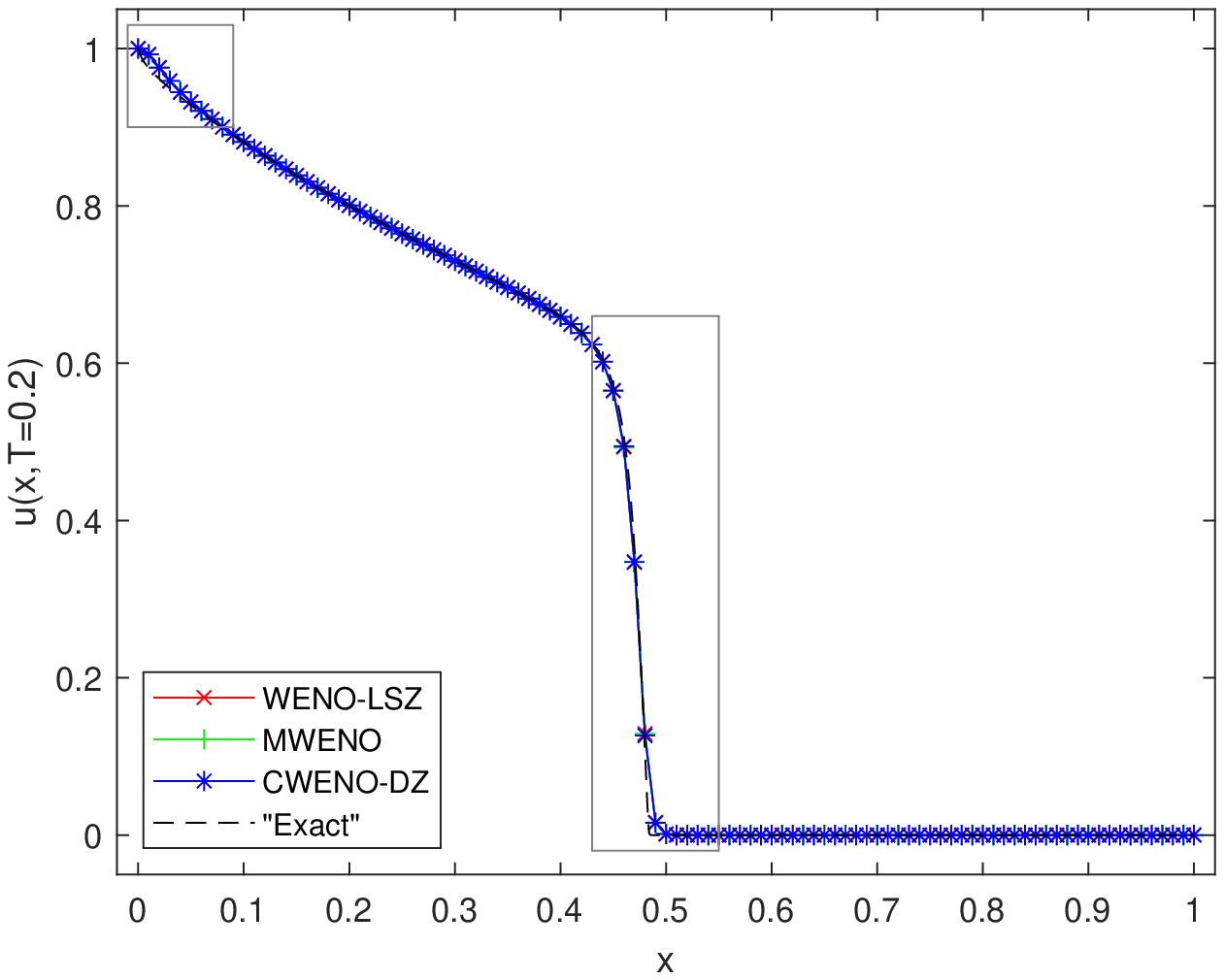}
\includegraphics[width=0.32\textwidth]{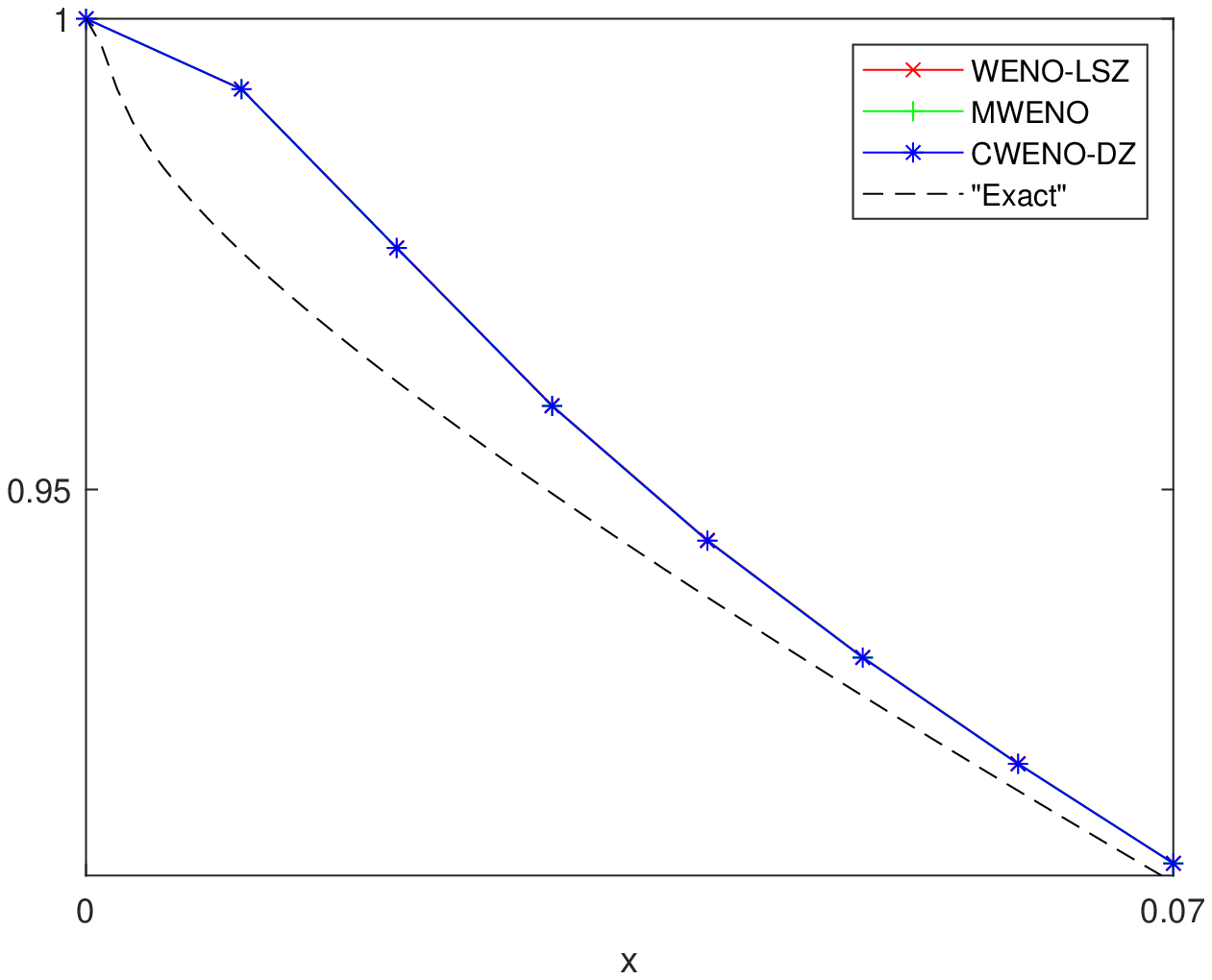}
\includegraphics[width=0.32\textwidth]{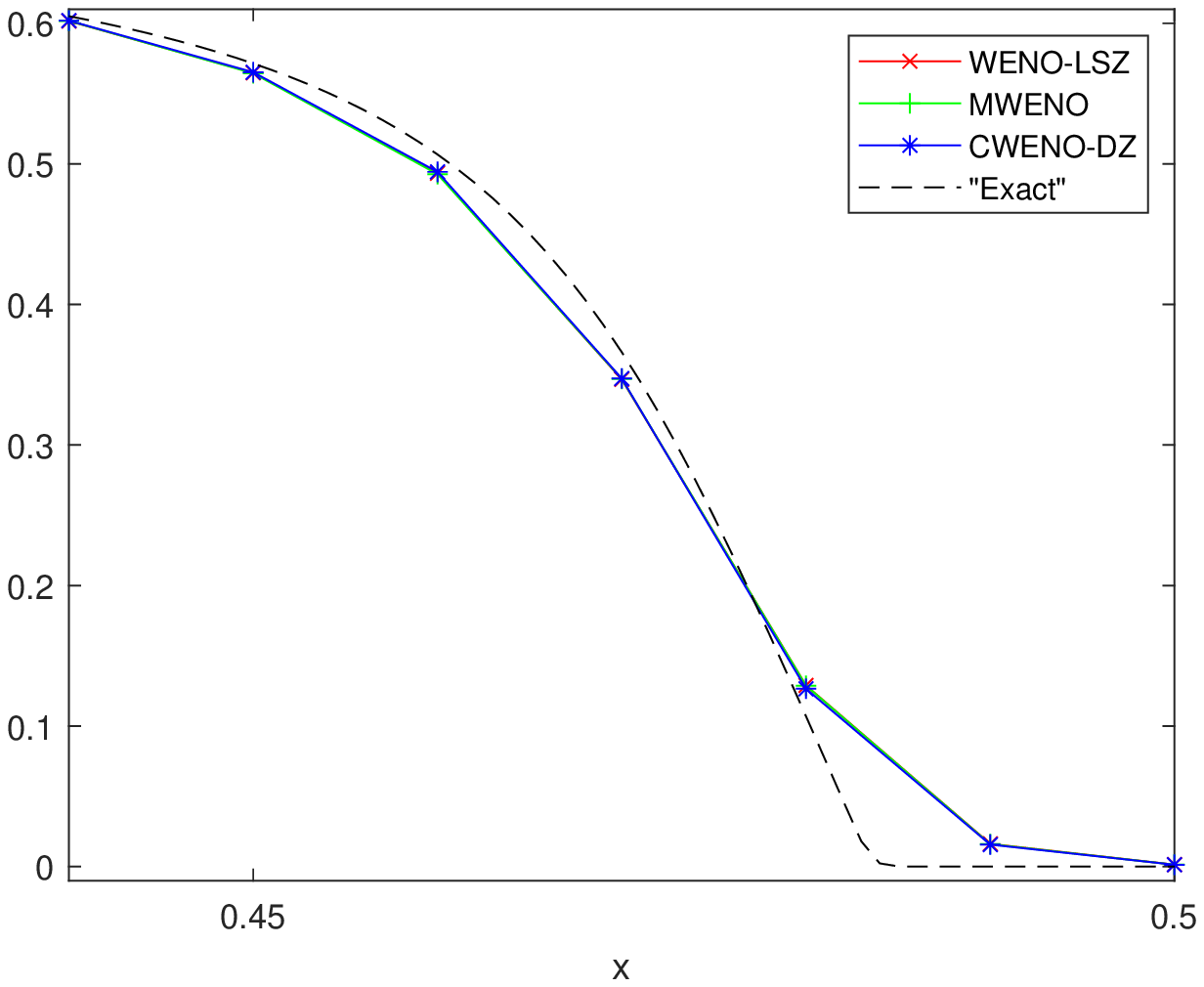}
\captionof{figure}{Solution profiles for Buckley-Leverett equation in Example \ref{ex:Buckley_Leverett} at $T = 0.2$ (left), close-up view of the solutions in the boxes on the left/right (middle/right) computed by WENO-LSZ (red), MWENO (green) and CWENO-DZ (blue) with $N = 100$. 
The dashed black lines are generated by WENO-M and MWENO with $N = 1000$.}
\label{fig:BLE}
\end{figure}

\begin{example} \label{ex:Buckley_Leverett}
The Buckley-Leverett equation \cite{Buckley} is of the form
\begin{equation} \label{eq:Buck}
 u_t + f(u)_x = \epsilon \left( \nu(u) u_x \right)_x,~~\epsilon \nu(u) \geqslant 0,
\end{equation}
which is a prototype model for oil reservoir simulation.
This is an example of degenerate parabolic equations since $\nu(u)$ vanishes at some values of $u$.
Following \cite{Kurganov}, the convection flux $f(u)$ is of the s-shaped form  
\begin{equation} \label{eq:nongravitaional_flux}
 f(u) = \frac{u^2}{u^2+(1-u)^2},
\end{equation}   
$\epsilon=0.01$, and
\begin{equation} \label{eq:diffusion_coeff_BLE}
 \nu(u) = \begin{cases}
          4u(1-u), & 0 \leqslant u \leqslant 1, \\
          0,       & \mbox{otherwise}.
          \end{cases}
\end{equation} 
The diffusion term $\epsilon \left( \nu(u) u_x \right)_x$ can be written in the form of $g(u)_{xx}$, where
$$
   g(u) = \begin{cases}
          0, & u < 0, \\
          \epsilon \left( - \frac{4}{3} u^3 + 2 u^2 \right), & 0 \leqslant u \leqslant 1, \\
          \frac{2}{3} \epsilon, & u > 1.
         \end{cases}
$$
The initial condition is given by
$$
   u(x,0) = \begin{cases}
            1-3x, & 0 \leqslant x \leqslant 1/3, \\ 
            0,    & 1/3 < x \leqslant 1. 
           \end{cases}
$$
and the Dirichlet boundary condition is $u(0,t)=1$.
The computational domain $[0,~1]$ is divided into $N = 100$ uniform cells and the time step is $\Delta t = \cfl \cdot \Delta x^2$.
The numerical solution computed by CWENO-DZ at $T=0.2$ is very close to those by WENO-LSZ and MWENO. 
This results in the overlapping in Figure \ref{fig:BLE}.
\end{example}

\begin{figure}[h]
\centering
\includegraphics[width=0.32\textwidth]{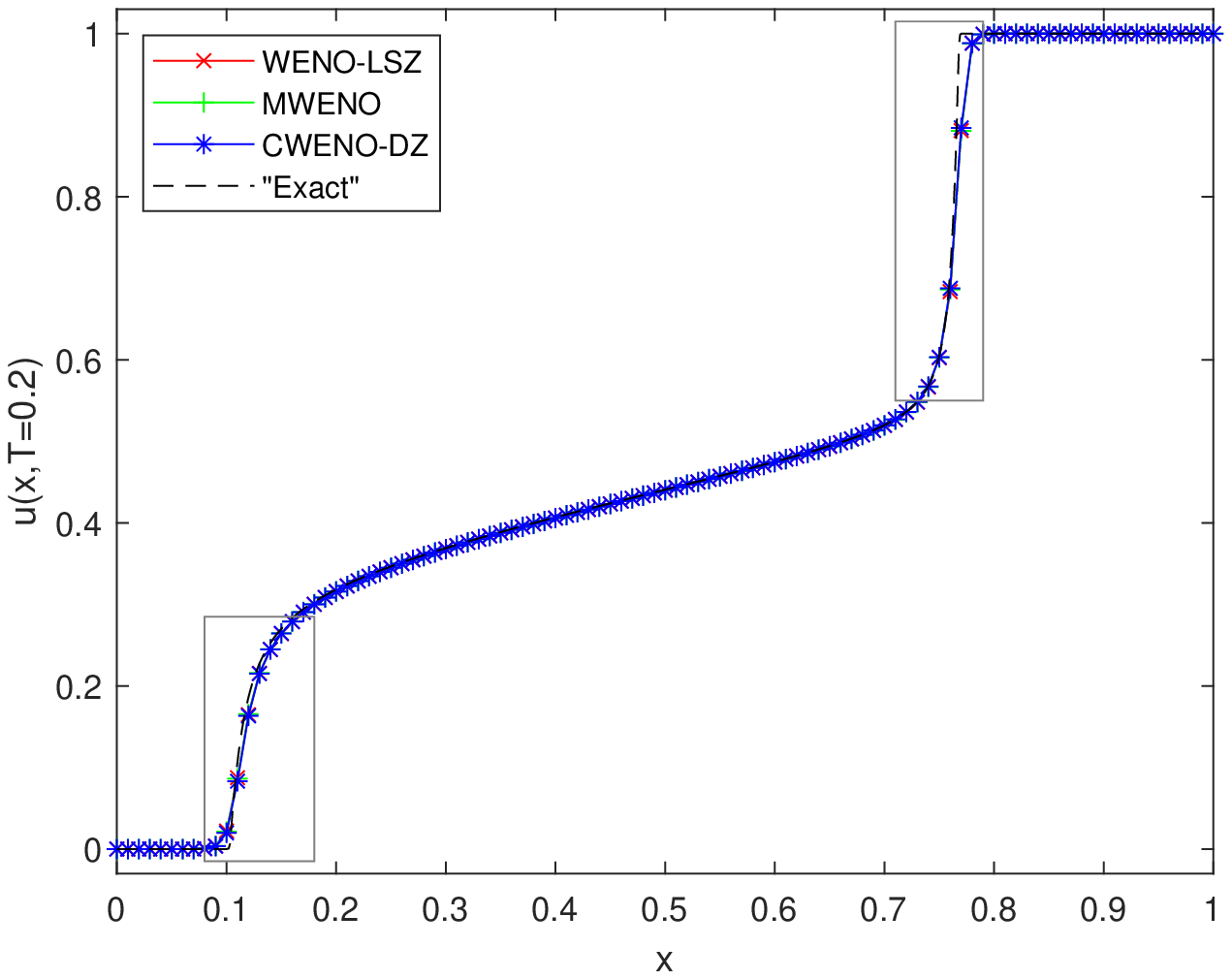}
\includegraphics[width=0.32\textwidth]{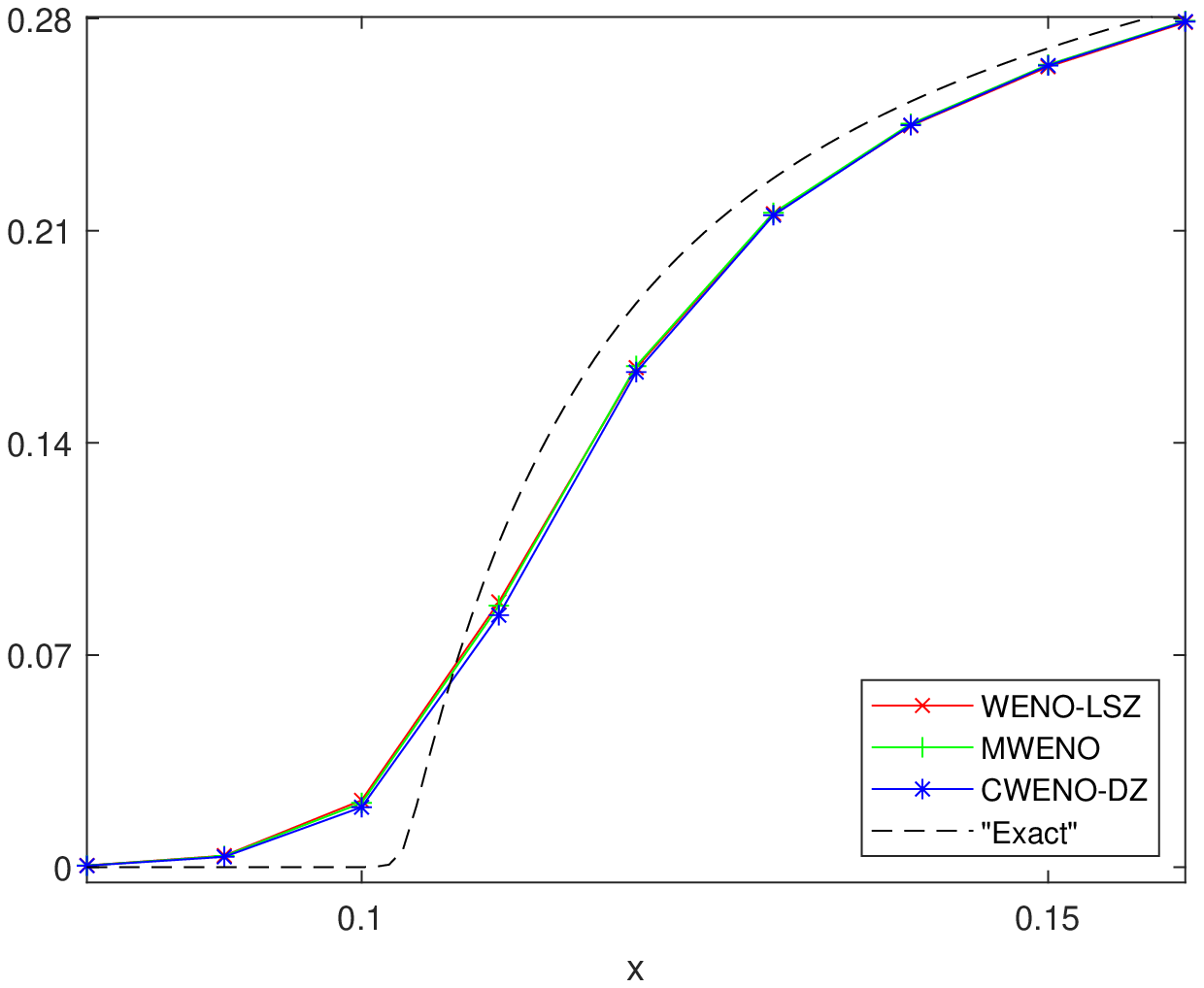}
\includegraphics[width=0.32\textwidth]{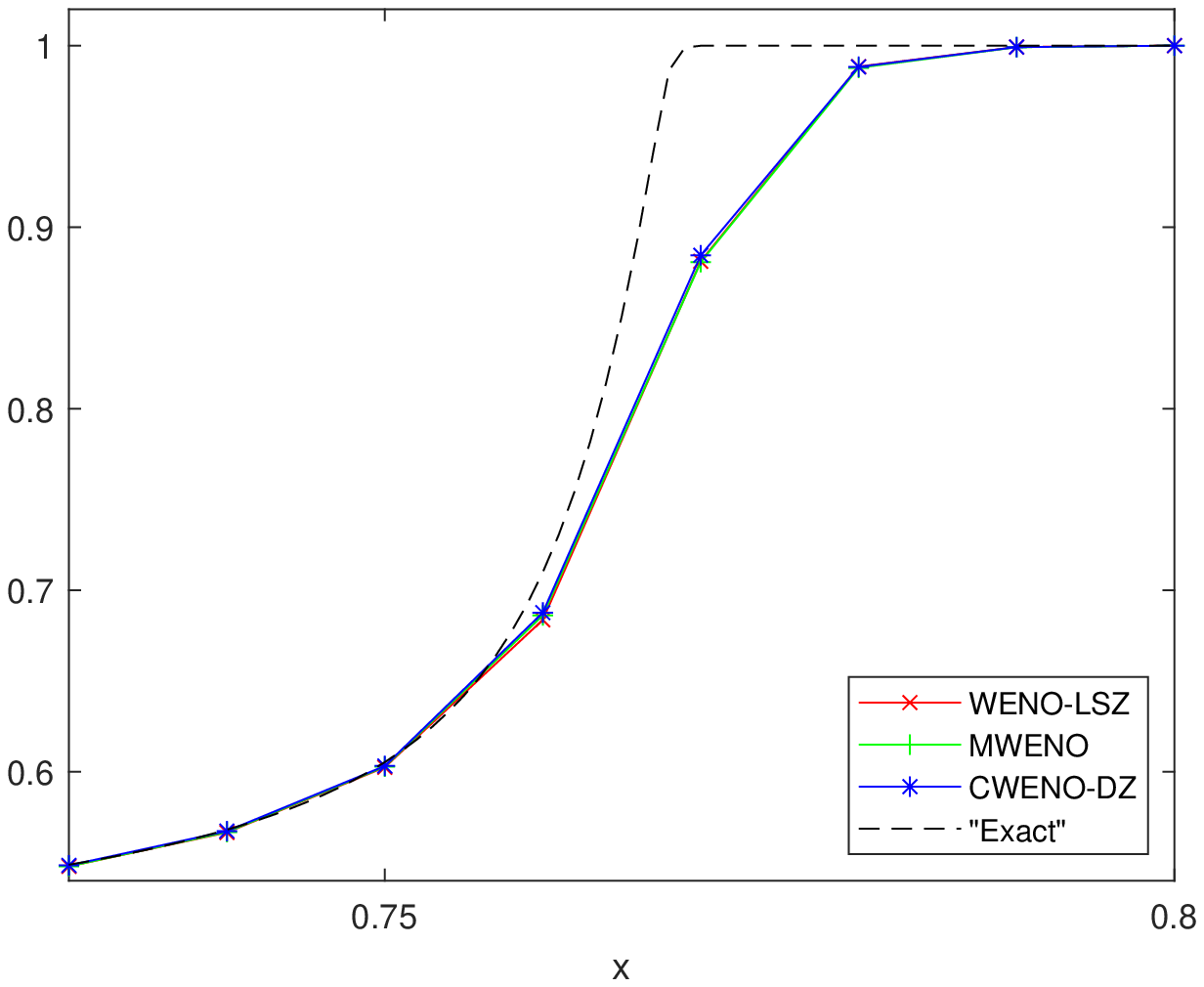}
\captionof{figure}{Solution profiles for Buckley-Leverett equation \eqref{eq:Buck} with gravitation in Example \ref{ex:Buckley_Leverett_gravitation} at $T=0.2$ (left), close-up view of the solutions in the boxes on the left/right (middle/right) computed by WENO-LSZ (red), MWENO (green) and CWENO-DZ (blue) with $N = 100$. The dashed black lines are generated by WENO-M and MWENO with $N = 1000$.}
\label{fig:BLE_grav}
\end{figure}

\begin{figure}[h]
\centering
\includegraphics[width=0.32\textwidth]{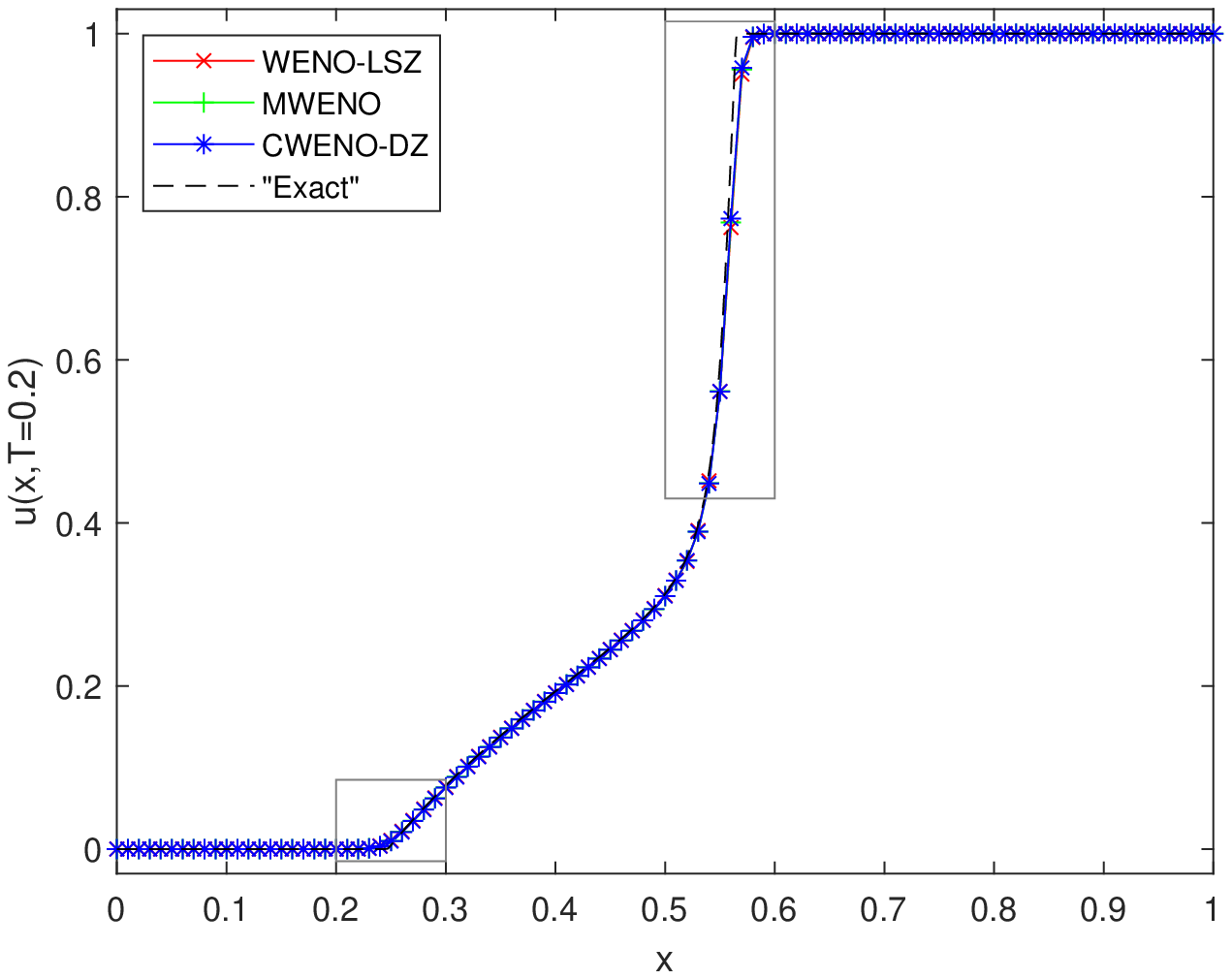}
\includegraphics[width=0.32\textwidth]{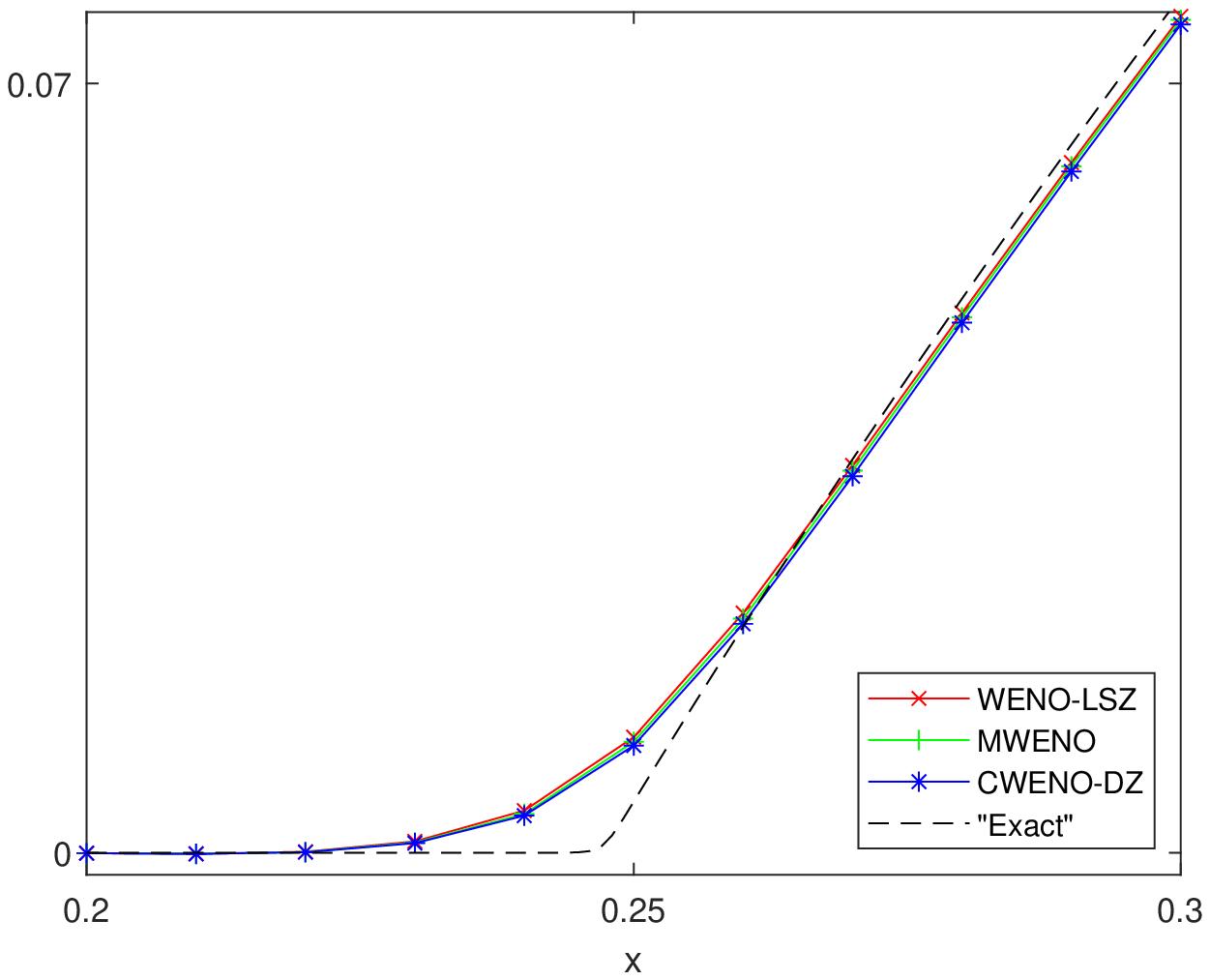}
\includegraphics[width=0.32\textwidth]{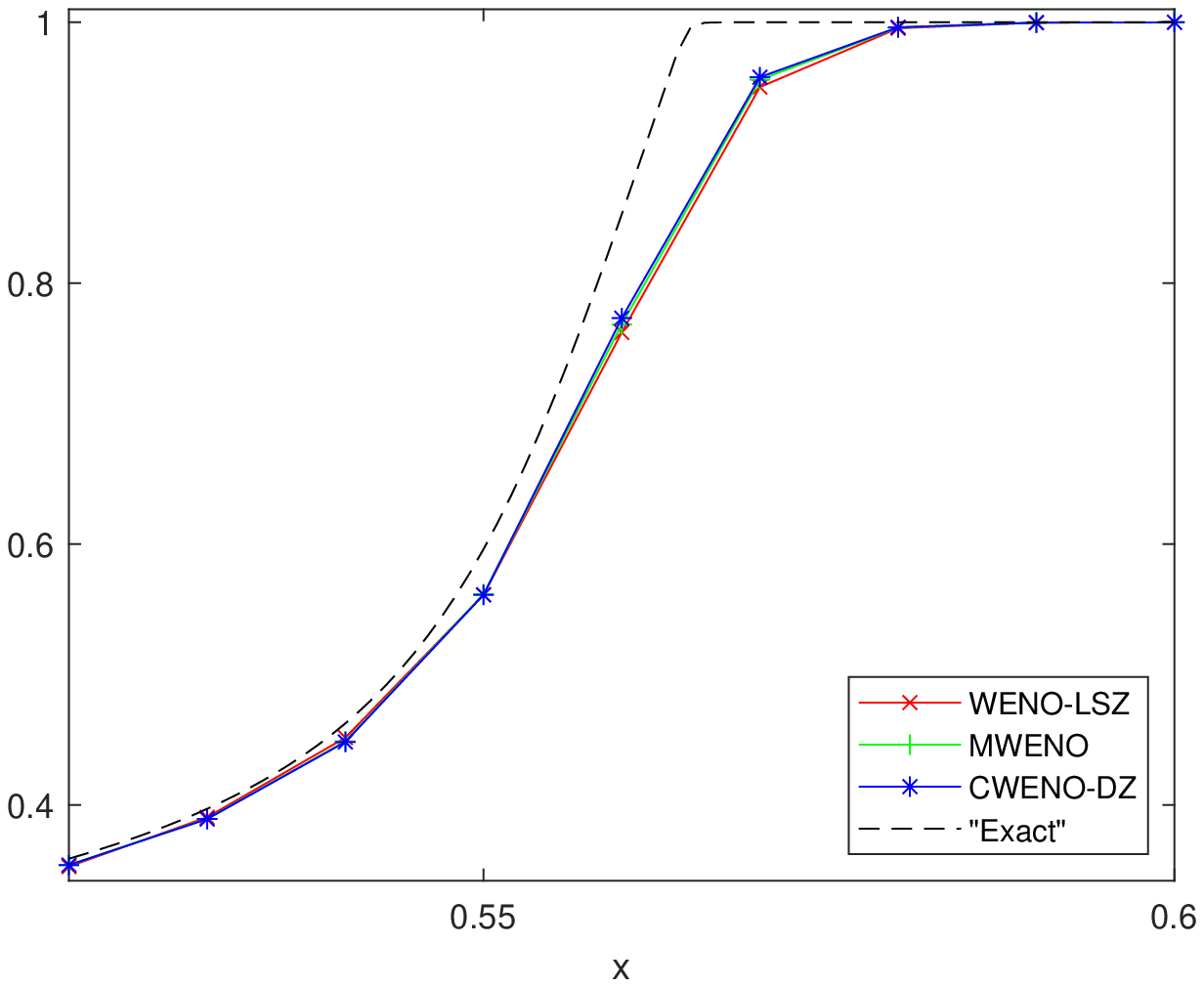}
\captionof{figure}{Solution profiles for Buckley-Leverett equation \eqref{eq:Buck} without gravitation in Example \ref{ex:Buckley_Leverett_gravitation} at $T=0.2$ (left), close-up view of the solutions in the boxes on the left/right (middle/right) computed by WENO-LSZ (red), MWENO (green) and CWENO-DZ (blue) with $N = 100$. The dashed black lines are generated by WENO-M and WENO-LSZ with $N = 1000$.}
\label{fig:BLE_nograv}
\end{figure}

\begin{example} \label{ex:Buckley_Leverett_gravitation}
We continue to consider the Buckley-Leverett equation \eqref{eq:Buck} with the same $\epsilon=0.01$ and $\nu(u)$ \eqref{eq:diffusion_coeff_BLE} as in Example \ref{ex:Buckley_Leverett}. 
The flux function $f(u)$ with gravitational effects is
\begin{equation} \label{eq:gravitaional_flux}
 f(u)=\frac{u^2}{u^2+(1-u)^2}(1-5(1-u)^2),
\end{equation}
where the sign of $f'(u)$ changes in $[0,~1]$.
The Riemann initial condition is
$$
   u(x,0) = \begin{cases}
             0, & 0 \leqslant x < 1 - 1/\sqrt{2}, \\
             1, & 1 - 1/\sqrt{2} \leqslant x \leqslant 1.
            \end{cases}
$$
We divide the computational domain $[0,~1]$ into $N=100$ uniform cells and the time step is $\Delta t = \cfl \cdot \Delta x^2$.
Figure \ref{fig:BLE_grav} shows the numerical solutions at $T=0.2$ for the convection flux $f(u)$ \eqref{eq:gravitaional_flux} with gravitational effects while Figure \ref{fig:BLE_nograv} presents the ones for $f(u)$ \eqref{eq:nongravitaional_flux} without gravitational effects.
In Figure \ref{fig:BLE_grav}, all WENO schemes yield comparable results, while in Figure \ref{fig:BLE_nograv}, CWENO-DZ produces the numerical solution slightly closer to the reference solution than WENO-LSZ and MWENO around the shock.
\end{example}

\begin{figure}[h]
\centering
\includegraphics[width=0.45\textwidth]{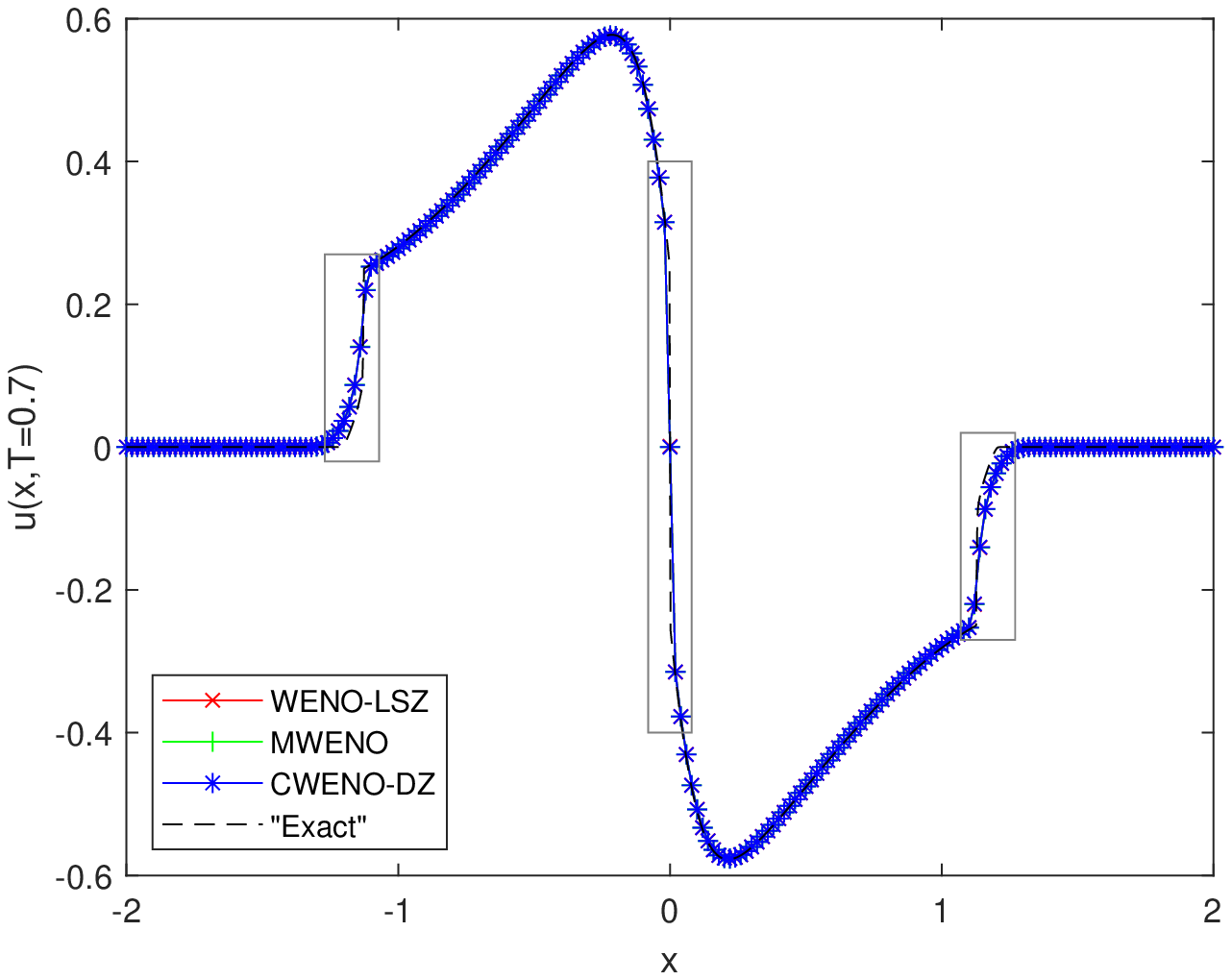}
\includegraphics[width=0.45\textwidth]{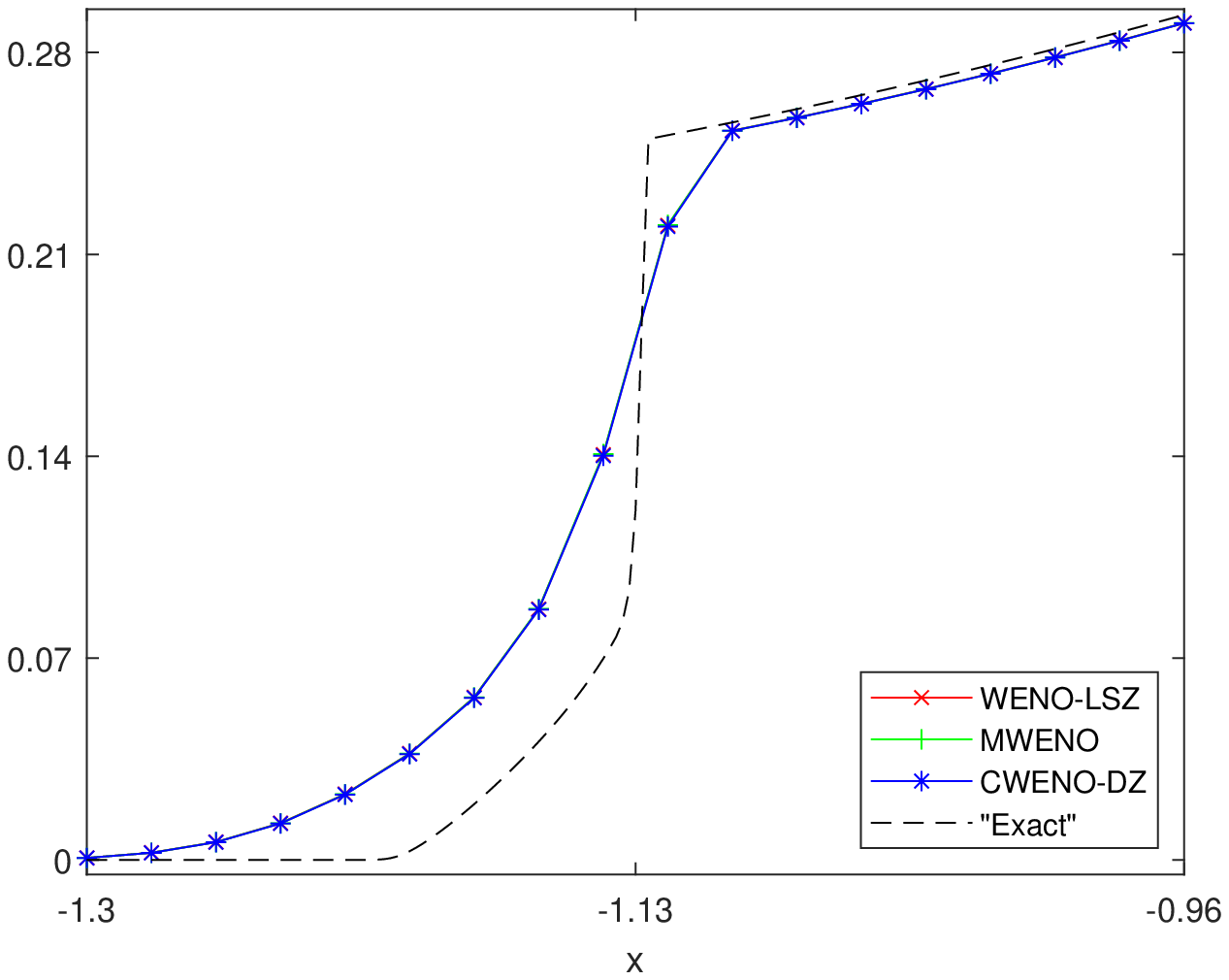}
\includegraphics[width=0.45\textwidth]{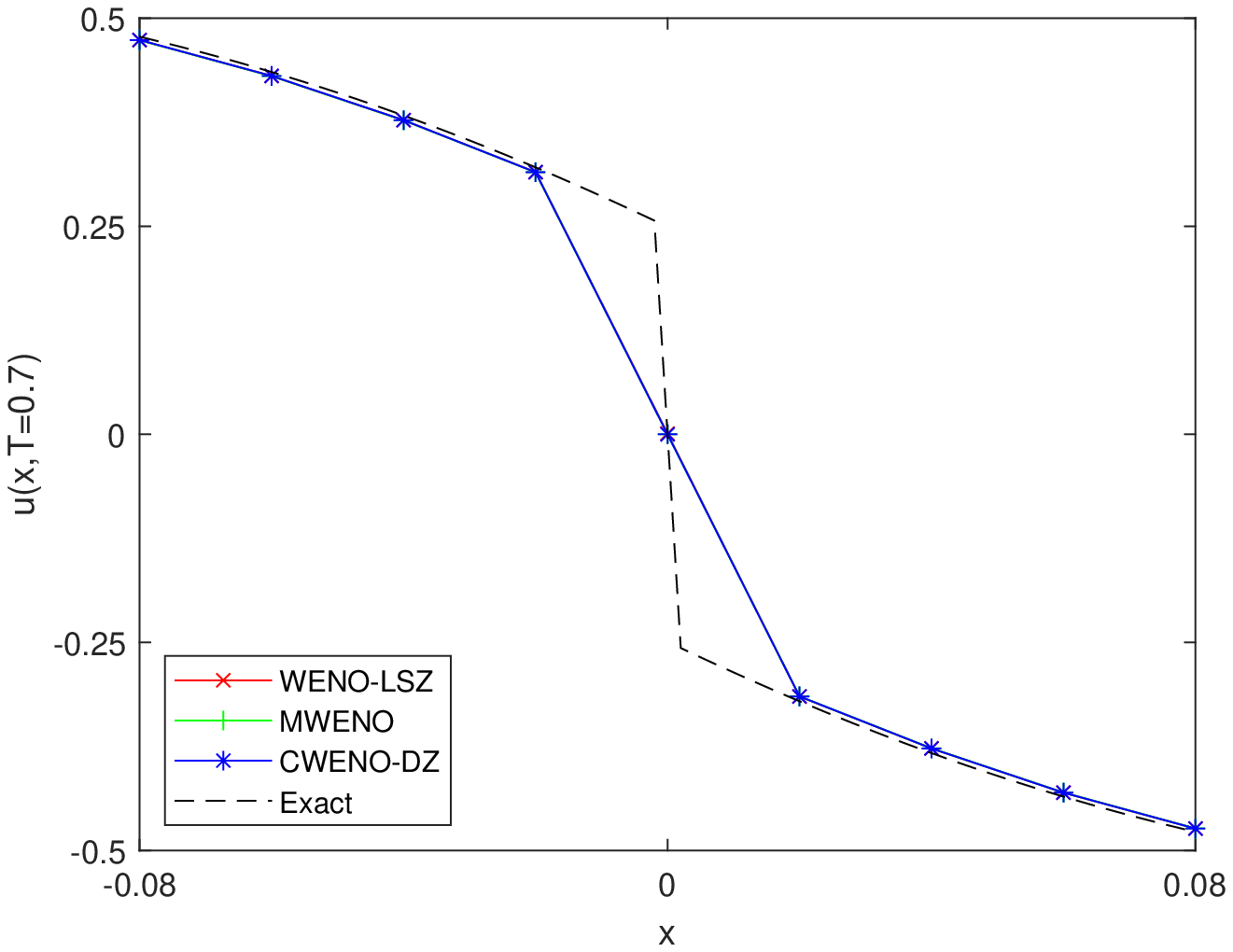}
\includegraphics[width=0.45\textwidth]{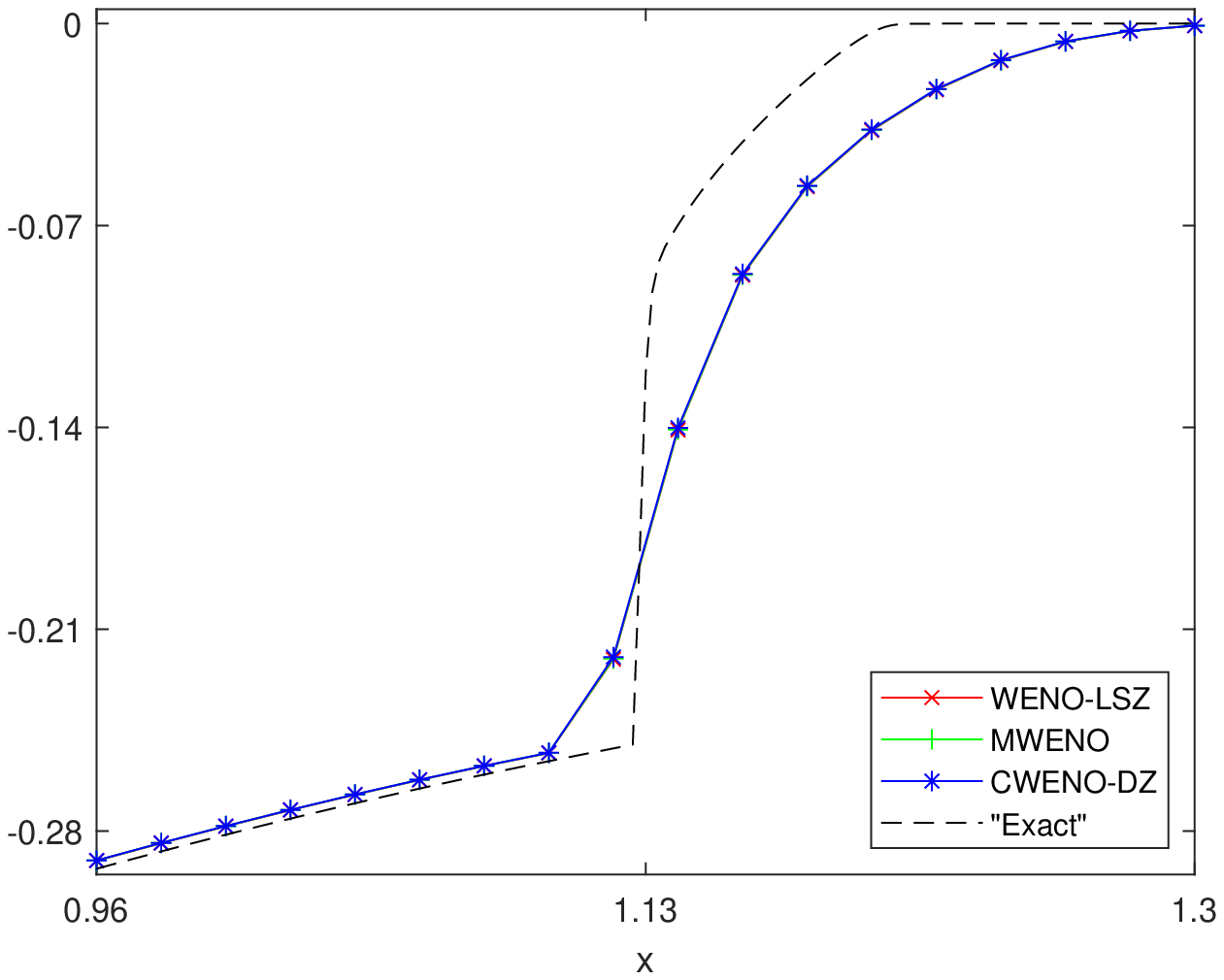}
\captionof{figure}{Solution profiles for Example \ref{ex:strongly_degenerate_cd_1d} at $T=0.7$ (top left), close-up view of the solutions in the boxes on the left/middle/right (top right/bottom left/bottom right) computed by WENO-LSZ (red), MWENO (green) and CWENO-DZ (blue) with $N = 200$. The dashed black lines are generated by WENO-M and WENO-LSZ with $N = 2000$.}
\label{fig:SDP}
\end{figure}

\begin{example} \label{ex:strongly_degenerate_cd_1d}
In this example, we consider the strongly degenerate parabolic convection-diffusion equation
\begin{equation} \label{eq:strongly_degenerate_cd_1d}
   u_t + f(u)_x = \epsilon \left( \nu(u) u_x \right)_x,~~\epsilon \nu(u) \geqslant 0.
\end{equation}
We take $\epsilon=0.1,~f(u)=u^2$, and
\begin{equation} \label{eq:diffusion_coeff_sdp}
 \nu(u) = \begin{cases}
           0, & |u| \leqslant 0.25, \\
           1, & |u| > 0.25.
          \end{cases}
\end{equation}
If $|u| \leqslant 0.25$, the equation \eqref{eq:strongly_degenerate_cd_1d} returns to the hyperbolic equation.
The diffusion term $\epsilon \left( \nu(u) u_x \right)_x$ can be written in the form of $g(u)_{xx}$, where
$$
   g(u) = \begin{cases}
          \epsilon (u+0.25), & u < -0.25, \\
          \epsilon (u-0.25), & u > 0.25, \\
          0,                 & \mbox{otherwise}. 
         \end{cases}
$$
The initial condition is given by 
$$
   u(x,0) = \begin{cases}
             1, & -1/\sqrt{2}-0.4 < x < -1/\sqrt{2}+0.4, \\ 
            -1, &  1/\sqrt{2}-0.4 < x < 1/\sqrt{2}+0.4, \\ 
             0, & \mbox{otherwise}.
            \end{cases}
$$
We divide the computational domain $[-2,~2]$ into $N=200$ uniform cells and the time step is $\Delta t = \cfl \cdot \Delta x^2$.
The simulations at $T=0.7$ are presented in Figure \ref{fig:SDP}, where the numerical result with CWENO-DZ is comparable to those with WENO-LSZ and MWENO.
\end{example}

\subsection{Two-dimensional numerical examples}

\begin{example} \label{ex:heat_2d}
We test the accuracy of those WENO schemes for the two-dimensional heat equation
$$
   u_t = u_{xx} + u_{yy},~~-\pi \leqslant x, y \leqslant \pi,~~t>0
$$
subject to the initial data 
$$
   u(x,y,0) = \sin(x+y),
$$
and the periodic boundary conditions in both directions.
The exact solution is 
$$
   u(x,t)= e^{-2t} \sin(x+y).
$$
The numerical solutions are computed at the final time $T=2$ with the time step $\Delta t = 0.2 \cdot \min (\Delta x, \Delta y)^2$.
The $L_1, L_2$, and $L_{\infty}$ errors, along with the orders of accuracy, are provided in Tables \ref{tab:heat_2d_L1}, \ref{tab:heat_2d_L2} and \ref{tab:heat_2d_Linf}, respectively. 
All WENO schemes exhibit sixth order accuracy overall.
As in Examples \ref{ex:heat_1d}, the errors produced by CWENO-DZ are larger than WENO-LSZ for $N=10$, but we see that the proposed CWENO-DZ scheme performs the best in terms of accuracy subsequently.
\end{example}

\begin{table}[h!]
\centering
\caption{$L_1$ errors and order of convergence for Example \ref{ex:heat_2d}.}      
\begin{tabular}{clcrlcrlc} 
\hline 
N & \multicolumn{2}{l}{WENO-LSZ} & & \multicolumn{2}{l}{MWENO} & & \multicolumn{2}{l}{CWENO-DZ}  \\ 
    \cline{2-3}                      \cline{5-6}                     \cline{8-9}                    
  & Error & Order                & & Error & Order               & & Error & Order              \\
\hline 
10 $\times$ 10   & 1.83E-6  & --     & & 9.20E-6   & --     & & 1.20E-5   & --      \\  
20 $\times$ 20   & 3.97E-8  & 5.5240 & & 6.10E-8   & 7.2381 & & 3.18E-9   & 11.8756 \\  
40 $\times$ 40   & 6.30E-10 & 5.9785 & & 6.55E-10  & 6.5412 & & 5.40E-10  & 2.5608  \\
80 $\times$ 80   & 9.71E-12 & 6.0189 & & 9.73E-12  & 6.0719 & & 9.51E-12  & 5.8271  \\ 
160 $\times$ 160 & 1.55E-13 & 5.9665 & & 1.55E-13  & 5.9696 & & 1.55E-13  & 5.9402  \\   
\hline
\end{tabular}
\label{tab:heat_2d_L1}
\end{table}

\begin{table}[h!]
\centering
\caption{$L_2$ errors and order of convergence for Example \ref{ex:heat_2d}.}      
\begin{tabular}{clcrlcrlc} 
\hline
N & \multicolumn{2}{l}{WENO-LSZ} & & \multicolumn{2}{l}{MWENO} & & \multicolumn{2}{l}{CWENO-DZ}  \\ 
    \cline{2-3}                      \cline{5-6}                     \cline{8-9}                    
  & Error & Order                & & Error & Order               & & Error & Order              \\
\hline
10 $\times$ 10   & 2.09E-6  & --     & & 1.06E-5  & --     & & 1.37E-5   & --     \\  
20 $\times$ 20   & 4.44E-8  & 5.5563 & & 6.83E-8  & 7.2737 & & 4.16E-9   & 11.6802 \\  
40 $\times$ 40   & 7.01E-10 & 5.9869 & & 7.29E-10 & 6.5513 & & 6.04E-10  & 2.7848 \\
80 $\times$ 80   & 1.08E-11 & 6.0212 & & 1.08E-11 & 6.0740 & & 1.06E-11  & 5.8355 \\ 
160 $\times$ 160 & 1.72E-13 & 5.9672 & & 1.72E-13 & 5.9702 & & 1.72E-13  & 5.9419 \\   
\hline
\end{tabular}
\label{tab:heat_2d_L2}
\end{table}

\begin{table}[h!]
\centering
\caption{$L_\infty$ errors and order of convergence for Example \ref{ex:heat_2d}.}      
\begin{tabular}{clcrlcrlc} 
\hline 
N & \multicolumn{2}{l}{WENO-LSZ} & & \multicolumn{2}{l}{MWENO} & & \multicolumn{2}{l}{CWENO-DZ} \\ 
    \cline{2-3}                      \cline{5-6}                     \cline{8-9}                    
  & Error & Order                & & Error & Order               & & Error & Order              \\
\hline 
10 $\times$ 10   & 2.76E-6  & --     & & 1.41E-5   & --     & & 1.78E-5  & --     \\  
20 $\times$ 20   & 6.26E-8  & 5.4641 & & 9.61E-8   & 7.1989 & & 7.46E-9  & 11.2224 \\  
40 $\times$ 40   & 9.91E-10 & 5.9814 & & 1.03E-9   & 6.5445 & & 8.61E-10 & 3.1154 \\
80 $\times$ 80   & 1.53E-11 & 6.0212 & & 1.53E-11  & 6.0728 & & 1.50E-11 & 5.8460 \\ 
160 $\times$ 160 & 2.44E-13 & 5.9673 & & 2.44E-13  & 5.9701 & & 2.43E-13 & 5.9433 \\   
\hline
\end{tabular}
\label{tab:heat_2d_Linf}
\end{table}

\begin{figure}[h!]
\centering
\includegraphics[width=0.32\textwidth]{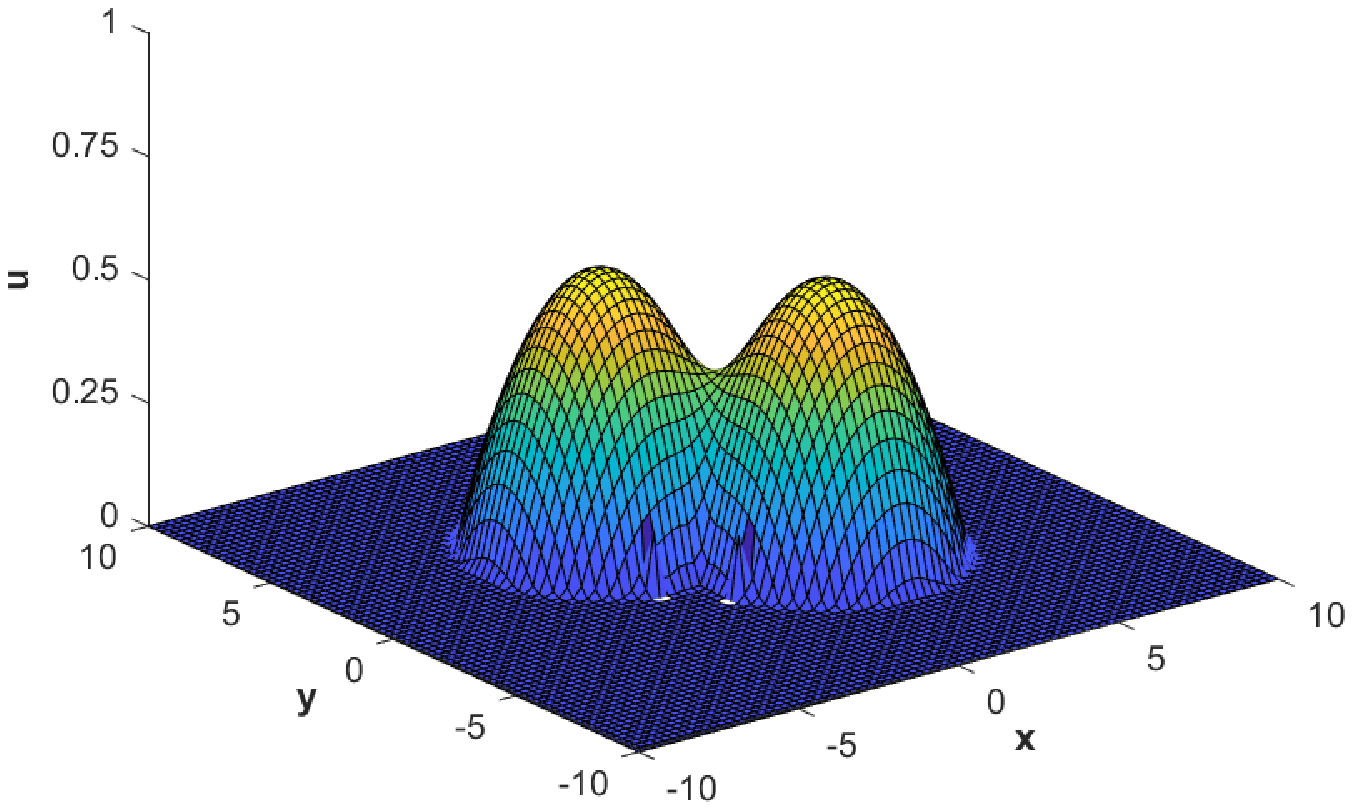}
\includegraphics[width=0.32\textwidth]{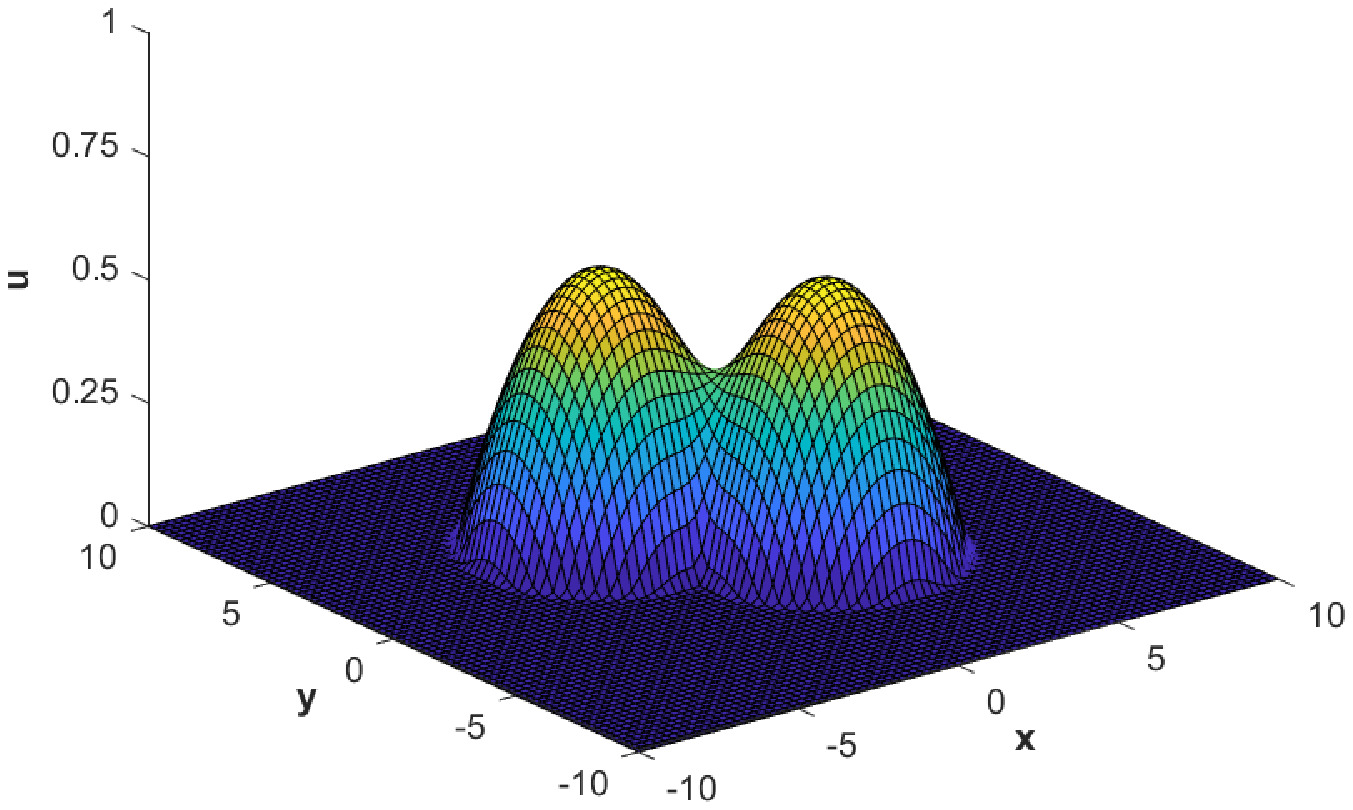}
\includegraphics[width=0.32\textwidth]{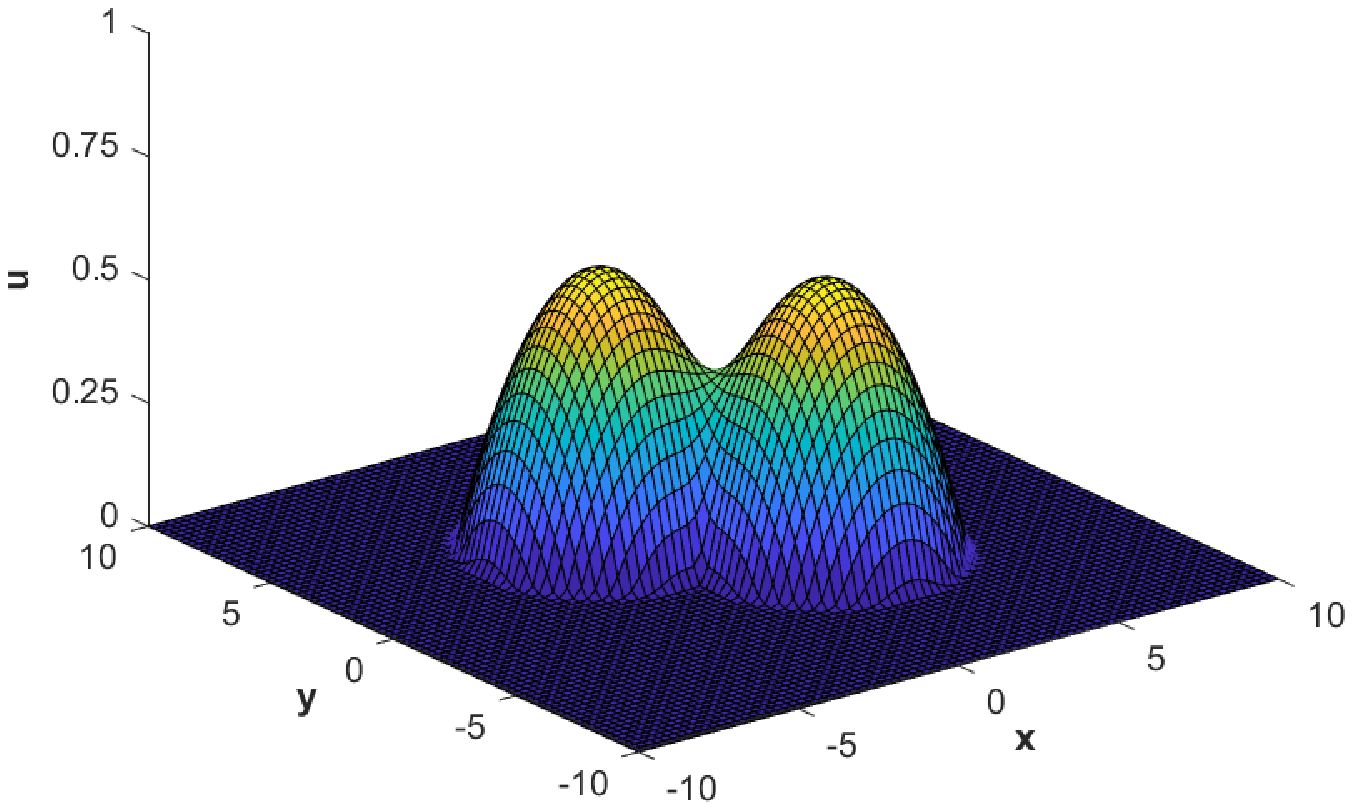}
\includegraphics[width=0.32\textwidth]{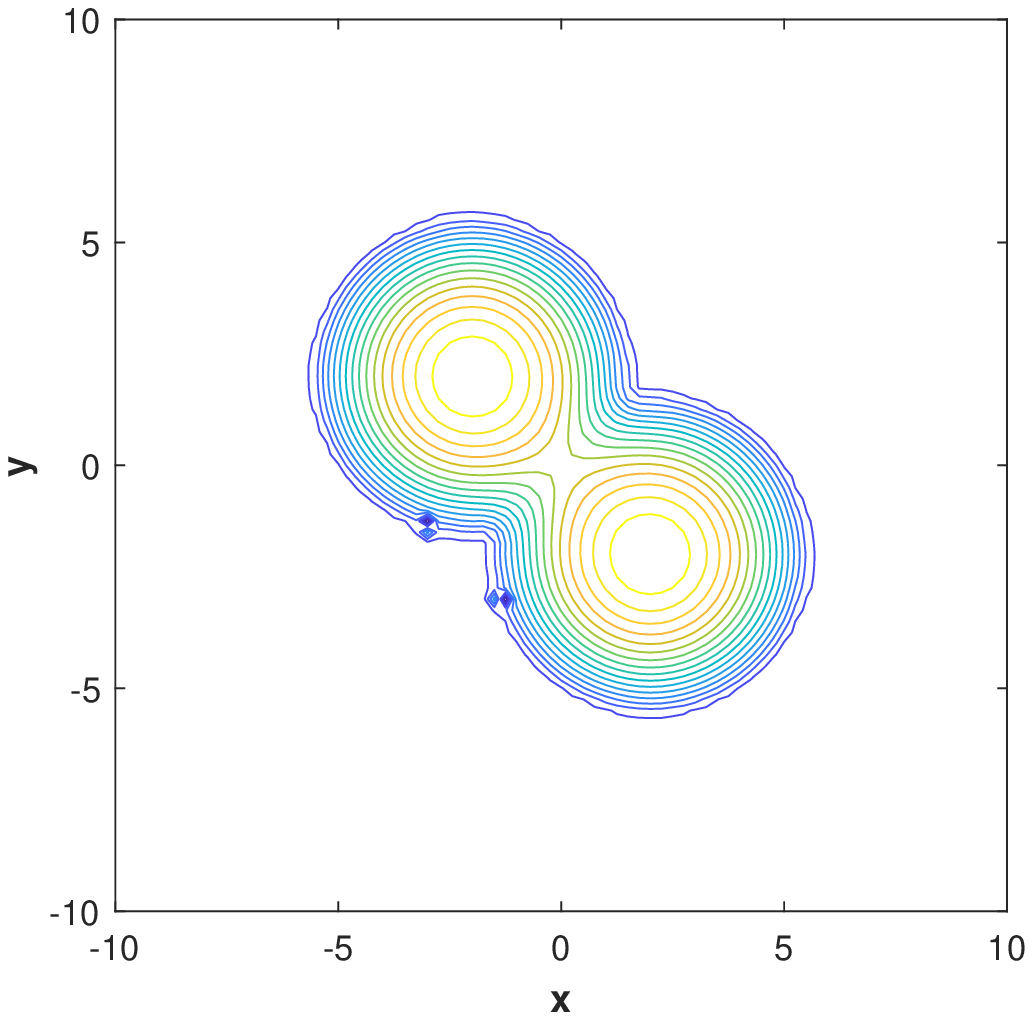}
\includegraphics[width=0.32\textwidth]{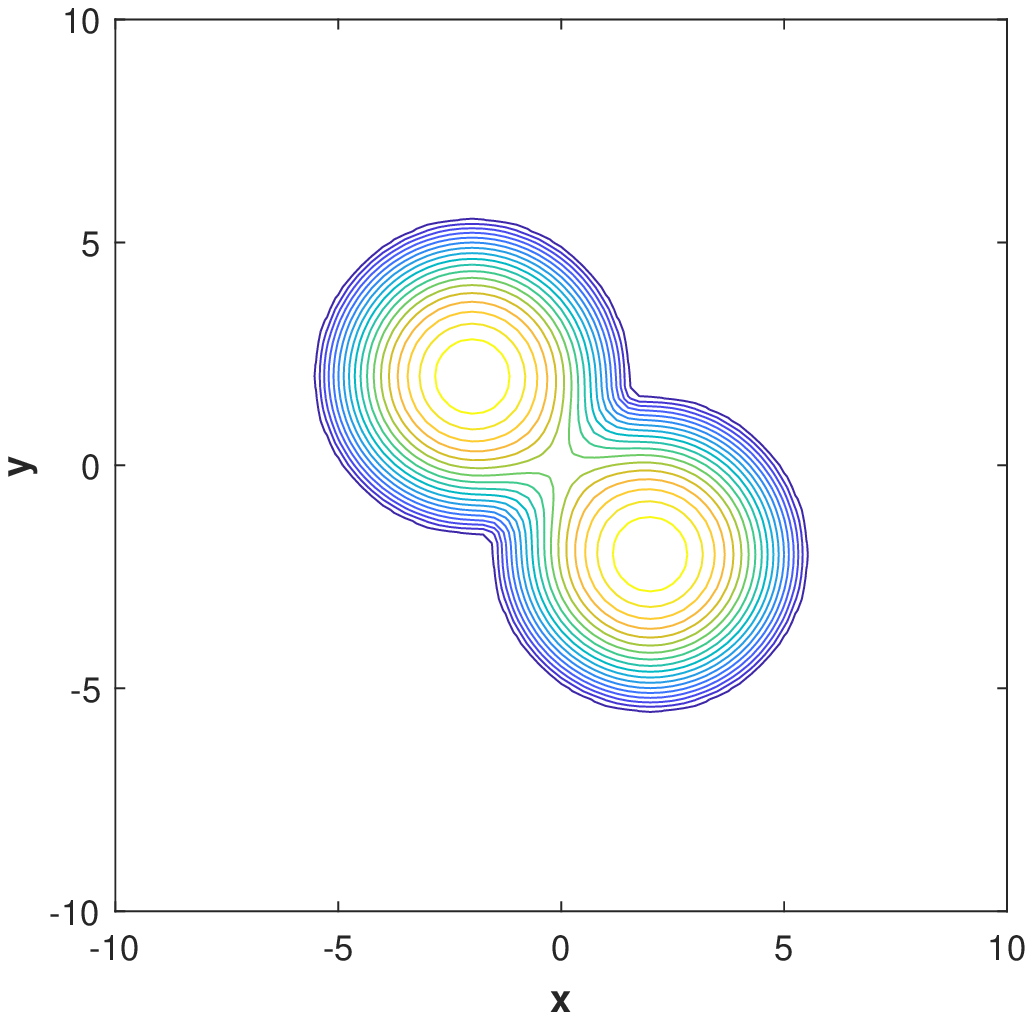}
\includegraphics[width=0.32\textwidth]{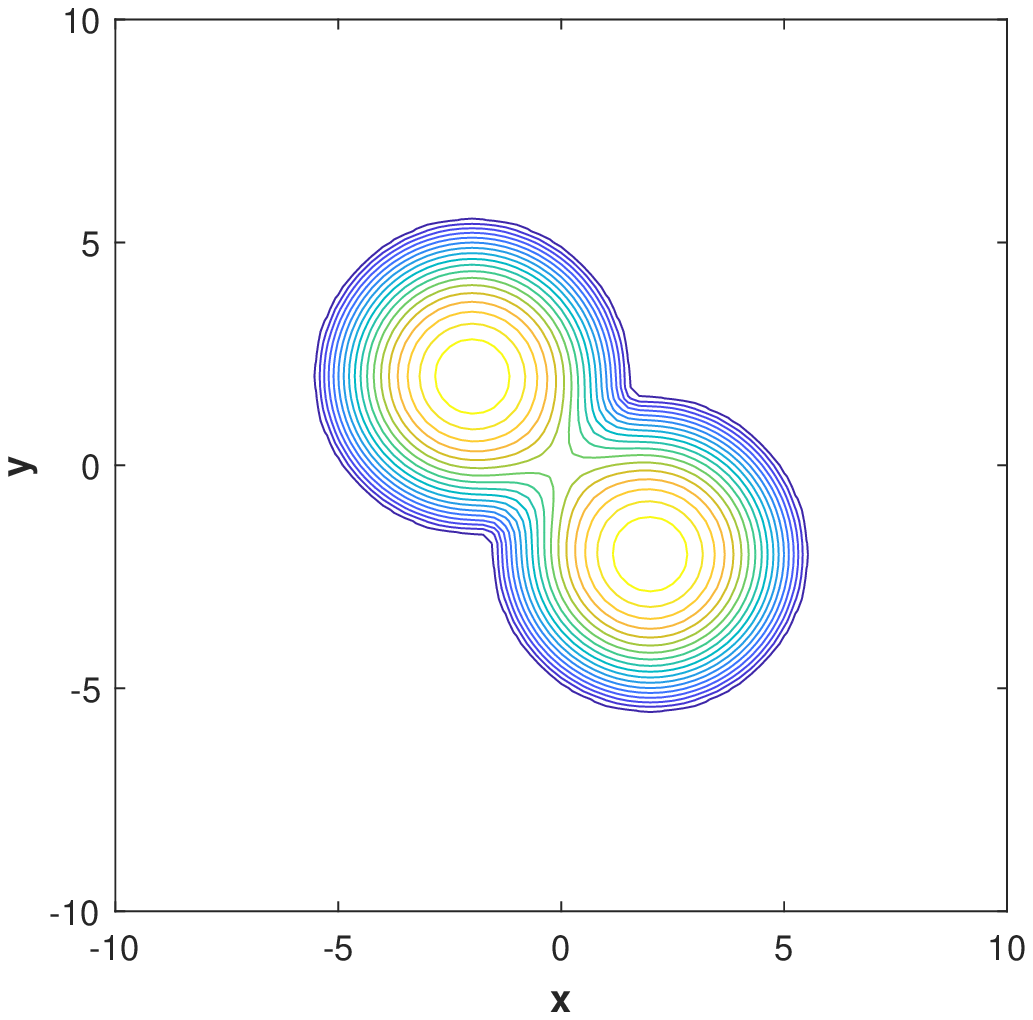}
\captionof{figure}{Solutions in the surface (top) and contour (bottom) plots for Example \ref{ex:PME_2d} at $t = 1$ by WENO-LSZ (left), MWENO (middle) and CWENO-DZ (right) with $N_x \times N_y = 80 \times 80$. Each contour plot includes 18 contours of $u$.}
\label{fig:PME_2d_T1}
\end{figure}

\begin{figure}[h!]
\centering
\includegraphics[width=0.32\textwidth]{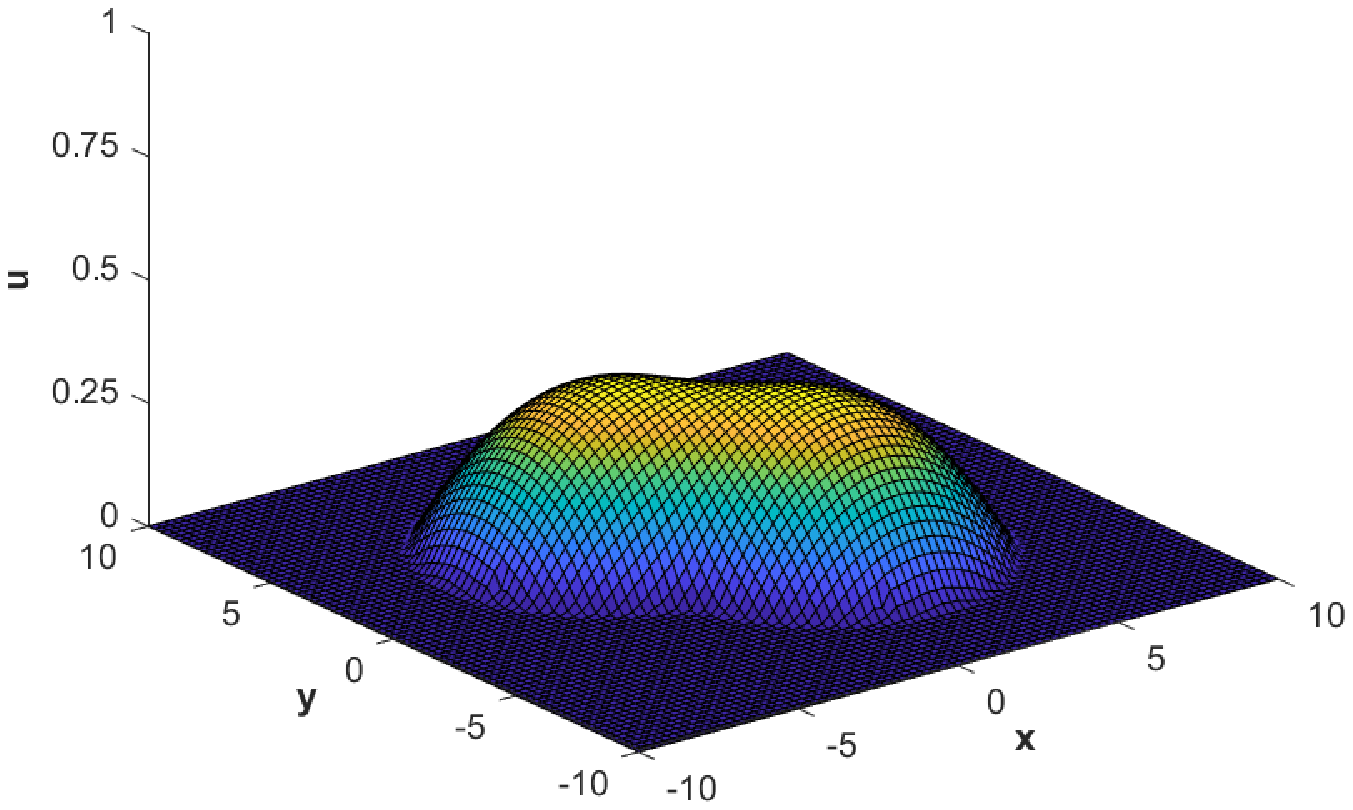}
\includegraphics[width=0.32\textwidth]{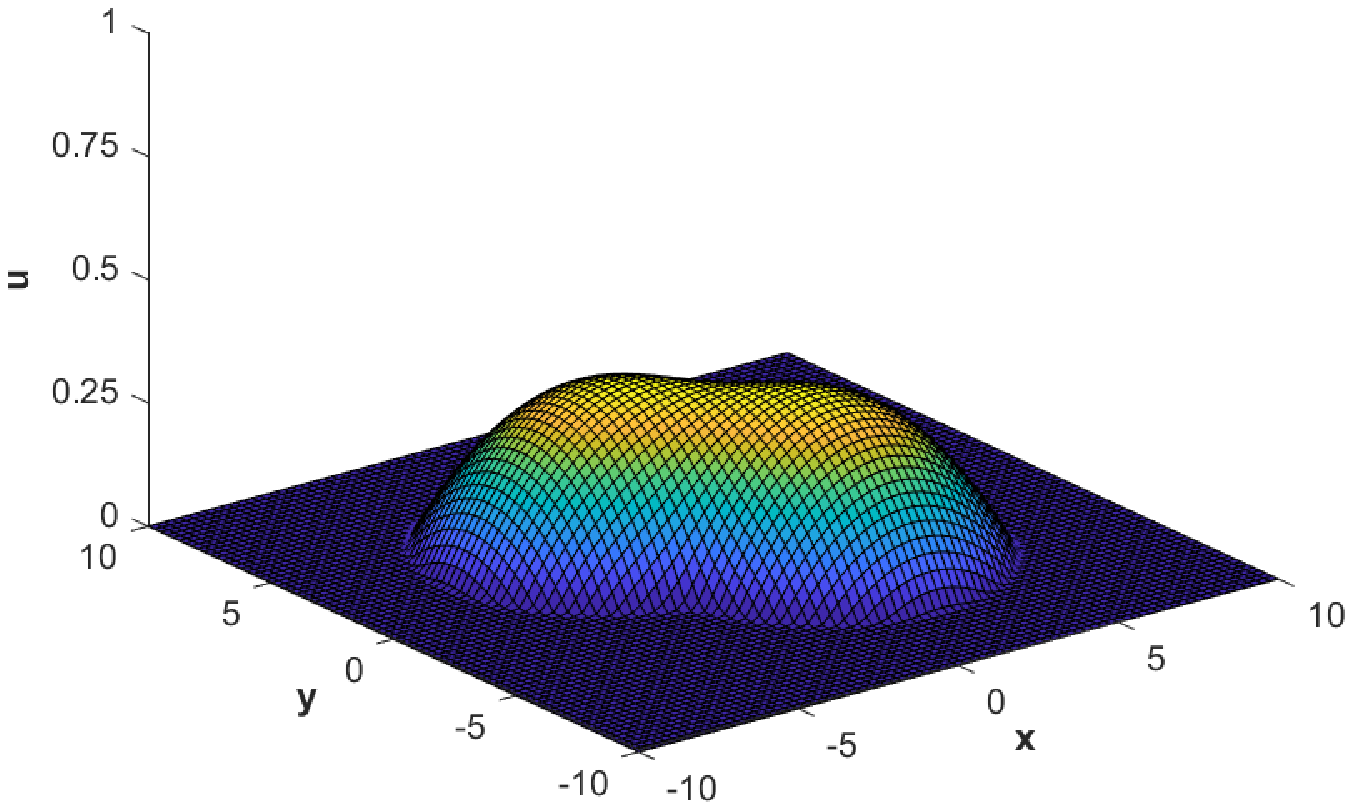}
\includegraphics[width=0.32\textwidth]{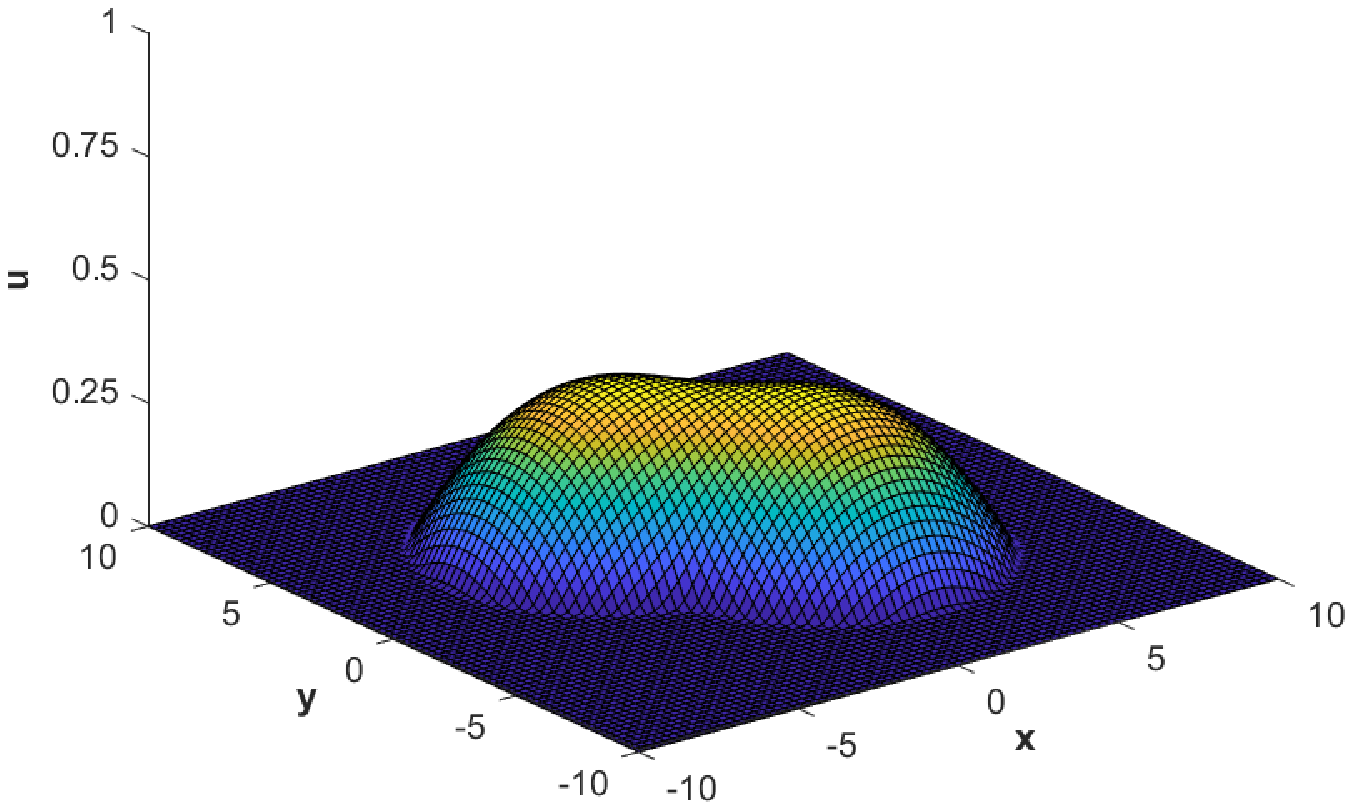}
\includegraphics[width=0.32\textwidth]{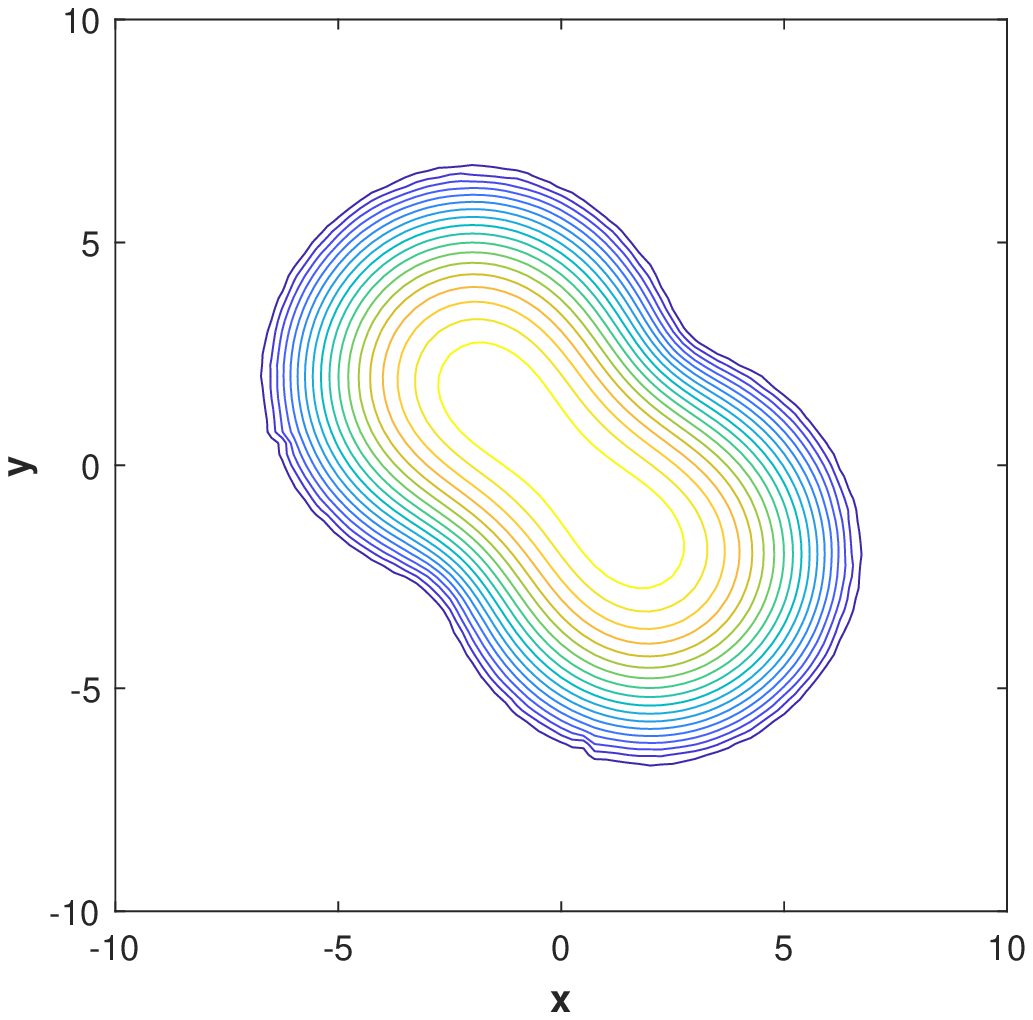}
\includegraphics[width=0.32\textwidth]{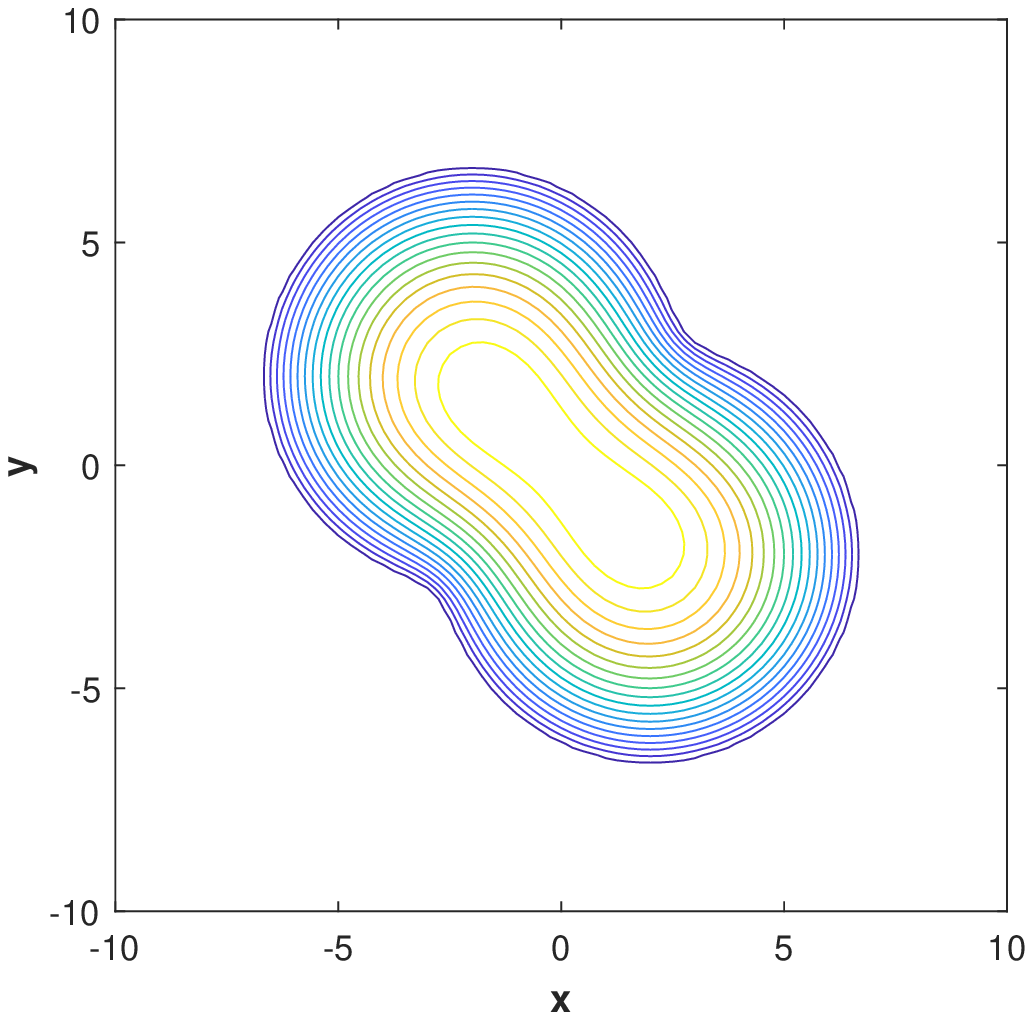}
\includegraphics[width=0.32\textwidth]{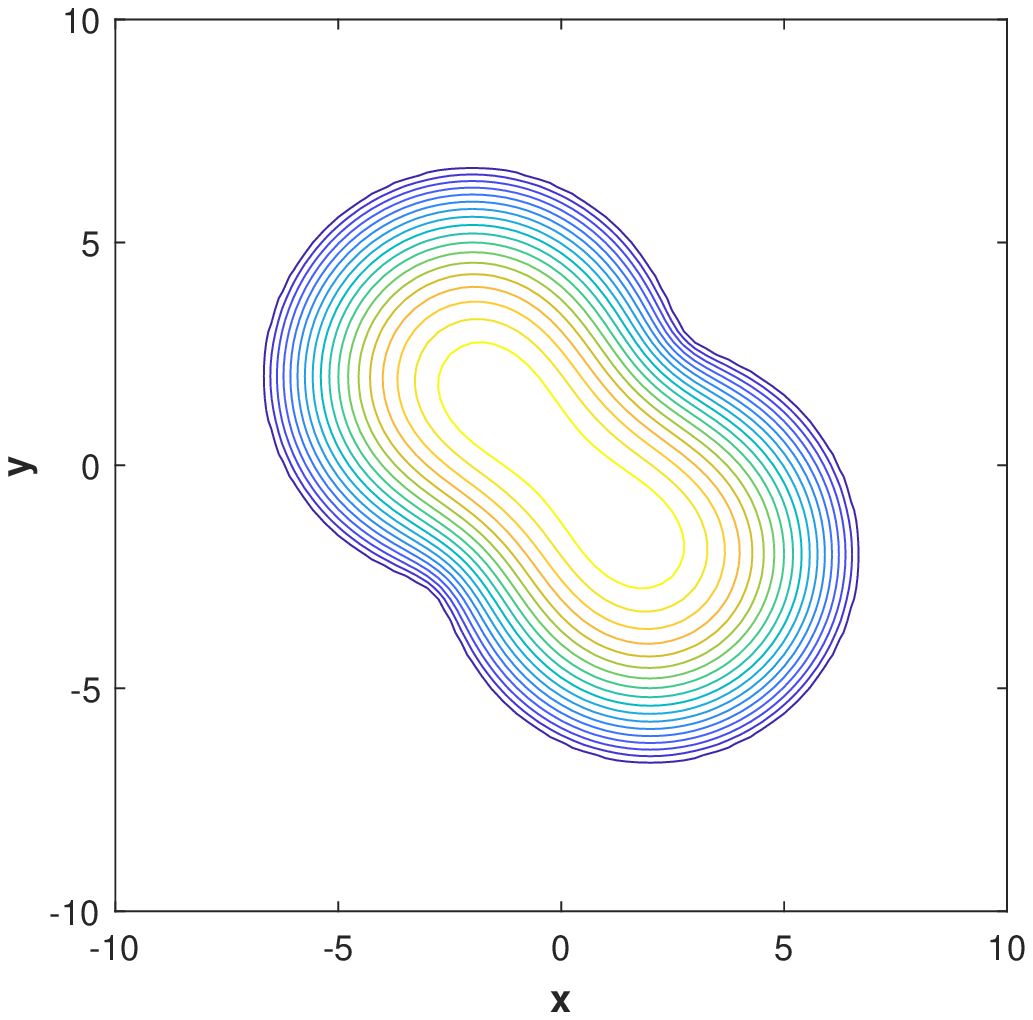}
\captionof{figure}{Solutions in the surface (top) and contour (bottom) plots for Example \ref{ex:PME_2d} at $t = 4$ by WENO-LSZ (left), MWENO (middle) and CWENO-DZ (right) with $N_x \times N_y = 80 \times 80$. Each contour plot includes 18 contours of $u$.}
\label{fig:PME_2d_T4}
\end{figure}

\begin{example} \label{ex:PME_2d}
Consider the two-dimensional PME given by
$$
   u_t = \left( u^2 \right)_{xx} + \left( u^2 \right)_{yy}, 
$$
with the initial condition 
$$
   u(x,y,0) = \begin{cases}
              \exp \left( - \frac{1}{6 - (x-2)^2 - (y+2)^2} \right), & (x-2)^2 + (y+2)^2 < 6, \\
              \exp \left( - \frac{1}{6 - (x+2)^2 - (y-2)^2} \right), & (x+2)^2 + (y-2)^2 < 6, \\
              0,                                                     & \mbox{otherwise},
              \end{cases}
$$
and the periodic boundary condition in each direction. 
We divide the square computational domain $[-10,~10] \times [-10,~10]$ into $N_x \times N_y = 80 \times 80$ uniform cells and the time step $\Delta t = \cfl \cdot \min (\Delta x, \Delta y)^4/2$.
The numerical solutions at $t=1$ and $t=4$ are shown in Figures \ref{fig:PME_2d_T1} and \ref{fig:PME_2d_T4}, respectively.
At the time $t=1$, there are some small-scale oscillations around the free boundaries in the solution by WENO-LSZ, which are implied by the white spots in the surface plot on the top left of Figure \ref{fig:PME_2d_T1}.
The oscillations are largely damped by MWENO and CWENO-DZ as there is no obvious white spot in the surface plot on the top middle and right, respectively.
However, at the time $t=4$, all WENO schemes are able to capture the free boundaries without noticeable oscillation, as shown in Figure \ref{fig:PME_2d_T4}.
Table \ref{tab:PME_2d} shows the minimum value of every numerical solution, which agrees with our observation above.
\end{example}

\begin{table}[h!]
\centering
\caption{Minimum values of the numerical solutions at $t=1$ and $t=4$ for Example \ref{ex:PME_2d}.}      
\begin{tabular}{cccc} 
\hline  
t & WENO-LSZ & MWENO & CWENO-DZ \\ 
\hline 
1 & -9.0127E-2 & -1.1547e-16 & -4.5836e-22 \\
\hline 
4 & -2.0504E-8 & -2.3381e-16 & -9.6261e-22 \\
\hline
\end{tabular}
\label{tab:PME_2d}
\end{table}

Finally, we use WENO schemes to solve the two-dimensional scalar convection-diffusion equations.
The WENO-JS scheme for the convection term is combined with WENO-LSZ for the diffusion term, while WENO-ZR \cite{Gu}, which gives sharper approximations around the shocks, is applied with both MWENO and CWENO-DZ.

\begin{figure}[h!]
\centering
\includegraphics[width=0.32\textwidth]{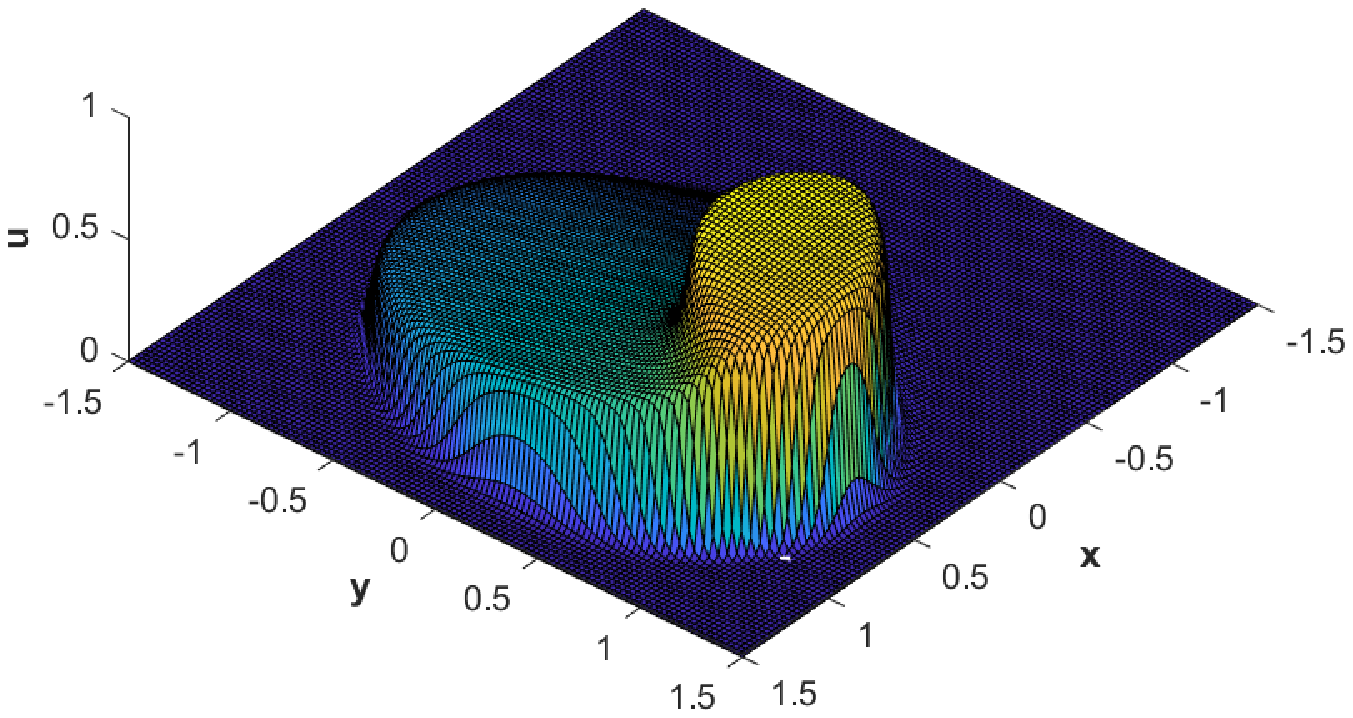}
\includegraphics[width=0.32\textwidth]{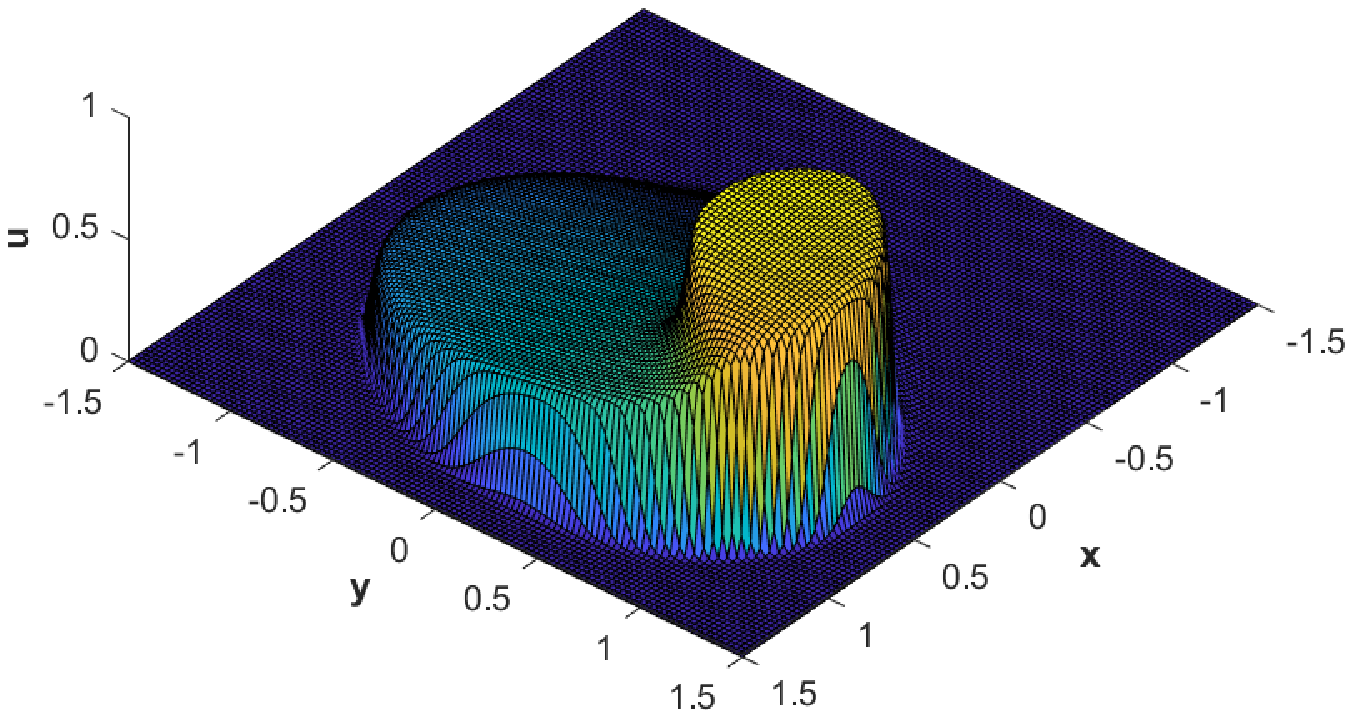}
\includegraphics[width=0.32\textwidth]{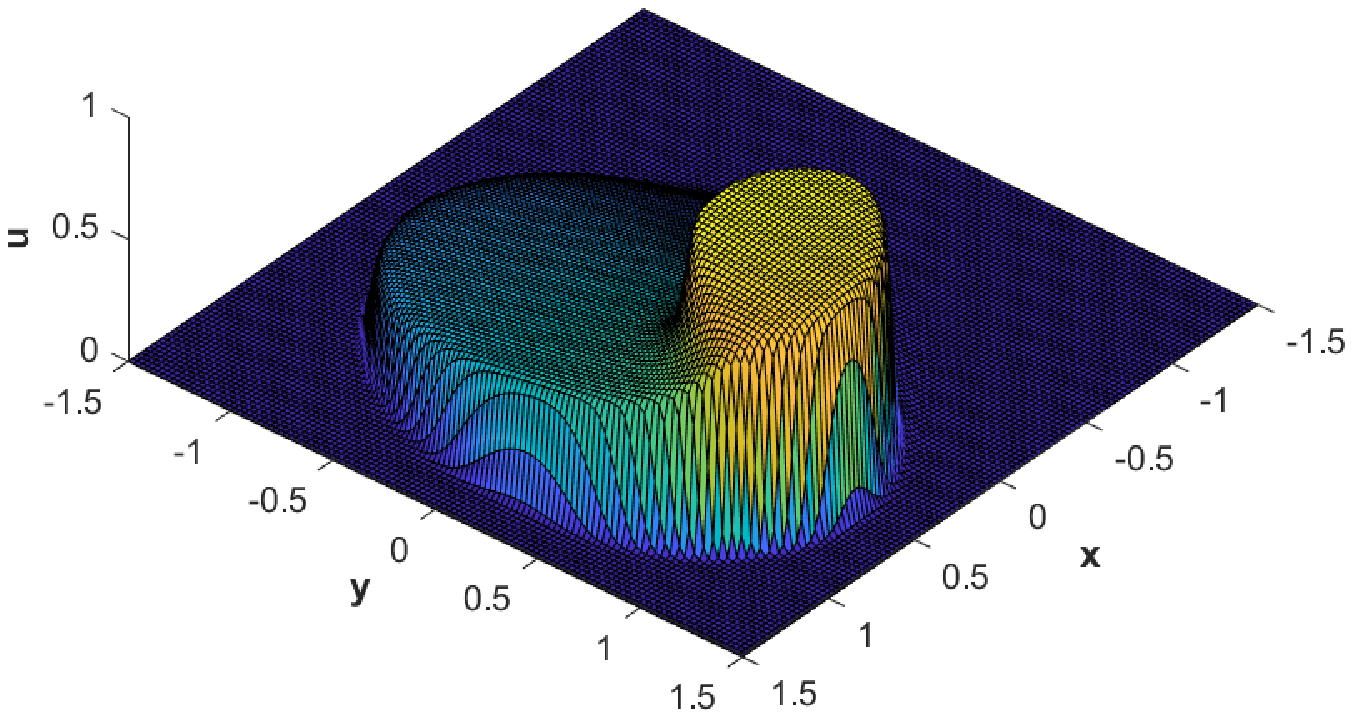}
\includegraphics[width=0.32\textwidth]{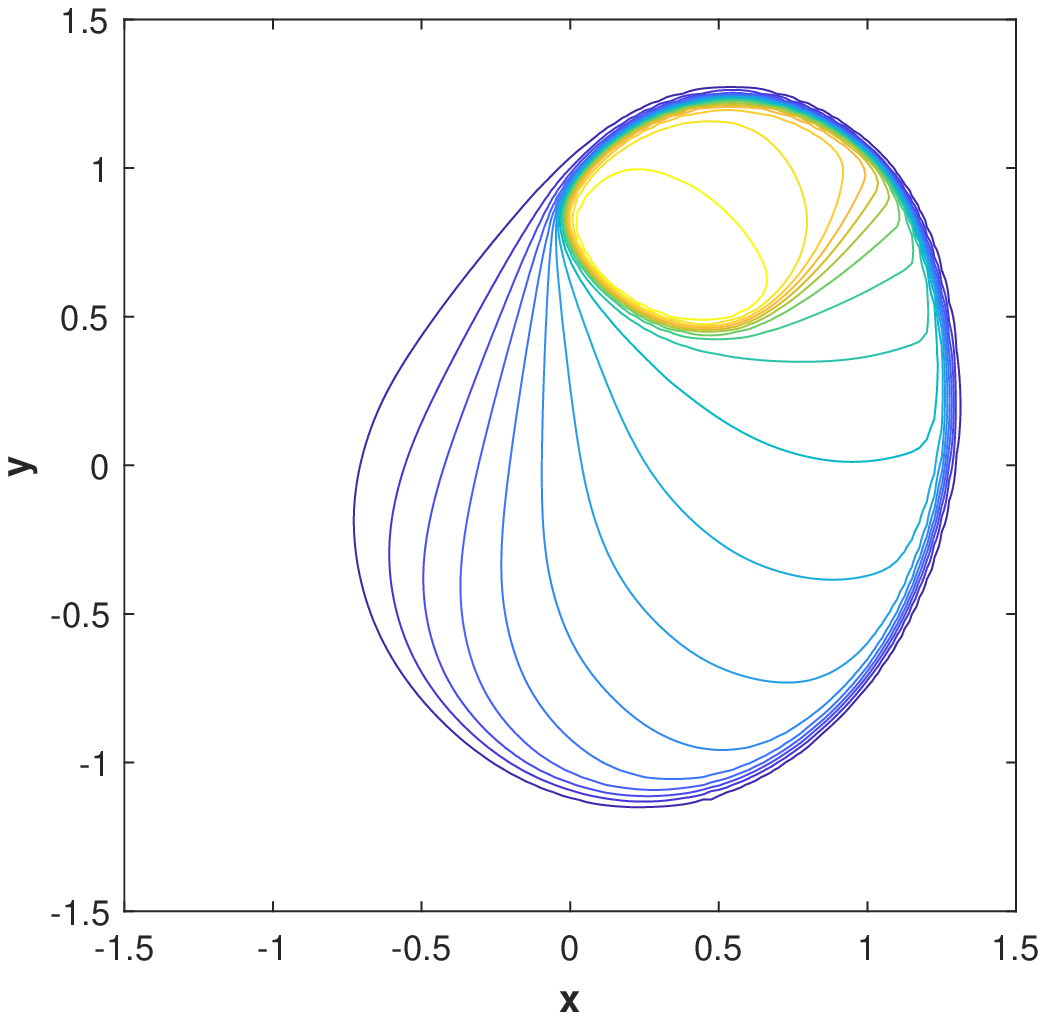}
\includegraphics[width=0.32\textwidth]{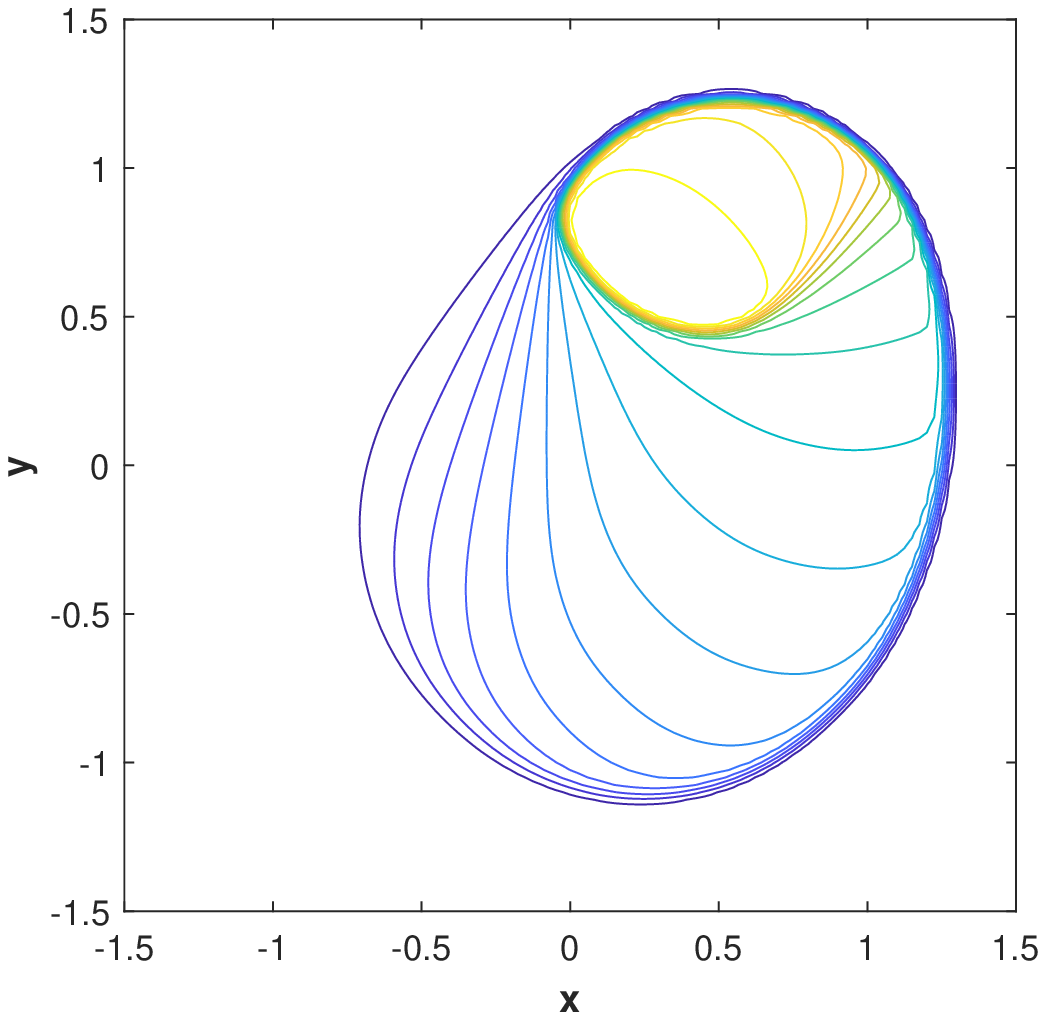}
\includegraphics[width=0.32\textwidth]{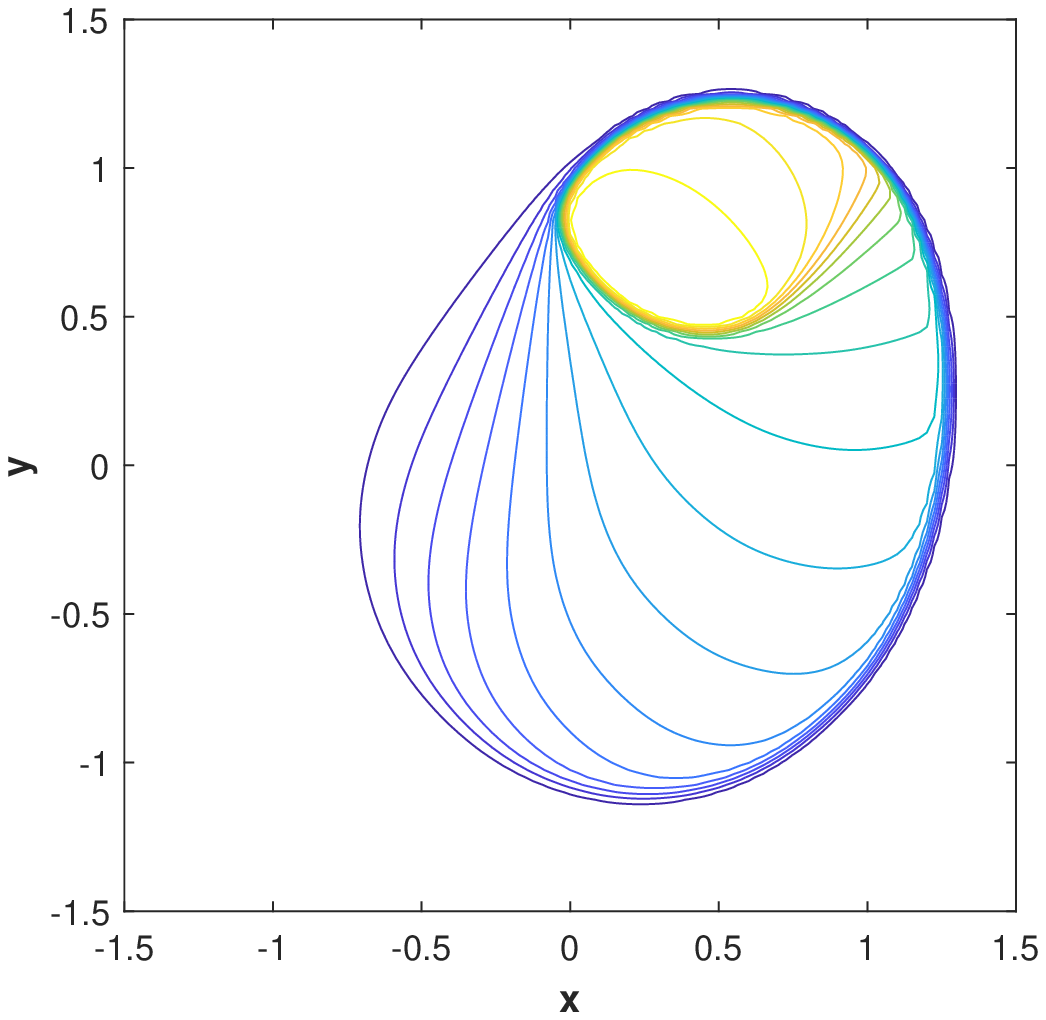}
\captionof{figure}{Solutions in the surface (top) and contour (bottom) plots for Example \ref{ex:Buckley_Leverett_2d} at $T = 0.5$ by WENO-JS/WENO-LSZ (left), WENO-ZR/MWENO (middle) and WENO-ZR/CWENO-DZ (right) with $N_x \times N_y = 120 \times 120$. Each contour plot includes 18 contours of $u$.}
\label{fig:BLE_2d}
\end{figure}

\begin{example} \label{ex:Buckley_Leverett_2d} 
We consider the two-dimensional Buckley-Leverett equation of the form
$$
   u_t + f_1(u)_x + f_2(u)_y = \epsilon \left( u_{xx} + u_{yy} \right),
$$
with $\epsilon = 0.01$ and the flux functions given by
$$
   f_1(u) = \frac{u^2}{u^2+(1-u)^2},~\quad f_2(u) = \left( 1 - 5(1-u)^2 \right) f_1(u).
$$
Then the equation includes gravitational effects only in the y-direction.
The initial condition is
$$
   u(x,y,0) = \begin{cases}
              1, & x^2 + y^2 < 0.5, \\
              0, & \mbox{otherwise}.
              \end{cases}
$$ 
The square computational domain $[-1.5,~1.5] \times [-1.5,~1.5]$ is divided into $N_x \times N_y = 120 \times 120$ uniform cells and the time step is $\Delta t = \cfl \cdot \min (\Delta x, \Delta y)^2$.
The solutions at $T = 0.5$ are plotted in Figure \ref{fig:BLE_2d}.
The white spot in the surface plot on the top left indicates the small-scale oscillations around the discontinuities in the solution by WENO-JS and WENO-LSZ.
Those oscillations are smoothed by WENO-ZR with both MWENO and CWENO-ZR, corresponding to the surface plot on the top middle and right, respectively.
We also provide Table \ref{tab:BLE_2d} showing the minimum value of each solution.
\end{example}

\begin{table}[h!]
\centering
\caption{Minimum values of the numerical solutions at $T = 0.5$ for Example \ref{ex:Buckley_Leverett_2d}.}      
\begin{tabular}{cccc} 
\hline  
T & WENO-LSZ & MWENO & CWENO-DZ \\ 
\hline 
0.5 & -6.2550E-3 & 1.5645E-39 & 8.0662E-39 \\
\hline 
\end{tabular}
\label{tab:BLE_2d}
\end{table}

\begin{figure}[h!]
\centering
\includegraphics[width=0.32\textwidth]{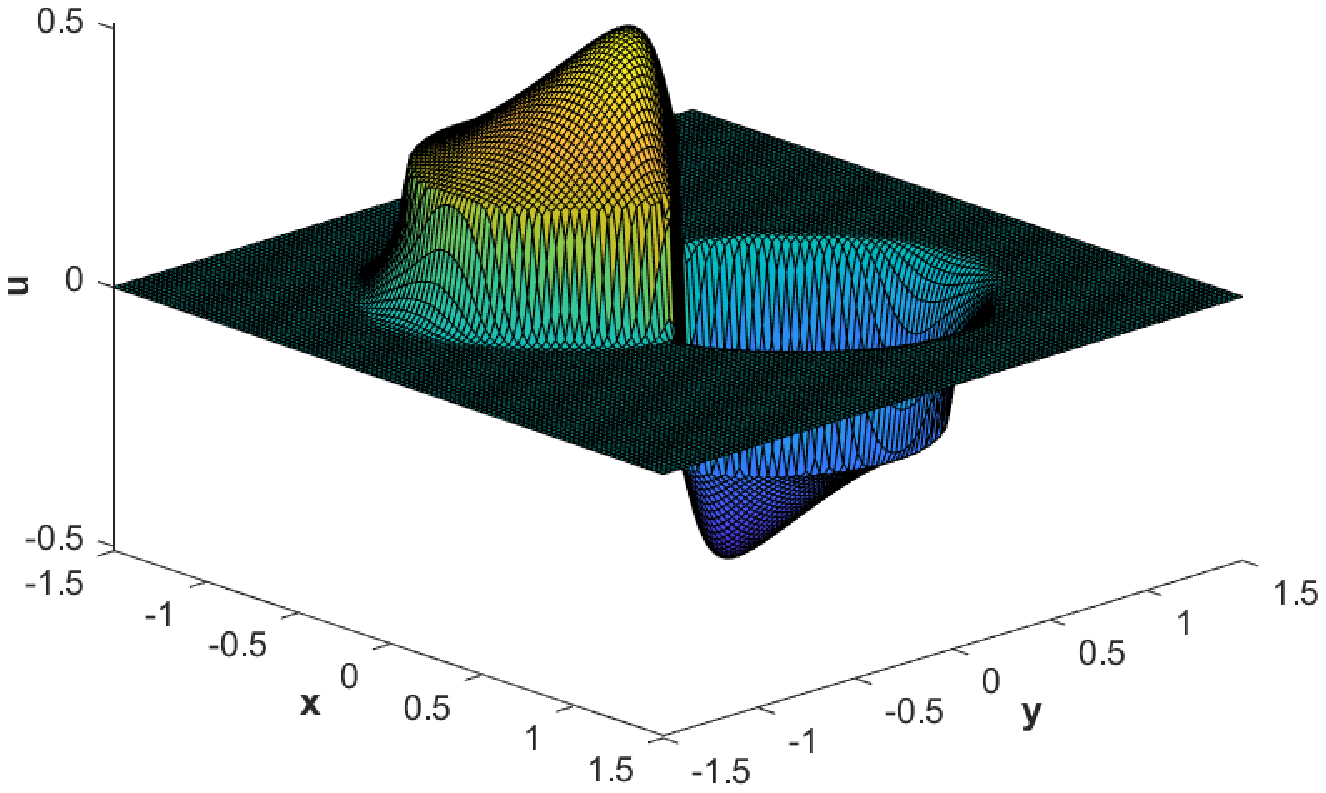}
\includegraphics[width=0.32\textwidth]{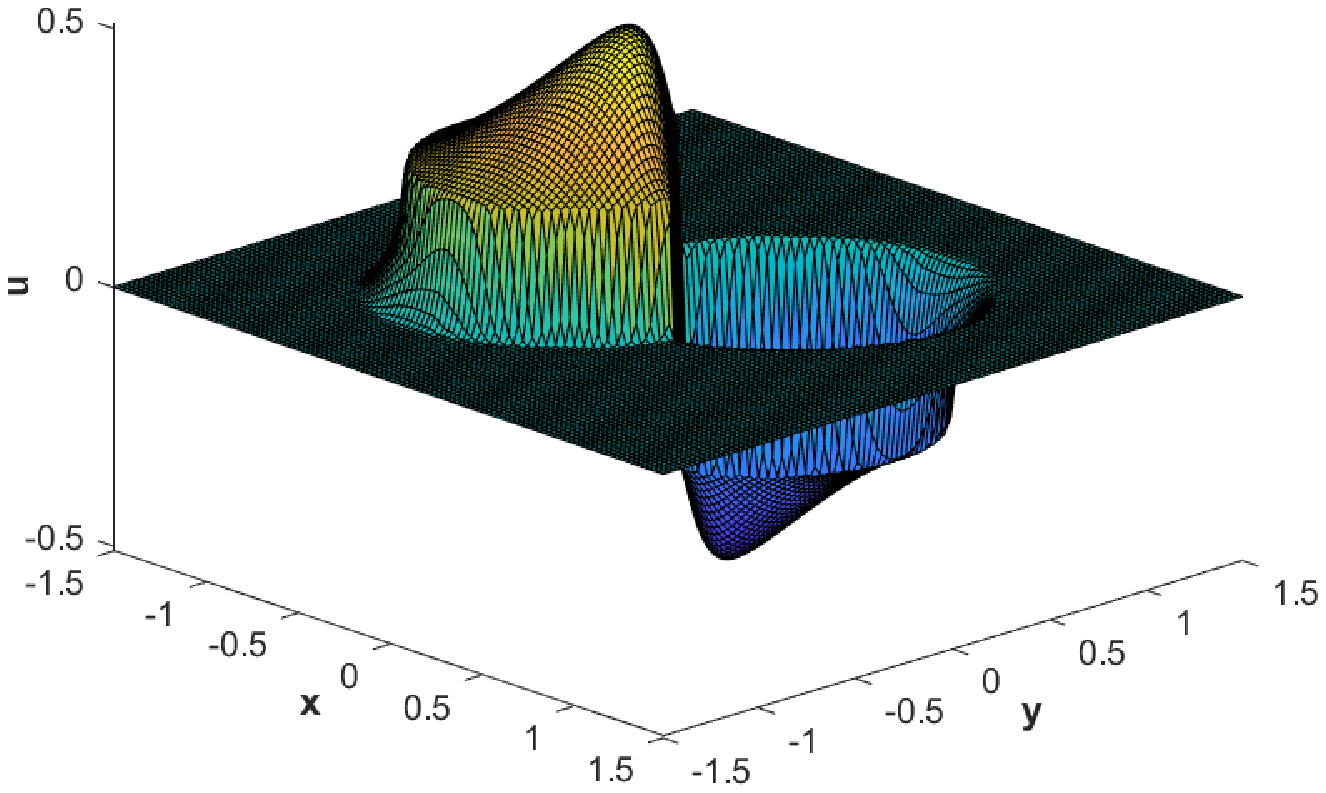}
\includegraphics[width=0.32\textwidth]{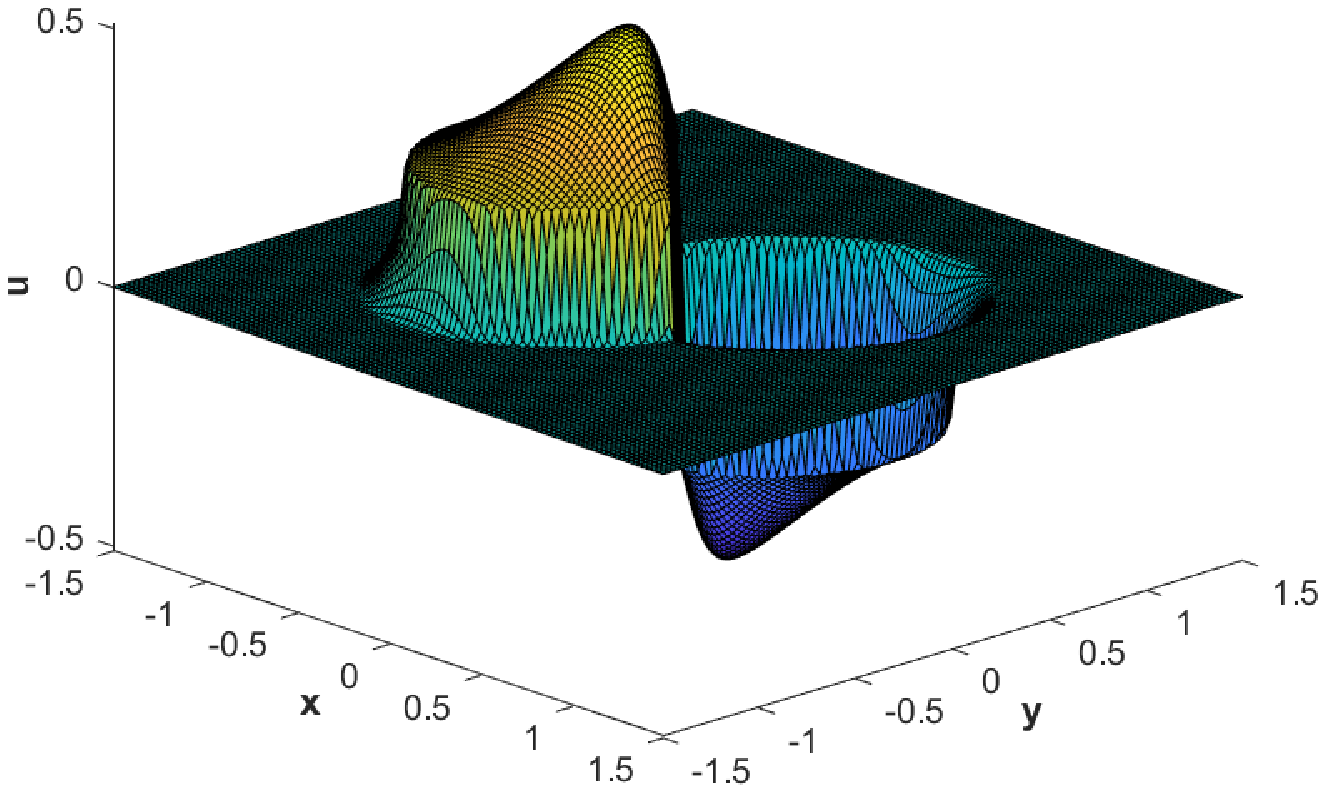}
\includegraphics[width=0.32\textwidth]{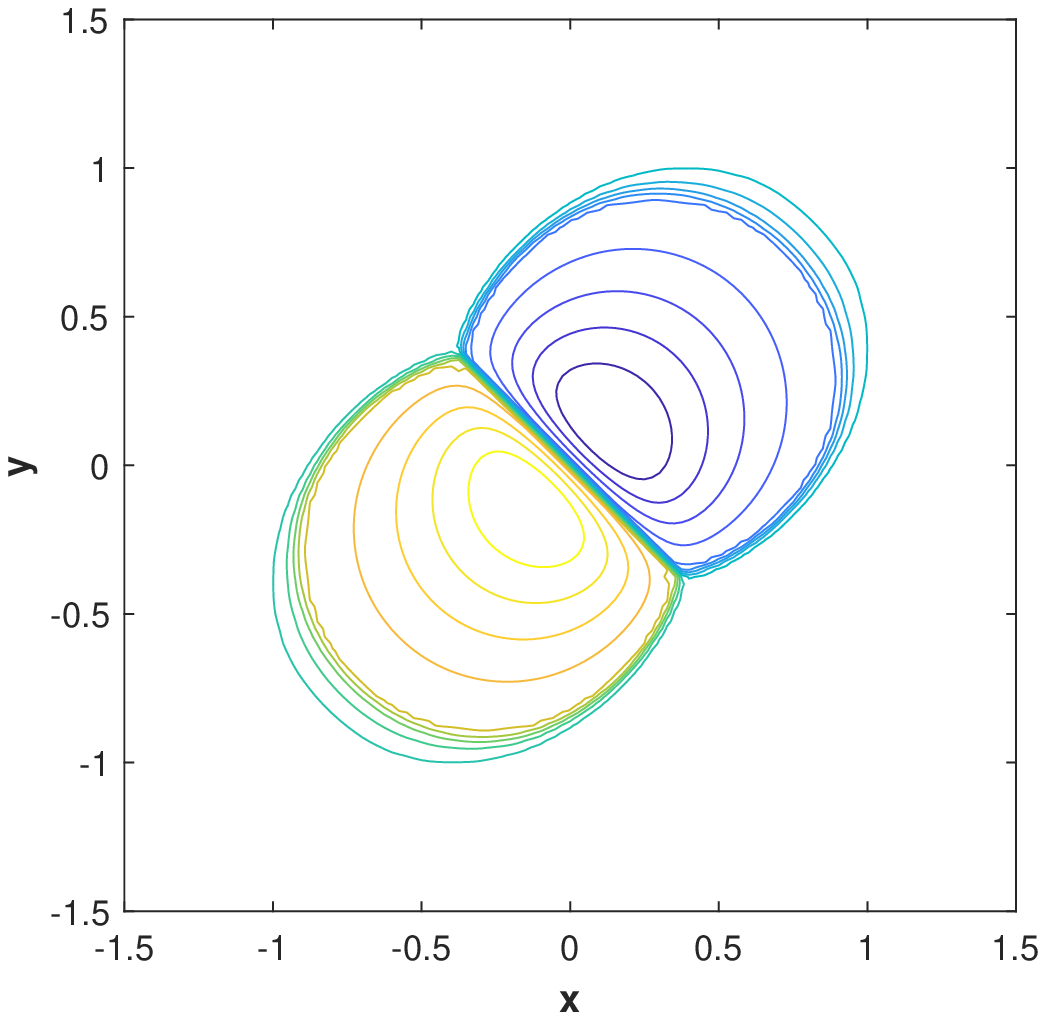}
\includegraphics[width=0.32\textwidth]{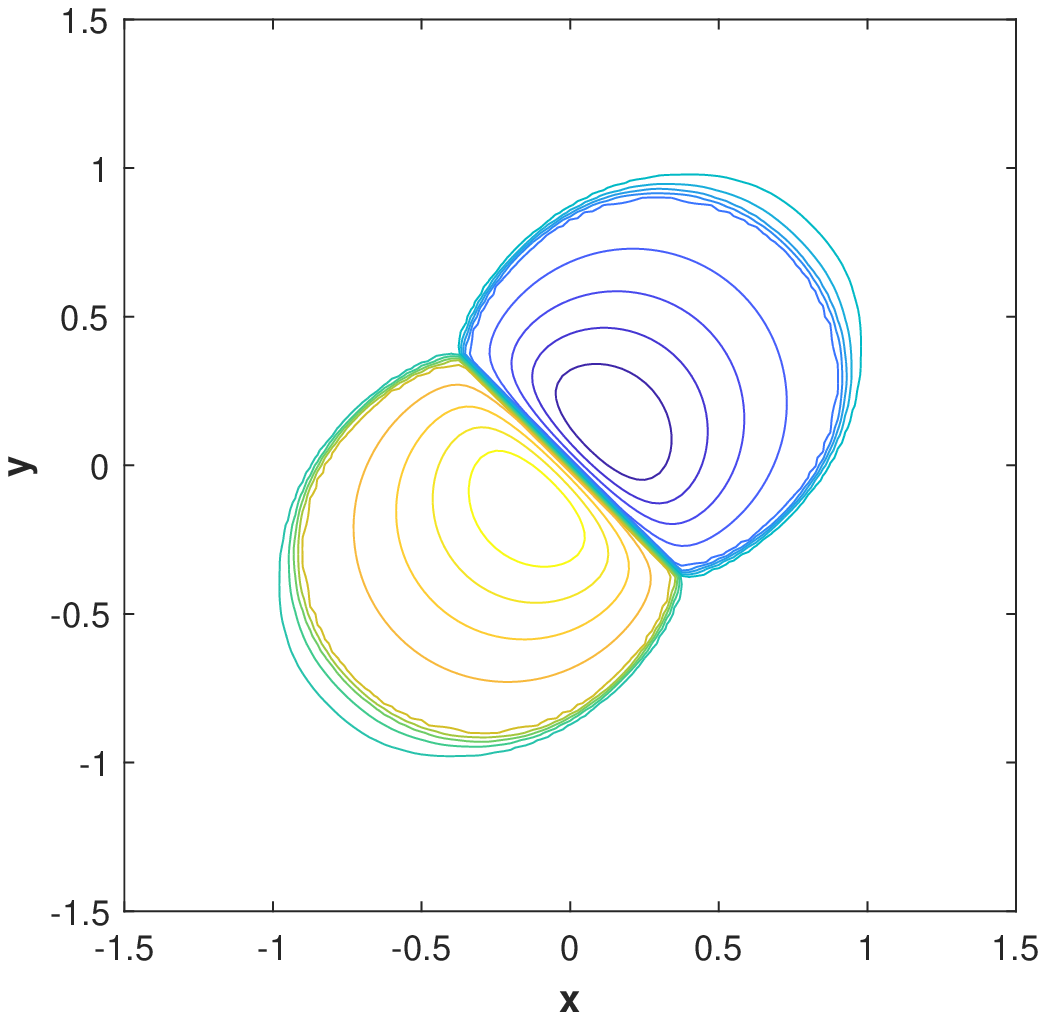}
\includegraphics[width=0.32\textwidth]{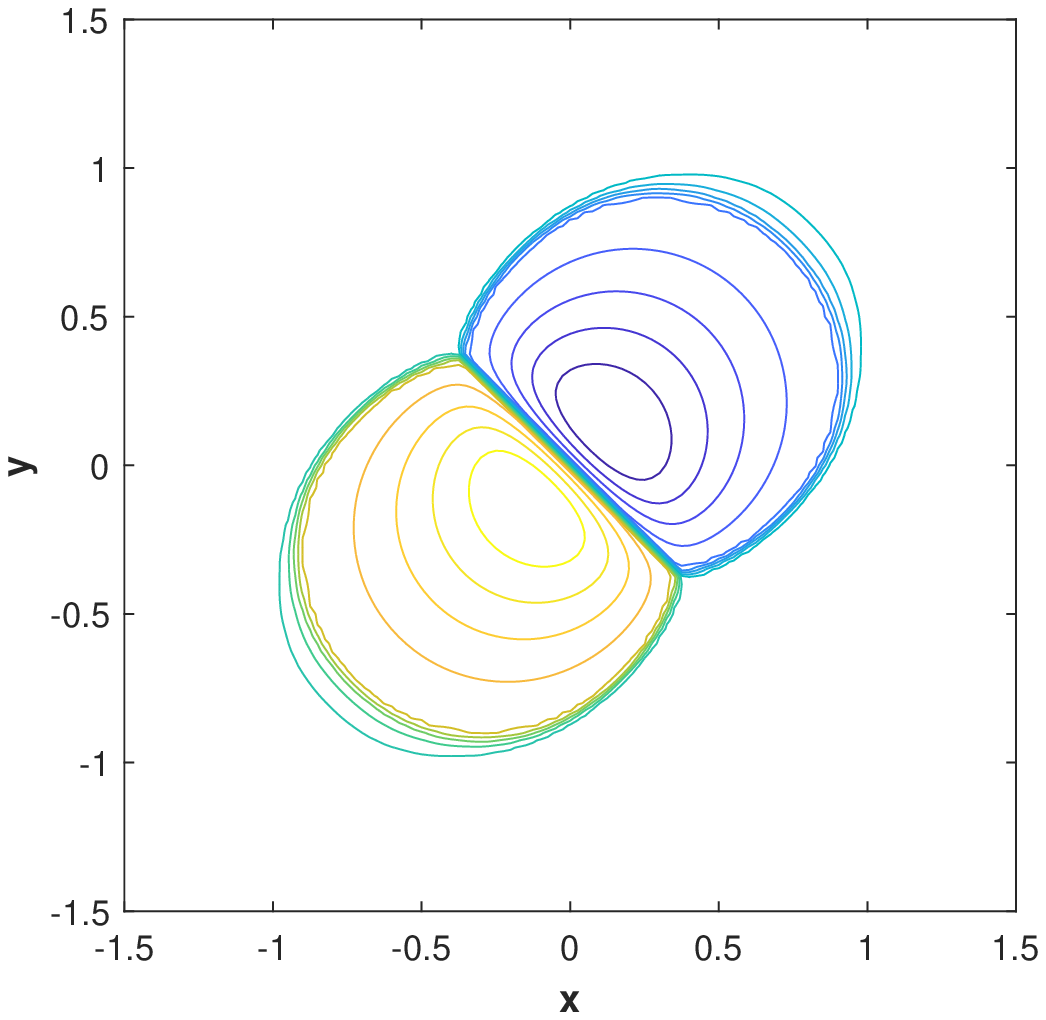}
\captionof{figure}{Solutions in the surface (top) and contour (bottom) plots for Example \ref{ex:strongly_degenerate_cd_2d} at $T = 0.5$ by WENO-JS/WENO-LSZ (left), WENO-ZR/MWENO (middle) and WENO-ZR/CWENO-DZ (right) with $N_x \times N_y = 120 \times 120$. Each contour plot includes 18 contours of $u$.}
\label{fig:SDP_2d}
\end{figure}

\begin{example} \label{ex:strongly_degenerate_cd_2d}
We conclude this section with the two-dimensional strongly degenerate parabolic convection-diffusion equation
$$
   u_t + f(u)_x + f(u)_y  = \epsilon \left( \nu(u) u_x \right)_x + \epsilon \left( \nu(u) u_y \right)_y,
$$
where $\epsilon=0.1,~f(u)=u^2$, and $\nu(u)$ \eqref{eq:diffusion_coeff_sdp} are the same as in Example \ref{ex:strongly_degenerate_cd_1d}.
The initial condition is
$$
   u(x,y,0) = \begin{cases}
               1, & (x+0.5)^2+(y+0.5)^2<0.16, \\
              -1, & (x-0.5)^2+(y-0.5)^2<0.16, \\
               0, & \mbox{otherwise}.
              \end{cases}
$$
We divide the computational domain $[-1.5,~1.5] \times [-1.5,~1.5]$ into $N_x \times N_y = 120 \times 120$ uniform cells and the time step is $\Delta t = \cfl \cdot \min (\Delta x, \Delta y)^2$.
The numerical solutions at $T = 0.5$, generated by those WENO schemes, look similar in Figure \ref{fig:SDP_2d}.
\end{example}

\section{Conclusion} \label{sec:conclusion}
In this paper, we proposed a six-order finite difference CWENO scheme to solve nonlinear degenerate parabolic equations.
The key idea is to introduce a centered polynomial such that the positivity of linear weights is guaranteed.
Numerical examples show that the proposed CWENO scheme achieves sixth order accuracy with smaller errors than WENO-LSZ and MWENO, and inhibits the small-scale oscillations introduced by WENO-LSZ.

\section*{Acknowledgments}
The first author is supported by IIPE, Visakhapatnam, India, under the IRG grant number $\text{IIPE/DORD/IRG/001}$ and NBHM, DAE, India (Ref. No. $\text{02011/46/2021 NBHM(R.P.)/R \& D II/14874}$). 
The second author is supported by POSTECH Basic Science Research Institute under the NRF grant number $\text{NRF2021R1A6A1A1004294412}$.


\end{document}